\documentclass[10pt]{amsart}

\usepackage{graphicx}
\usepackage{psfrag}

\newtheorem{theorem}{Theorem}[section]

\newtheorem{lemma}[theorem]{Lemma}           
\newtheorem{cor}[theorem]{Corollary}

\def\newtext{}
\newcommand{\HOX}[1]{}
\def\newtekstt{}
\def\newtekst{}

\def \ba {\begin {eqnarray*} }
\def \ea {\end {eqnarray*} }
\def \beq {\begin {eqnarray}}    
\def \eeq {\end {eqnarray}}

\def \det {\hbox{det}}

\def \e {\varepsilon}
\def \p {\partial}

\def\R{\mathbb R}

\def\C{\mathbb C}

\def\O{{\mathcal O}}

\theoremstyle{definition}
\newtheorem{definition}[theorem]{Definition}

\theoremstyle{remark}

\numberwithin{equation}{section}

\begin{document}
\baselineskip 14pt

\author{Hiroshi Isozaki}
\address{Hiroshi Isozaki, Institute of Mathematics \\
University of Tsukuba,
Tsukuba, 305-8571, Japan}
\author{Yaroslav Kurylev}
\address{Yaroslav Kurylev, Department of Mathematics \\
University College of London, United Kingdom}
\author{Matti Lassas}
\address{Matti Lassas, Department of Mathematics and Statistics\\
University of Helsinki, Finland}
\title[manifolds with cylindrical ends]{Forward and inverse scattering on manifolds with 
asymptotically cylindrical ends}
\date{May 11,  2009}

\maketitle
\begin{abstract}
We study an inverse problem for a non-compact Riemannian manifold whose ends have the following properties : On each end, the Riemannian metric is assumed to be a short-range perturbation of the metric of the form $(dy)^2 + h(x,dx)$, $h(x,dx)$ being the metric of some compact manifold of codimension 1. Moreover  one end  is exactly cylindrical, i.e. the metric is equal to $(dy)^2 + h(x,dx)$. Given two such manifolds having the same scattering matrix on that exactly cylindrical end for all energy, we show that these two manifolds are isometric.   
\end{abstract}


\section{Introduction}
The aim of this paper is to study spectral properties and related inverse problems for a connected, non-compact Riemannian manifold $\Omega$ of dimension $n \geq 2$ with or without boundary. We assume that $\Omega$ is split into $N+1$ parts
\begin{equation}
\Omega = \mathcal K \cup \Omega_1 \cup \cdots \cup \Omega_N,
\label{S1:Omega}
\end{equation}
where $\mathcal K$ is a bounded open set, and $\Omega_i$, called an {\it end} of $\Omega$, is diffeomorphic to $M_i\times(0,\infty)$, $M_i$ being a compact manifold of dimension $n-1$. (See the figure 1.) Denoting the local coordinates on $M_i$ by $x$, we assume that $M_i$ is equipped with a Riemannian metric $h_i(x,dx) = \sum_{p,q=1}^{n-1}h_{i,pq}(x)dx^pdx^q$. Letting $y$ be the coordinate on $(0,\infty)$, we denote the local coordinates on $\Omega_i$  by $X = (x,y)$.
We assume that the Riemannian metric $G$ on $\Omega$, which  is denoted by $G_i =$ $\sum_{p,q=1}^ng_{i,pq}(X)dX^pdX^q$ on $\Omega_i$,  has the following property
\begin{equation}
|\partial_X^{\alpha}(g_{i,pq}(X) - h_{i,pq}(x))| \leq C_{\alpha}(1 + y)^{-1-\epsilon_0}, \quad \forall \alpha,
\label{eq:Sec1DecayAssumSR}
\end{equation}
 where $h_{i,pn}(x) = h_{i,np}(x) = 0$ if $1 \leq p \leq n-1$ and $h_{i,nn}(x) = 1$, and $C_{\alpha}$ is a constant. The metric $h_i(x,dx)$ on $M_i$ is allowed to be different for different ends. We shall assume either $\Omega$ has  no boundary or each $M_i$, consequently $\Omega$ itself, has a boundary. In the latter case, we impose Dirichlet or Neumann boundary condition on $\partial\Omega$.
Let $H = - \Delta_G$, where $\Delta_G$ is the Laplace-Beltrami operator associated with the metric $G$.  One can then define a scattering operator $\widehat S(\lambda) = \big(\widehat S_{ij}(\lambda)\big)$, which is a bounded operator on $L^2(M_1)\oplus \cdots \oplus L^2(M_N)$, where $\lambda \in (E_0,\infty)\setminus{\mathcal E}(H)$ is the energy parameter, $E_0 = \inf\sigma_{ess}(H)$, and ${\mathcal E}(H)$ is the set of exceptional points to be defined in (\ref{S3:ExcepPoints}).
Our goal is the following.

\begin{figure}[htbp]
\begin{center}
\psfrag{1}{$\Omega_1$}
\psfrag{2}{}
\psfrag{3}{$\Omega_2$}
\psfrag{4}{$\Omega_3$}
\psfrag{5}{\hspace{-5mm}$\mathcal K$}
\includegraphics[width=8cm]{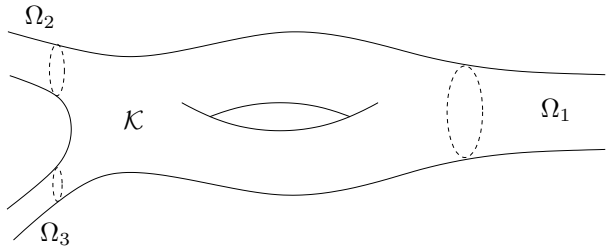} 
\label{pic 1}
\end{center}
\caption{Manifold $\Omega$ has ends $\Omega_j$, $j=1,2,\dots,N$.
}
\end{figure}


\begin{theorem}\label{main thm}
Suppose we are given two manifolds $\Omega^{(r)}, \, r = 1,2$, of the form (\ref{S1:Omega}) having $N_r$ ends, $\Omega_i^{(r)}$, $i = 1,\cdots,N_r$,
equipped with the metric $G^{(r)}$  satisfying the assumption (\ref{eq:Sec1DecayAssumSR}). 
Assume that $\Omega_{1}^{(1)} = \Omega_{1}^{(2)}$ and 
\begin{equation}
G^{(1)}_1 = G^{(2)}_1 = (dy)^2 + h_1(x,dx),
\quad h_1(x,dx)=\sum_{j,k=1}^{n-1} h_{1,jk}(x)dx^jdx^k
\label{eq:equalmetriconOmega1}
\end{equation}
 on $\Omega_{1}^{(1)} = \Omega_{1}^{(2)}$, moreover $\widehat S^{(1)}_{11}(\lambda) = \widehat S^{(2)}_{11}(\lambda)$ for all $\lambda \in (E',\infty)\setminus({\mathcal E}^{(1)}\cup{\mathcal E}^{(2)})$, where ${\mathcal E}^{(r)}$ is the set of exceptional points for $H^{(r)}$, and $E' = \max\,(E_0^{(1)},E_0^{(2)})$. Then $\Omega^{(1)}$ and $\Omega^{(2)}$ are isometric as Riemannian manifolds with metrics $G^{(1)}$, $G^{(2)}$.
\end{theorem}

This means that if we observe waves coming in and going out of one end $\Omega_{1}$, which is assumed to be non-perturbed, we can identify the whole manifold $\Omega$.
Note that in Theorem \ref{main thm}, neither the number of ends of each $\Omega^{(r)}$ nor the metric on the manifold $M_i^{(r)}$ are assumed to be known a-priori.  The key idea of the proof is to introduce generalized eigenfunctions of the Laplace-Beltrami operator which are exponentially growing at infinity, and define the associated {\it non-physical} scattering amplitude. The crucial fact is that this non-physical scattering amplitude is the analytic continuation of the {\it physical} scattering amplitude. Then the physical scattering amplitude determines the non-physical scattering amplitude, which further determines the Neumann-Dirichlet map of the {\it interior} domain. By the boundary control method (see \cite{AKKLT,Be87,Be97,KKL01,KL00,KL02}), one can determine the metric inside. 

In this paper, we exclusively deal with the Neumann boundary condition. The other cases are treated similarly and in fact more easily. The forward problem, i.e., the existence and completeness of wave operators, the eigenfunction expansion and the representation of S-matrices are well-known for short-range perturbations (see e.g. \cite{Go69}, \cite{Go74}, \cite{Lyf75}, \cite{Lyf76}, \cite{Wil76}, \cite{GuWi78}, \cite{Wil83}, see also \cite{Me95}). The new issue we have to discuss in this paper is the difference of conormal derivatives on the boundary associated with unperturbed and perturbed metrics. Therefore, focusing on this point, we only explain the outline of the proof of the forward problem under the assumption (\ref{eq:Sec1DecayAssumSR}) following the approach in \cite{IsKu08}, where spectral theory and inverse problems on hyperbolic spaces are developed in an elementary way. 

In the cylindrical ends, the physical generalized eigenfunction of the Laplace-Beltrami operator admits the analytic continuation with respect to the energy parameter, and this analytically continued eigenfunction is exponentially growing as $y \to \infty$. This sort of non-physical exponentially growing generalized eigenfunction was first introduced by Faddeev to develop the multi-dimensional Gel'fand-Levitan theory (\cite{Fa74}). The exponentially growing solutions of Schr{\"o}dinger equation was rediscovered in 1980's, where the new feature is its utility in the inverse scattering from a fixed energy as well as in the inverse boundary value problem (see \cite{NaAb84}, \cite{SU87}, \cite{Na88}, \cite{HN87} and the review \cite{GKLU09}). The interesting fact is that this apparently mysterious exponentially growing generalized eigenfunctions appear naturally in the cylindrical domain. 
Using these exponentially growing eigenfunctions, it is possible to obtain, from the scattering operator $\widehat S_{11}(\lambda)$, the Gel'fand spectral data on a part of the boundary $\Gamma = M_1\times\{1\}$ of the non-compact manifold
$\Omega^1 = \Omega\setminus\left(M_1\times(1,\infty)\right)$. The Gel'fand boundary data for this case is the family of the Neumann-Dirichlet map, $\Lambda(z)$, $\Lambda(z)f = u\big|_{\Gamma}$, where $u$ is the solution to the boundary value problem
\begin{equation}
\left\{
\begin{split}
& ( - \Delta_{G^1}-z)u = 0 \quad {\rm in} \quad \Omega^1, \\
& \partial_{\nu}u = 0 \quad {\rm on} \quad \partial\Omega\cap\partial\Omega^1, \\
& \partial_{\nu}u = \partial_{\nu}y = f \quad {\rm on} \quad \Gamma.
\end{split}
\right.
\nonumber
\end{equation}
To solve this problem, we use the boundary control (BC) method (see \cite{Be87} for the pioneering work and \cite{Be97}, \cite{KKL01} for the detailed exposition).  We note that typically the BC method deals with inverse problem on compact manifolds. The case of non-compact manifold considered here requires substantial modifications into the method, since the spectrum is no more discrete and it is also impossible  to use eigenfunctions as coordinate functions. A short description of the BC-method for non-compact manifolds was given in \cite{KKL04}. Here we provide detailed constructions for the considered case of a manifold with asymptotically cylindrical ends.

The structure of this paper is as follows. Sections 2, 3, 4 are devoted to a detailed analysis of scattering on manifolds with asymptotically cylindrical ends. After some preliminary estimates for the case of a half-cylinder with a product metric in \S 2, we discuss the spectral properties of the Laplacin
in $\Omega$ in \S 3. Using these properties, we develop the scattering theory for such manifolds in \S 4. The remaining part of the paper, Sections 5, 6 are devoted to the inverse scattering. In \S 4, we show that $\widehat S_{11}(\lambda)$ determines the Neumann-Dirichelt map $\Lambda(z)$.
An important step, which at the moment requires the product structure of the metric on $M_1\times(0,\infty)$, is the recovery, from physical scattering matrix $\widehat S_{11}(\lambda)$, the non-physical scattering amplitude. At last, \S 6 is devoted to the development of the B-method for non-compact manifolds. For the convenience of the reader, interested predominantly in the inversion methods, we make this section independent of the previous ones.

 Our manifold $\Omega$ is a mathematical model of compound waveguides, e.g. settings of optical and electric cables, oil, gas and water pipelines, etc, which are the most typical geometric constructions encountered in
the every-day life.  {As for the inverse problem, many works have been devoted so far to the distribution of resonances for the waveguides (\cite{Chr95}, \cite{ChrZw95}, \cite{BGRS97}, \cite{DP98}, \cite{DEM98}, \cite{APV00}, \cite{Chr02}). Identification or reconstruction of the domain or the medium for grating, layers or waveguides are studied by
\cite{CGI02}, \cite{IMN04}, \cite{NaSi04}, \cite{ERY08}. In particular, a similar inverse problem for waveguides was considered by Eskin-Ralston-Yamamoto \cite{ERY08} when $\Omega$ is a slab, $(0,B)\times\R$, with the variable sound speed $c(x,y)$, where $c(x,y) = c(x)$ for large $|y|$. The present paper deals with the forward and inverse scattering problems for waveguide in a full generality.}

The notation in this paper is standard. For a self-adjoint operator $A$, $\sigma(A)$, $\sigma_p(A)$ and $\sigma_{ess}(A)$ mean its spectrum, point spectrum and essential spectrum, respectively. For two Banach spaces $\mathcal H_1, \mathcal H_2$, ${\bf B}({\mathcal H_1};{\mathcal H_2})$ means the space of all bounded operators from $\mathcal H_1$ to $\mathcal H_2$. For an operator $A$ on a Hilbert space $\mathcal H$, $D(A)$ denotes its domain of definition. For a Riemannian manifold $\mathcal M$, $H^m(\mathcal M)$ denotes the usual Sobolev space of order $m$ on $\mathcal M$. For a domain $D$ and a Hilbert pace $\mathcal H$, $L^2(D;{\mathcal H};d\mu)$ means the space of $\mathcal H$-valued $L^2$-functions on $D$ with respect to the measure $d\mu$. If $\mathcal H = {\C}$, we omit it. {For a differentiable manifold $M$ and $p \in M$, $T_p(M)$ denotes the tangent space of $M$ at $p$.} A simplified version  of our results is given in \cite{IsKuLa07}.


\section{A-priori estimates in half-cylinders}


As a preliminary, let us begin with proving some a-priori estimates for the operator $- \partial_y^2 - \Delta_h$
on $\Omega_0 = M \times {\R_+}$ with Neumann boundary condition, where $y \in R_+ = (0,\infty)$, $M$ is a compact Riemannian manifold, and $\Delta_h$ is the Laplace-Beltrami operator associated with metric $h(x,dx)$ equipped on $M$. 

\subsection{Besov type spaces on cylinder}
We define an abstract Besov type space, which was introduced  by H{\"o}rmander \cite{AgHo76} in the case of ${\R}^n$. 
Let $M$ be the above mentioned compact manifold, and $(\, , \,)_{M}$, $\|\cdot\|_{L^2(M)}$ be inner product and norm of $L^2(M)$, respectively.
 We define intervals $I_n$ by
\begin{equation}
I_n =
\left\{
\begin{split} 
& (2^{n-1},2^n], \quad n \geq 1, \\
& (0,1], \quad n = 0.
\end{split}
\right.
\nonumber
\end{equation}
 Let $\mathcal B$ be the Banach space of $L^2(M)$-valued functions on $(0,\infty)$ equipped with norm 
\begin{equation}
\|f\|_{\mathcal B} = \sum_{n =0}^{\infty}2^{n/2}
\left(\int_{I_n}\|f(y)\|_{L^2(M)}^2dy \right)^{1/2}.
\nonumber
\end{equation}
Its dual space is the set of $L^2(M)$-valued functions $u(y)$ satisfying
\begin{equation}
\|u\|_{\mathcal B^{\ast}} = \sup_{n\geq 0}2^{-n/2}
\left(\int_{I_n}\|v(y)\|_{L^2(M)}^2dy\right)^{1/2} < \infty.
\nonumber
\end{equation}
It is easy to see that there exists a constant $C > 0$ such that\begin{equation}
\begin{split}
C^{-1}\sup_{n\geq 0}2^{-n/2}
\left(\int_{I_n}\|v(y)\|_{L^2(M)}^2dy\right)^{1/2} & \leq 
\left(\sup_{R>1}\frac{1}{R}
\int_{0}^R\|u(y)\|_{L^2(M)}^2dy\right)^{1/2} \\
& \leq C\sup_{n\geq0}2^{-n/2}
\left(\int_{I_n}\|v(y)\|_{L^2(M)}^2dy\right)^{1/2}.
\end{split}
\nonumber
\end{equation}
Therefore, we identify $\mathcal B^{\ast}$
 with the space equipped with norm
\begin{equation}
\|u\|_{{\mathcal B}^{\ast}} = \left(\sup_{R>1}\frac{1}{R}
\int_{0}^R\|u(y)\|_{L^2(M)}^2dy\right)^{1/2}.
\nonumber
\end{equation}
We also use the following weighted $L^2$ space and weighted Sobolev space: For $s \in {\R}$,
\begin{equation}
L^{2,s} \ni f \Longleftrightarrow 
\|f\|_s^2 = \int_{0}^{\infty}(1 + y)^{2s}\|f(y)\|_{L^2(M)}^2dy < \infty,
\nonumber
\end{equation}
\begin{equation}
H^{m,s} \ni u \Longleftrightarrow \|u\|_{H^{m,s}} = \|(1 + y)^su\|_{H^m(M\times(0,\infty))} < \infty.
\nonumber
\end{equation}
In the following, $\|\cdot\|$ means $\|\cdot\|_0$ and $(\cdot,\cdot)$ denotes the inner product of $L^2(M\times{\R}_+)$. It often denotes the coupling of two functions $f \in L^{2,s}$ and $g \in L^{2,-s}$ or $f \in \mathcal B$ and $g \in \mathcal B^{\ast}$.
The following inclusion relations can be shown easily, and the proof is omitted.


\begin{lemma} For $s > 1/2$, we have
\begin{equation}
L^{2,s} \subset {\mathcal B} \subset L^{2,1/2} \subset L^2 \subset 
L^{2,-1/2} \subset {\mathcal B}^{\ast} \subset L^{2,-s}.
\nonumber
\end{equation}
\end{lemma}

We often make use of the following lemma, whose proof is also elementary and omitted.


\begin{lemma}
Suppose $u \in \mathcal B^{\ast}$. Then 
\begin{equation}
\lim_{R\to\infty}\frac{1}{R}\int_{0}^R\|u(y)\|_{L^2(M)}^2dy = 0,
\label{eq:Integralmean0}
\end{equation}
if and only if 
\begin{equation}
\lim_{R\to\infty}\frac{1}{R}\int_0^{\infty} \rho\Big(\frac{y}{R}\Big)\|u(y)\|_{L^2(M)}^2dy = 0, \quad \forall \rho \in C_0^{\infty}((0,\infty)).
\label{eq:Integralmean0byrho}
\end{equation}
\end{lemma}


\subsection{A-priori estimates}
Let us consider the following equation in $\Omega_0 = M\times{\R}_+$:
\begin{equation}
\left\{
\begin{split}
& (- \partial_y^2 - \Delta_h - z)u = f \quad {\rm in} \quad \Omega_0\\
& \partial_{\nu}u = 0 \quad {\rm on} \quad \partial\Omega_0,
\end{split}
\right.
\label{eq:Sect2Helmholtzeq}
\end{equation}
$z$ being a complex parameter, and $\partial_{\nu}$ conormal differentiation on the boundary.  In the following, we often denote by $\|\partial_x^{\alpha}u\|$ the norm of derivatives of $|\alpha|$-th order of $u$ without mentioning local coordinates.


\begin{lemma}
Let $z \in {\C}$ be given. Then : \\
\nonumber
 (1) If $u, f \in L^{2,s}$ for some $s \in {\R}$, we 
have
\begin{equation}
\sum_{|\alpha| + l \leq 2}\|\partial_x^{\alpha}\partial_y^l u\|_{s} 
\leq C(\|u\|_s + \|f\|_s).
\nonumber
\end{equation}
(2) If $u, f \in \mathcal B^{\ast}$, then we have 
\begin{equation}
\|\partial_xu\|_{\mathcal B^{\ast}} + \|\partial_yu\|_{\mathcal B^{\ast}} \leq C (\|u\|_{\mathcal B^{\ast}} + \|f\|_{\mathcal B^{\ast}}).
\nonumber
\end{equation}
\end{lemma}
\noindent \noindent {\bf Proof.} We shall prove (2). Pick $\chi(y) \in C_0^{\infty}({\R})$ such that $\chi(y) = 1 \ (|y| < 1)$, $\chi(y) = 0 \ (|y| > 2)$ and put $\chi_R(y) = \chi(y/R)$. We take the inner product in $L^2(\Omega_0)$ of (\ref{eq:Sect2Helmholtzeq}) and $\chi_R^2(y)u$.
We then have
\begin{equation}
\|\chi_R\partial_yu\|^2 + \big(\chi_R\partial_yu,\frac{2}{R}\chi'\big(\frac{y}{R}\big)u\big) + \|\chi_R\partial_xu\|^2 - z\|\chi_Ru\|^2 = (f,\chi_R^2u),
\nonumber
\end{equation}
which implies
\begin{equation}
\|\chi_R\partial_yu\|^2 + \|\chi_R\partial_xu\|^2 
\leq C\left(\frac{1}{R^2}\|\chi'\big(\frac{y}{R}\big)u\|^2 + \|\chi_Ru\|^2 + 
 \|\chi_Rf\|^2\right).
\nonumber
\end{equation}
Then we have for $R > 1$
\begin{eqnarray*}
& & \int_0^R\|\partial_yu\|^2_{L^2(M)}dy + \int_0^R\|\partial_xu\|_{L^2(M)}^2dy \\
&\leq& C\left(\int_0^{2R}\|u\|^2_{L^2(M)}dy + 
 \int_0^{2R}\|f\|^2_{L^2(M)}dy\right).
\end{eqnarray*}
Dividing by $R$ and taking the supremum with respect to $R$, we obtain (2).

Let us prove (1). The 1st order derivatives are dealt with in the same way as above. We put $v = (1 + y)^{s}u$. Then $v$ satisfies $(- \partial_y^2 - \Delta_h - z)v = g$, where $g \in L^2(\Omega)$. By the a-priori estimates for elliptic operators, we have $v \in H^2(\Omega)$, which proves (1).
 \qed

\bigskip
Let $\lambda_1 < \lambda_2 \leq \cdots \to \infty$ be the eigenvalues of $- \Delta_h$, and $P_n$ the associated eigenprojection. Then
\begin{equation}
\sqrt{z + \Delta_h} = \sum_{n=1}^{\infty}\sqrt{z - \lambda_n}\,P_n,
\nonumber
\end{equation}
where for $\zeta = re^{i\theta},\ ( r > 0, 0 < \theta < 2\pi)$, we define $\sqrt{\zeta} = \sqrt{r}e^{i\theta/2}$. 

Our next aim is to derive some a-priori estimates for solutions to the equation (\ref{eq:Sect2Helmholtzeq}). We use the method of integration by parts due to Eidus (\cite{Eid65}). 
We put 
\begin{equation}
P(z) = \sqrt{z + \Delta_h},
\nonumber
\end{equation}
\begin{equation}
D_{\pm}(z) = \partial_y \mp iP(z).
\nonumber
\end{equation}
Then the equation (\ref{eq:Sect2Helmholtzeq}) is rewritten as
\begin{equation}
\partial_y D_{\pm}(z)u = \mp i P(z)D_{\pm}(z)u - f.
\label{eq:EquationofDpm(z)u}
\end{equation}


\begin{lemma}
Let $\varphi(y) \in C^{\infty}({\R})$ be such that $\varphi(y) \geq 0$. For a solution $u$ of the equation (\ref{eq:Sect2Helmholtzeq}), we put $w = D_+(z)u$. Then if ${\rm Im}\,z \geq 0$ we have for any $0 < a < b < \infty$
\begin{equation}
\int_a^b\varphi'(y)\|w(y)\|_{L^2(M)}^2dy \leq 2\int_a^b\varphi(y)\big|(f,w)_{L^2(M)}\big|dy + 
\Big[\varphi\|w\|^2_{L^2(M)}\Big]_{y=a}^{y=b}.
\nonumber
\end{equation}
\end{lemma}
\noindent {\bf Proof.} Since $w$ satisfies $\partial_y w = - i P(z)w - f$, we have
\begin{equation}
\int_a^b\varphi(y)(\partial_y w,w)_{L^2(M)}dy = - i\int_a^b\varphi(y)(P(z)w,w)_{L^2(M)}dy - 
\int_a^b\varphi(y)(f,w)_{L^2(M)}dy.
\nonumber
\end{equation}
Integrating  by parts and taking the real part, we have
\begin{eqnarray*}
& & \Big[\varphi\|w\|_{L^2(M)}^2\Big]_a^b - \int_a^b\varphi'(y)\|w(y)\|_{L^2(M)}^2dy \\
&=& 2\int_a^b\varphi(y)({\rm Im}\,P(z)w,w)_{L^2(M)}dy -
2 {\rm Re}\int_a^b\varphi(y)(f,w)_{L^2(M)}dy.
\end{eqnarray*}
Taking notice of ${\rm Im}\,P(z) \geq 0$ for ${\rm Im}\,z \geq 0$, we get the lemma. \qed

\medskip
Let ${\C}_+ = \{z \in {\C}\, ; {\rm Im}\,z \geq 0\}$.


\begin{lemma}
Let $w$ be as in Lemma 2.4 and suppose that
\begin{equation}
\lim_{R\to\infty}\frac{1}{R}\int_1^R\|w(y)\|_{L^2(M)}^2dy = 0.
\label{eq:limRtoinftywy}
\end{equation}
Then there exists a constant $C > 0$ independent of $z \in {\C}_+$ such that
$$
\|w(y)\|_{L^2(M)}^2 \leq C\|f\|_{\mathcal B}\|w\|_{{\mathcal B}^{\ast}}, \quad
\forall y \in {\R}.
$$
\end{lemma}
\noindent \noindent {\bf Proof.} Taking $\varphi(y) = 1$ in Lemma 2.4, we have
\begin{equation}
\begin{split}
\|w(a)\|_{L^2(M)}^2 & \leq \|w(b)\|_{L^2(M)}^2 + 2\int_a^b|(f,w)_{L^2(M)}|dy \\
&\leq \|w(b)\|^2_{L^2(M)} + C\|f\|_{\mathcal B}\|w\|_{{\mathcal B}^{\ast}}.
\end{split}
\nonumber
\end{equation}
The assumption of the lemma implies 
$\displaystyle{\liminf_{b\to\infty}\|w(b)\|_{L^2(M)} = 0}$, 
which proves the lemma. \qed


\begin{cor} Under the assumption of Lemma 2.5,
there exists a constant $C > 0$ such that
\begin{equation}
\|w\|_{{\mathcal B}^{\ast}} \leq C\|f\|_{\mathcal B}, \quad
\forall z \in {\C}_+.
\nonumber
\end{equation}
\end{cor}
\noindent \noindent {\bf Proof.} Lemma 2.5 implies that 
\begin{equation}
\|w\|^2_{\mathcal B^{\ast}} = \sup_{R>1}\frac{1}{R}\int_{0}^R\|w(y)\|^2_{L^2(M)}dy \leq C\|f\|_{\mathcal B}\|w\|_{\mathcal B^{\ast}},
\nonumber
\end{equation}
which proves this corollary. \qed

\begin{theorem}
For a small $\delta > 0$, let 
$$
J_{\delta} = \{z \in {\C}_+\;; \, {\rm dist}\,\big({\rm Re}\,z,\sigma(-\Delta_h)\big) > \delta\}.
$$
Let $u$ be a solution to (\ref{eq:Sect2Helmholtzeq}) such that $w = D_+(z)u$ satisfies (\ref{eq:limRtoinftywy}). 
Then there exists a constant $C > 0$ such that
$$
\|u\|_{\mathcal B^{\ast}} \leq C\|f\|_{\mathcal B}
$$
holds uniformly for $z \in J_{\delta}$.
\end{theorem}
\noindent \noindent {\bf Proof.} 
Let $A(z) = {\rm Re}\,P(z) = (P(z) + P(z)^{\ast})/2$. By the equation (\ref{eq:EquationofDpm(z)u}), we have
$$
 \partial_y(w,u)_{L^2(M)} 
 = - i(P(z)w,u)_{L^2(M)} - (f,u)_{L^2(M)} + (w,\partial_yu)_{L^2(M)}.
$$
In view of the formula
\begin{equation}
\begin{split}
-i(P(z)w,u)_{L^2(M)} & = - 2i(A(z)w,u)_{L^2(M)} + i(P(z)^{\ast}w,u)_{L^2(M)} \\
& = - 2i(w,A(z)u)_{L^2(M)} + i(w,P(z)u)_{L^2(M)}, 
\end{split}
\nonumber
\end{equation}
we then have
$$
 \partial_y(w,u)_{L^2(M)} = - 2i(w,A(z)u)_{L^2(M)}  - (f,u)_{L^2(M)} + \|w\|_{L^2(M)} ^2. 
$$
Using $w = \partial_yu - iP(z)u$, we compute
$$
2i(w,A(z)u)_{L^2(M)} = 2i(\partial_yu,A(z)u)_{L^2(M)} + \|P(z)u\|^2_{L^2(M)} 
+ (P(z)^2u,u)_{L^2(M)}.
$$
Summing up, we have arrived at
\begin{equation}
\begin{split}
 \partial_y(w,u)_{L^2(M)} 
 & = -2i(\partial_yu,A(z)u)_{L^2(M)}  - \|P(z)u\|^2_{L^2(M)}  \\
 & \ \ \ \ - ((z + \Delta_h)u,u)_{L^2(M)} 
 - (f,u)_{L^2(M)}  + \|w\|^2_{L^2(M)} .
\end{split}
\nonumber
\end{equation}
Taking the imaginary part and integrating in $y$, we have
\begin{equation}
\begin{split}
{\rm Im}\,\Big[(w,u)_{L^2(M)} \Big]_{y=a}^{y=b} & = 
- 2{\rm Re}\,\int_a^b(\partial_yu,A(z)u)_{L^2(M)} \\
& \ \ \ \  - {\rm Im}\,z\int_a^b\|u\|^2_{L^2(M)}dy
- {\rm Im}\,\int_a^b(f,u)_{L^2(M)} dy.
\end{split}
\nonumber
\end{equation}
Since $A(z)$ is self-adjoint, we have by integration by parts
\begin{equation}
2{\rm Re}\,\int_a^b(\partial_yu,A(z)u)_{L^2(M)} dy = \Big[(A(z)u,u)_{L^2(M)} \Big]_{y=a}^{y=b}.
\nonumber
\end{equation}
Using ${\rm Im}\,z \geq 0$, we obtain
\begin{equation}
{\rm Im}\,\Big[(w,u)_{L^2(M)} \Big]_{y=a}^{y=b} + \Big[(A(z)u,u)_{L^2(M)} \Big]_{y=a}^{y=b}
\leq C\|f\|_{\mathcal B}\|u\|_{\mathcal B^{\ast}},
\label{eq:wufidenitity}
\end{equation}
where $C$ is independent of $z \in {\C}_+$. 
We renumber the eigenvalues of $- \Delta_h$ in the increasing order $\mu_1 < \mu_2 < \cdots$ without counting multiplicities and put $\mu_0 = - \infty$, i.e. $\{\lambda_n \, ; n = 1, 2, \cdots\}$ and $\{\mu_n\, ; n = 1,2\cdots\}$ are the same as subsets of ${\R}$.
For a sufficiently small $\delta > 0$, we put
\begin{equation}
J_{n,\delta} = \{z \in {\C}_+ \; ; \ \mu_{n-1} + \delta < {\rm Re}\,z < \mu_n - \delta\}.
\nonumber
\end{equation}
Assume $z \in J_{n,\delta}$ and split $u$ as $u = u_{<} + u_{>}$, where
\begin{equation}
u_{<} = \sum_{\lambda_j \leq \mu_{n-1}} P_ju, \quad
u_{>} = \sum_{\lambda_j \geq \mu_n}P_ju,
\nonumber
\end{equation}
Recall that $P_j$ is the eigenprojection associated with $\lambda_j$.
We also define $w_<$, $w_>$, $f_<$, $f_>$ similarly. Note that $w_{<} = D_+(z)u_{<}$. Let us remark that (\ref{eq:Sect2Helmholtzeq}) and therefore (\ref{eq:wufidenitity}) hold with $w, u, f$ replaced by $w_<, u_<, f_<$ and $w_>, u_>, f_>$, respectively.
For eigenvalues $\lambda_j \leq \mu_{n-1}$, we have ${\rm Re}\,\sqrt{z-\lambda_j} \geq \sqrt{\delta}$. Therefore
\begin{equation}
(A(z)u_<,u_<)_{L^2(M)}  \geq \sqrt{\delta}\|u_<\|_{L^2(M)} ^2.
\label{eq:(A(z)u<u<)}
\end{equation}
Since $\partial_yu(0) = 0$, we have $w_{<}(0) = - iP(z)u_{<}(0)$. Therefore
\begin{equation}
\begin{split}
- {\rm Im}\,(w_{<}(0),u_{<}(0))_{L^2(M)} & = {\rm Re}\,(P(z)u_{<}(0),u_{<}(0))_{L^2(M)} \\
&= (A(z)u_{<}(0),u_{<}(0))_{L^2(M)}.
\end{split}
\nonumber
\end{equation}
Letting $a = 0$, $b = t$ in (\ref{eq:wufidenitity}), we then have
$$
{\rm Im}\,(w_{<}(t),u_{<}(t))_{L^2(M)} + 
(A(z)u_{<}(t),u_{<}(t))_{L^2(M)} \leq C\|f\|_{\mathcal B}\|u\|_{\mathcal B^{\ast}}.
$$
Using (\ref{eq:(A(z)u<u<)}), we have
\begin{equation}
\|u_<(t)\|_{L^2(M)} ^2 \leq C\big(\|w_<(t)\|_{L^2(M)} ^2 + \|f_<\|_{\mathcal B}\|u_<\|_{\mathcal B^{\ast}}\big).
\nonumber
\end{equation}
Using Corollary 2.6, we then have
\begin{equation}
\frac{1}{R}\int_0^R\|u_<(y)\|_{L^2(M)} ^2 dy \leq C\left(\|f_<\|_{\mathcal B}^2 + \|f_<\|_{\mathcal B}\|u_<\|_{\mathcal B^{\ast}}\right),
\nonumber
\end{equation}
which implies
\begin{equation}
\|u_<\|_{\mathcal B^{\ast}} \leq C\|f_<\|_{\mathcal B}.
\label{eq:Estimateofu<}
\end{equation}
On the other hand, if $\lambda_j \geq \mu_n$, we have ${\rm Re}\,(\lambda_j - z) \geq \delta$. Therefore
\begin{equation}
(- \partial_y^2 - \Delta_h - z)u_{>} = (-\partial_y^2 + B_z - i {\rm Im}\,z)u_{>} = f_{>},
\label{S2:Estimtateu>}
\end{equation}
where $B_z$ is a uniformly, with respect to $z$, strictly positive operator on $L^2(M)$.
Hence, we have 
\begin{equation}
\|u_>\|_{L^2} \leq C\|f_>\|_{L^2},
\label{Estimatesofu>inL2}
\end{equation}
which by Lemma 2.1 implies
\begin{equation}
\|u_>\|_{\mathcal B^{\ast}} \leq C\|f_>\|_{\mathcal B}.
\label{eq:Estimateofu>}
\end{equation}
The above two inequalities (\ref{eq:Estimateofu<}) and (\ref{eq:Estimateofu>}) prove the theorem. \qed


\section{Manifolds with cylindrical ends}


\subsection{Resolvent equation}
We return to the manifold $\Omega = \mathcal K\cup\Omega_1\cup\cdots\cup\Omega_N$ introduced in \S 1.
Fix a point $P_0 \in \mathcal K$ arbitrarily, and let ${\rm dist}(P,P_0)$ be the geodesic distance with respect to the metric $G$ from $P_0$ to $P$. We put 
$$
\Omega_0(R) = \{P \in \Omega \, ;\, {\rm dist}\,(P,P_0) < R\}, \quad
\Omega_{\infty}(R) = \{P \in \Omega \, ;\, {\rm dist}\,(P,P_0) \geq R\}.
$$
For $R > 0$ large enough, take $\chi_0 \in C_0^{\infty}(\Omega)$ such that $\chi_0 = 1$ on $\Omega_0(R)$, $\chi_0 = 0$ on $\Omega_{\infty}(R+1)$. Define $\chi_j = 1 - \chi_0$ on $\Omega_j$, $\chi_j = 0$ on $\Omega\setminus\Omega_j$. Then $\{\chi_j\}_{j=0}^{N}$ is a partition of unity on $\Omega$. 

Let $\Delta_G$ be the Laplace-Beltrami operator for the metric $G$ on $\Omega$ endowed with Neumann boundary condition on $\partial\Omega$. The conormal differentiation with respect to $G$ is denoted by $\partial_{\nu}$. We put
 $$
 H = - \Delta_G, \quad R(z) = (H - z)^{-1}.
 $$
As in \S 1, we identify $\Omega_j$ with $M_j\times(0,\infty)$, and let $h_j(x,dx)$ be the metric on $M_j$.
We compare $G$ with the unperturbed metric $G_{j}^{(0)} = (dy)^2 + h_j(x,dx)$ on $\Omega_j$. Let $\Delta_{G_{j}^{(0)}}$ be the Laplace-Beltrami operator for $G_{j}^{(0)}$ with  Neumann boundary condition on $\partial\Omega_j$.  The associated conormal differentiation is denoted by $\partial_{\nu^{(0)}_j}$. We put 
 $$
 H_{j}^{(0)} = - \Delta_{G_{j}^{(0)}}, \quad R_{j}^{(0)}(z) = (H_{j}^{(0)}- z)^{-1}.
 $$

Our next concern is the difference between the boundary conditions for $H$ and $H_j^{(0)}$. We put for large $R > 0$
\begin{equation}
\partial\Omega_j(R) = \partial\Omega \cap \Omega_j \cap \Omega_{\infty}(R).
\nonumber
\end{equation}


\begin{lemma} There exists a real function $w(x,y) \in C^{\infty}(\Omega_j)$  such that 
\begin{equation}
\left\{
\begin{split}
& \partial_{\nu}w(x,y) = 0 \quad  {\rm on} \quad \partial\Omega_j(R)\\
& w(x,y) = y + O(y^{-1-\epsilon_0}) \quad {\rm as} \quad y \to \infty.
\end{split}
\right. 
\label{S3:Propertyw}
\end{equation}
\end{lemma}
\noindent \noindent {\bf Proof.} By the decay assumption (\ref{eq:Sec1DecayAssumSR}),  letting $w(x,y) = y + \widetilde w(x,y)$, we should have
$\partial_{\nu}\widetilde w = - \partial_{\nu}y = O(y^{-1-\epsilon_0}) \
{\rm on} \ \partial\Omega_j(R)$. 
Extending the vector field $\nu$ near the boundary and integrating along it, we get $\widetilde w = O(y^{-1-\epsilon_0})$. \qed

\medskip
 For $m \geq 0$ and $s \in {\R}$, we define the weighted Sobolev space on the boundary by
$$
\psi \in H^{m,s}(\partial\Omega_j(R)) \Longleftrightarrow (1 + y)^s\psi \in H^{m}(\partial\Omega_j(R)).
$$


\begin{lemma}
There exists an operator of extension $\widetilde{\mathcal E_j}$ such that
for $m \geq 1/2$ and $\psi \in H^{m}(\partial\Omega_j(R))$
\begin{equation}
\partial_{\nu}\widetilde{\mathcal E_j}\psi = 
\left\{
\begin{split}
& \psi \quad {\rm on} \quad \partial\Omega_j(R), \\
& 0 \quad {\rm on} \quad \Omega\setminus\left(\Omega_j\cap\Omega_{\infty}(R-1/2)\right),
\end{split} 
\right.
\label{S3:WidetildeEj}
\end{equation} 
\begin{equation}
{\rm supp}\,(\widetilde{\mathcal E_j}\psi) \subset \Omega_j\cap\Omega_{\infty}(R-1).
\label{S3:WitilEesupport}
\end{equation}
For $m \geq 1/2$ and $s \geq 0$, it satisfies
\begin{equation}
\widetilde{\mathcal E_j} \in {\bf B}(H^{m,s}(\partial\Omega_j(R));H^{m+3/2,s}(\Omega_j)).
\label{S3:Widetildeinweight}
\end{equation}
\end{lemma}
\noindent \noindent {\bf Proof.} Let $\mathcal M' =  \Omega_j\cap\Omega_{\infty}(R-2)$. We smoothly modify the corner of $\mathcal M'$, i.e. $\{P \in \Omega_j\cap\partial\Omega \, ; \, {\rm dist}(P,P_0) = R-2\}$, and let $\mathcal M$ be the resulting manifold. Let $\nu_{\mathcal M}$ be the unit outer normal to $\mathcal M$. By solving the elliptic boundary value problem
\begin{equation}
\left\{
\begin{split}
&(- \Delta_G + i)u = 0 \quad {\rm in} \quad \mathcal M, \\
&\partial_{\nu_{\mathcal M}}u = \psi \quad {\rm on} \quad \partial{\mathcal M},
\end{split}
\right.
\label{S3:EllipticBP}
\end{equation}
we define $\widetilde{\mathcal E}_j\psi = \widetilde \chi_ju$, where $\widetilde\chi_j \in C^{\infty}(\Omega_j)$ is  such that $\widetilde\chi_j$ = 1 on $\Omega_j\cap\Omega_{\infty}(R-1/4)$, $\widetilde\chi_j = 0$ on $\Omega\setminus\big(\Omega_j\cap\Omega_{\infty}(R-1/2)\big)$. It then satisfies (\ref{S3:WidetildeEj}), (\ref{S3:WitilEesupport}).
The property (\ref{S3:Widetildeinweight}) for $s = 0$ follows from the standard estimate for the elliptic boundary value problem. Let $ 0 < s \leq 1 + \epsilon_0$ and take $\psi \in H^{m,s}(\partial{\mathcal M})$. For the solution $u$ to the boundary value problem (\ref{S3:EllipticBP}), we put $u = (1 + w(x,y))^{-s}u'$ and $\psi = (1 + w(x,y))^{-s}\psi'$, where $w(x,y)$ is constructed in Lemma 3.1. Then $u'$ is a solution to the boundary value problem
\begin{equation}
\left\{
\begin{split}
&(- \Delta_G + L' + \kappa)u' = 0 \quad {\rm in} \quad \mathcal M, \\
&\partial_{\nu_{\mathcal M}}u' = \psi'' \quad {\rm on} \quad \partial{\mathcal M}, 
\end{split}
\right.
\nonumber
\end{equation}
where $\kappa > 0$ is sufficiently large, and $L'$ is a 1st order differential operator with bounded coefficients, and $\psi'' = \psi'$ on $\partial\Omega_j(R)$. Since the mapping $\psi'' \to u'$ is bounded from $H^m(\partial{\mathcal M})$ to $H^{m+3/2}(\mathcal M)$, we get (\ref{S3:Widetildeinweight}) with $0 < s \leq 1 + \epsilon_0$. Repeating this procedure, we can prove (\ref{S3:Widetildeinweight}) for all $s > 0$. \qed

\medskip
For $u \in H^2(\Omega_j)$ satisfying 
$\partial_{\nu_j^{(0)}}u = 0$ on $\partial\Omega_j(R)$, we have
\begin{equation}
\partial_{\nu}\left(\chi_ju\right) = w(x,y)^{-1-\epsilon_0}B_ju
 \quad {\rm on} \quad \partial\Omega_j(R),
\label{eq:Bj}
\end{equation}
where
\begin{equation}
B_j = w(x,y)^{1+\epsilon_0}\left(\chi_j\big(\partial_{\nu}-\partial_{\nu^{(0)}_j}\big) +
\big(\partial_{\nu}\chi_j\big)\right)
\label{S3:DiffOpBj}
\end{equation}
is a 1st order differential operator on $\partial\Omega_j(R)$ with bounded coefficients. 
We put
\begin{equation}
{\mathcal E}_j = w(x,y)^{-1-\epsilon_0}\widetilde{\mathcal E_j}.
\label{S3:DefineEj}
\end{equation}
 Then by  (\ref{S3:Propertyw}), (\ref{S3:WidetildeEj}) and (\ref{eq:Bj}), for $u \in H^2(\Omega)$ satisfying 
$\partial_{\nu_j^{(0)}}u = 0$ on $\partial\Omega_j(R)$ the following formula holds
\begin{equation}
\partial_{\nu}{\mathcal E}_jB_ju = \partial_{\nu}\left(\chi_ju\right) \quad 
{\rm on} \quad \partial\Omega_j(R).
\label{eq:ExtensionjbyEB}
\end{equation}
Moreover
\begin{equation}
y^{1+\epsilon_0}\mathcal E_jB_j \in {\bf B}(H^2(\Omega);H^2(\Omega))
\cap {\bf B}(H^{3/2}(\Omega);H^{3/2}(\Omega)).
\label{eq:EjBjbounded}
\end{equation}

Suppose $u$ satisfies 
\begin{equation}
\left\{
\begin{split}
& (- \Delta_{G^{(0)}_j} - z)u = f \quad {\rm in} \quad  \Omega_j, \\
& \partial_{\nu^{(0)}_j}u = 0 \quad {\rm on} \quad \Omega_j\cap\partial\Omega.
\end{split}
\right.
\nonumber
\end{equation}
Then by (\ref{eq:ExtensionjbyEB}), $v_j = \chi_ju - \mathcal E_jB_ju$ satisfies
\begin{equation}
\left\{
\begin{split}
&(- \Delta_G - z)v_j = \chi_jf + V_j(z)u \quad {\rm in} \quad \Omega,\\
& \partial_{\nu}v_j = 0 \quad {\rm on} \quad \partial\Omega,
\end{split}
\right.
\label{eq:H-zvjequation}
\end{equation}
where
\begin{equation}
V_j(z) = [-\Delta_{G},\chi_j] + \chi_j(\Delta_{G^{(0)}_j} - \Delta_{G}) + (\Delta_{G} + z)\mathcal E_jB_j.
\label{eq:Vjzformula}
\end{equation}


\begin{lemma} Let $\widetilde\chi_j \in C^{\infty}(\Omega)$ be such that $\widetilde\chi_j = 1$ on $\Omega_j\cap\Omega_{\infty}(R-1)$ and $\widetilde\chi_j = 0$ outside $\Omega_j\cap\Omega_{\infty}(R-2)$. Then for $z \not\in {\R}$, the following resolvent equations hold :
\begin{equation}
R(z)\chi_j = \Big(\chi_j - \mathcal E_jB_j - R(z)V_j(z)\Big)R_j^{(0)}(z)\widetilde\chi_j.
\label{eq:Resolventeq'}
\end{equation}
\begin{equation}
\chi_jR(z) = \widetilde\chi_jJ_j^{-1}R_{j}^{(0)}(z)J_j\Big(\chi_j - (\mathcal E_jB_j)^{\ast} - V_j(\overline{z})^{\ast}R(z)\Big),
\label{eq:Resolventeq}
\end{equation}
where $J_j = \big({\rm det}\,G/{\rm det}\, G^{(0)}_j\big)^{1/2}$, and the adjoint $\ast$ is taken with respect to the inner product of $L^2(\Omega)$ with volume element from the metric $G$. Moreover $R_j^{(0)}(z)J_j(\mathcal E_jB_j)^{\ast}$ and $R_j^{(0)}(z)J_jV_j(\overline{z})^{\ast}R(z)$ are compact on $L^2(\Omega)$.
\end{lemma}
\noindent \noindent {\bf Proof.} 
Let $u = R_j^{(0)}(z)\widetilde\chi_jf$ for $z \not\in {\R}$. Then checking the boundary condition by (\ref{eq:ExtensionjbyEB}), we have $v_j = \chi_jR_j^{(0)}(z)\widetilde\chi_jf - \mathcal E_jB_jR_j^{(0)}(z)\widetilde\chi_jf \in D(H)$, and by (\ref{eq:H-zvjequation}) $(H - z)v_j = \chi_j\widetilde\chi_jf + V_j(z)u = \chi_jf + V_j(z)u$, which implies
(\ref{eq:Resolventeq'}).

 By extending $f \in L^2(\Omega_j)$ to be 0 outside $\Omega_j$, we regard $L^2(\Omega_j)$ as a closed subspace of $L^2(\Omega)$.
The volume elements $dV$ and $dV_j^{(0)}$ of $G$ and $G_j^{(0)}$ satisfy $dV = J_jdV_j^{(0)}$. For $A \in {\bf B}(L^2(\Omega_j);L^2(\Omega_j))$, let $A^{\ast}$ and $A^{\ast(j)}$ denote their adjoint operators with respect to the volume element $dV$ and $dV_j^{(0)}$, respectively. Then it is easy to show that
$$
A^{\ast} = J_j^{-1}A^{\ast(j)}J_j.
$$
Taking $A = R_j^{(0)}(z)$, and noting that $R(z)^{\ast} = R(\overline{z})$ and 
$R_j^{(0)}(z)^{\ast(j)} = R_j^{(0)}(\overline z)$, we prove (\ref{eq:Resolventeq}). By (\ref{eq:EjBjbounded}) and (\ref{eq:Vjzformula}), $\mathcal E_jB_jR_j^{(0)}(z)$ and $R(z)V_j(z)J_jR_j^{(0)}(z)$ are compact on $L^2(\Omega)$, which implies the last assertion of the lemma. \qed


\subsection{Essential spectrum}


\begin{lemma} $\ $
$\sigma_{ess}(H) = [0,\infty)$.
\end{lemma}
\noindent \noindent {\bf Proof.}  
Lemma 3.3 implies
$\chi_jR(z) - \widetilde\chi_jJ_j^{-1}R_j^{(0)}(z)J_j\chi_j$ is compact. Therefore
\begin{equation}
R(z) = \sum_{j=1}^N\widetilde\chi_jJ_j^{-1}R_j^{(0)}(z)J_j\chi_j + K(z),
\label{S3Rzmodulocompact}
\end{equation}
where $K(z)$ is a compact operator and satisfies
\begin{equation}
\|K(z)\| \leq C|{\rm Im}\, z|^{-2}(1 + |z|), 
\label{S3Kzestimate}
\end{equation}
where $\|\cdot\|$ denotes the operator norm in $L^2(\Omega)$ and the constant $C$ is independent of $z$.
For $f(\lambda) \in C_0^{\infty}({\R})$, there exists $F(z) \in C_0^{\infty}({\C})$, called an almost analytic extension of $f$, such that $F(\lambda) = f(\lambda)$ for $\lambda \in {\R}$ and $|\overline{\partial_z}F(z)| \leq C_n|{\rm Im}\,z|^n$, $\forall n \geq 0$, and the following formula holds for any self-adjoint operator $A$ :
\begin{equation}
f(A) = \frac{1}{2\pi i}\int_{\C}\overline{\partial_z}F(z)(z - A)^{-1}dzd\overline{z}.
\label{S3:HelffeSjostrand}
\end{equation}
(See e.g. \cite{HeSj89} or  \cite{IsKu08}.) We replace $(z - A)^{-1}$ by $- R(z)$ and plug (\ref{S3Rzmodulocompact}). The inequality (\ref{S3Kzestimate}) implies $\|\overline{\partial_z}F(z)K(z)\| \leq C$, and the integral over $\C$ converges in the operator norm, hence it gives a compact operator. 
We then see that $\varphi(H) - \sum_{j=1}^N\widetilde\chi_jJ_j^{-1}\varphi(H_j^{(0)})J_j\chi_j$ is compact for any $\varphi(\lambda) \in C_0^{\infty}({\R})$. Since $\sigma(H_j^{(0)}) = [0,\infty)$, we have $\varphi(H_j^{(0)}) = 0$ if $\varphi(\lambda) \in C_0^{\infty}((-\infty,0))$. Therefore $\varphi(H)$ is compact if $\varphi(\lambda) \in C_0^{\infty}((-\infty,0))$, which implies that $(-\infty,0)\cap\sigma_{ess}(H) = \emptyset$. 
For $\lambda \in (0,\infty) = \sigma(H_j^{(0)})$, one can construct $u_n \in D(H_j^{(0)})$ such that $\|u_n\| = 1$, $\|(H_j^{(0)} - \lambda)u_n\| \to 0$, and ${\rm supp}\,u_n \subset \{y > R_n\}$ with $R_n \to \infty$. Then letting 
$v_n = \chi_ju_n - \mathcal E_jB_ju_n$, we have $v_n \in D(H)$, $\|(H - \lambda)v_n\| \to 0$, $v_n \to 0$ weakly and $\|v_n\| > C$ uniformly in $n$ with a constant $C > 0$. This implies $\lambda \in \sigma_{ess}(H)$. \qed

\bigskip
{\it The set of thresholds for  $H$} is defined by
\begin{equation}
\mathcal T(H) = {\mathop\bigcup_{j=1}^N}
 \sigma_p(- \Delta_{h_j}),
\label{S3Thresholds}
\end{equation}
where $\Delta_{h_j}$ is the Laplace-Beltrami operator on $M_j$.
Replacing $\Omega_0$ in \S 2 by $\Omega_j$ with $j = 1, \cdots, N$, we define the Besov type spaces $\mathcal B_j$, $\mathcal B_j^{\ast}$. We put
\begin{equation}
\|f\|_{\mathcal B} = \|\chi_0f\|_{L^2(\Omega)} + \sum_{j=1}^N\|\chi_jf\|_{\mathcal B_j},
\nonumber
\end{equation}
\begin{equation}
\|u\|_{\mathcal B^{\ast}} = \|\chi_0f\|_{L^2(\Omega)} + \sum_{j=1}^N\|\chi_jf\|_{\mathcal B_j^{\ast}}.
\nonumber
\end{equation}
The weighted $L^2$ space $L^{2,s}$ and the weighted Sobolev space $H^{m,s}$ are defined similarly.


\subsection{Radiation condition} A solution $u \in \mathcal B^{\ast}$ of the reduced wave equation
\begin{equation}
\left\{
\begin{split}
& (H - \lambda)u = f \quad {\rm in} \quad \Omega,  \quad \lambda > 0, \\
& \partial_{\nu}u = 0 \quad {\rm on} \quad \partial\Omega,
\end{split}
\right.
\nonumber
\end{equation}
 is said to satisfy the outgoing radiation condition if
\begin{equation} 
\lim_{R\to\infty}\frac{1}{R}\int_0^R
\|\chi_j\big(\partial_y - iP_j(\lambda)\big)u\|^2_{L^2(M_j)}dy = 0, \quad 1 \leq \forall j \leq N,
\label{S3:RadCond}
\end{equation}
where
\begin{equation}
P_j(z) = \sqrt{z + \Delta_{h_j}}.
\nonumber
\end{equation}
If $\partial_y - iP_j(\lambda)$ is replaced by $\partial_y + iP_j(\lambda)$, we say that $u$ satisfies the incoming radiation condition. In the following, $u$ is always assumed to satisfy the boundary condition $\partial_{\nu}u = 0$ on $\partial\Omega$.


\begin{lemma}
Let $\lambda \in (0,\infty)\setminus\mathcal T(H)$. If $u \in \mathcal B^{\ast}$ satisfies $(H - \lambda)u = 0$ and the outgoing (or incoming) radiation condition, it also satisfies
$$
\lim_{R\to\infty}\frac{1}{R}\int_0^R\|\chi_j u\|^2_{L^2(M_j)}dy = 0, \quad \ 1 \leq j \leq N.
$$
\end{lemma}
\noindent \noindent {\bf Proof.} We take $\rho(t) \in C_0^{\infty}((0,\infty))$ such that $\rho(t) \geq 0$, ${\rm supp}\,\rho(t) \subset (1,2)$ and $\int_{0}^{\infty}\rho(t)dt = 1$, and put
$$
\varphi_R(y) = \chi\left(\frac{y}{R}\right), \quad
\chi(t) = \int_{t}^{\infty}\rho(s)ds.
$$
Then $\varphi_R(y) = 1$ for $y < R$ and $\varphi_R(y) = 0$ for $y > 2R$. We next construct $\psi_R \in C_0^{\infty}(\Omega)$ such that $\psi_R = 1$ on $\mathcal K$ and $\psi_R = \varphi_R$ on $\Omega_j$ for $1 \leq j \leq N$. Then we have
$$
(i[H,\psi_R]u,u) = (i[H - \lambda,\psi_R]u,u) = 0.
$$
By the construction of $\psi_R$, $[H,\psi_R] = 0 $ on $\mathcal K$. By the assumption (\ref{eq:Sec1DecayAssumSR}), on $\Omega_j$ the commutator has the form
\begin{equation}
i[H,\psi_R] = \dfrac{2i}{R}\rho(\dfrac{y}{R})\partial_y + L_{j,R},
\label{S3:[H,psiR]}
\end{equation}
where $L_{j,R}$ is a 1st order differential operator whose coefficients have the form
$$
\frac{1}{R}\widetilde\chi\big(\frac{y}{R}\big)O(y^{-\epsilon_0})
$$
and $\widetilde\chi(y)$ is either $\rho(y)$ or $\rho'(y)$. 
Let $v = (1 + y)^{-\epsilon_0}u$. Then by Lemma 2.1 and Lemma 2.3 (1) (which also holds for $\Delta_{G_j}$), $\partial_xv, \partial_y v \in L^{2,-\delta}$ for some $0 < \delta < 1/2$. Therefore 
$$
\frac{1}{R}\int_0^R\left(\|\partial_xv\|_{L^2(M_j)}^2 + \|\partial_yv\|_{L^2(M_j)}^2\right)dy
\leq \frac{C}{R^{1-2\delta}}\left(\|\partial_xv\|_{-\delta}^2 + \|\partial_yv\|^2_{-\delta}\right),
$$
which tends to 0 as $R \to \infty$. Therefore by Lemma 2.2
\begin{equation}
\lim_{R\to\infty}(L_{j,R}u,u)_{L^2(\Omega_j)} = 0.
\nonumber
\end{equation}
Hence we have by using (\ref{S3:[H,psiR]}),
\begin{equation}
\lim_{R\to\infty}\sum_{j=1}^N\frac{1}{R}\int_0^{\infty}\rho\left(\frac{y}{R}\right)
\big(\partial_y\chi_ju,\chi_j u\big)_{L^2(M_j)}dy = 0.
\label{eq:limRrhochiu}
\end{equation}
Assume that $u$ satisfies the outgoing radiation condition. Using the inequality\begin{equation}
\begin{split}
& \left|\frac{1}{R}\int_0^{\infty}\rho\big(\frac{y}{R}\big)\big((\partial_y - iP_j(\lambda))\chi_ju,\chi_ju\big)_{L^2(M_j)}dy\right| \\
\leq \; & C\|u\|_{\mathcal B^{\ast}}\left(\frac{1}{R}\int_0^{\infty}\rho\big(\frac{y}{R}\big)\left\|(\partial_y - iP_j(\lambda))\chi_ju\right\|^2_{L^2(M_j)}dy\right)^{1/2}
\end{split}
\nonumber
\end{equation}
and (\ref{eq:limRrhochiu}), we then have
\begin{equation}
\lim_{R\to\infty}\sum_{j=1}^N\frac{1}{R}\int_0^{\infty}\rho\left(\frac{y}{R}\right)
\big(P_j(\lambda)\chi_ju,\chi_j u\big)_{L^2(M_j)}dy = 0.
\label{eq:limRrhochiuPj}
\end{equation}
As in the proof of Theorem 2.7, we split $\chi_ju$ into two parts, 
$$
\chi_ju_< = E_j((- \infty,\lambda))\chi_ju, \quad \chi_ju_> = E_j((\lambda,\infty))\chi_ju,
$$
where $E_j(\cdot)$ is the spectral projection associated with $- \Delta_{h_j}$.
Then by the short-range decay assumption of the metric,
$$
(- \partial_y^2 - \Delta_{h_j} - \lambda)\chi_ju_> =: f_j \in L^2(\Omega_j).
$$
Since $\lambda \not \in \sigma(- \Delta_{M_j})$, arguing in the same way as in the proof of (\ref{S2:Estimtateu>}), 
$$
(- \partial_y^2 - \Delta_{h_j} - \lambda)\chi_ju_{>} = 
(- \partial_y^2 + B_j)\chi_ju_{>},
$$
where $B_j$ is a self-adjoint operator on $L^2(M_j)$ such that
$B_j \geq \delta(1 - \Delta_{M_j})$,
 $\delta > 0$ being a constant. Therefore, 
$P_j(\lambda)\chi_ju_{>} \in L^2(\Omega_j)$, hence
\begin{equation}
\lim_{R\to\infty}\frac{1}{R}\int_0^{\infty}\rho\left(\frac{y}{R}\right)
\big(P_j(\lambda)\chi_ju_>,\chi_j u_>\big)_{L^2(M_j)}dy = 0.
\label{eq:limRchiu>}
\end{equation}
Since $P_j(\lambda)\chi_ju_{>} = iC_j(\lambda)\chi_ju_{>}$, where $C_j(\lambda)$ is a strictly positive operator on $L^2(M_j)$, this also implies
\begin{equation}
\lim_{R\to\infty}\frac{1}{R}\int_0^{\infty}\rho\left(\frac{y}{R}\right)
\|\chi_ju_>\|^2_{L^2(M_j)}dy = 0.
\label{eq:limRchiuPnasi>}
\end{equation}
We show that
\begin{equation}
\lim_{R\to\infty}\sum_{j=1}^N\frac{1}{R}\int\rho\left(\frac{y}{R}\right)
\big(P_j(\lambda)\chi_ju_<,\chi_j u_<\big)_{L^2(M_j)}dy = 0.
\label{eq:limRchiu<}
\end{equation}
In fact, in view of (\ref{eq:limRrhochiuPj}), splitting $u = u_{<} + u_{>}$ and using (\ref{eq:limRchiu>}), to prove (\ref{eq:limRchiu<}) we have only to show that
\begin{equation}
\lim_{R\to\infty}\frac{1}{R}\int\rho\left(\frac{y}{R}\right)
\big(P_j(\lambda)\chi_ju_>,\chi_j u_<\big)_{L^2(M_j)}dy = 0,
\label{eq:limRchiu<chiu>}
\end{equation}
and the same assertion with $u_{<}$ and $u_{>}$ exchanged. Let us note that
\begin{equation}
\begin{split}
|\big(P_j(\lambda)\chi_ju_>,\chi_j u_<\big)_{L^2(M_j)}| & = 
|\big(\chi_ju_>,\chi_j P_j(\lambda)^{\ast}u_<\big)_{L^2(M_j)}| \\
& \leq C\|\chi_ju_>\|\|\chi_ju_<\|.
\nonumber
\end{split}
\end{equation}
Therefore
\begin{equation}
\begin{split}
& \ \ \ \ \ \ \frac{1}{R}\int\rho\Big(\frac{y}{R}\Big)|\big(P_j(\lambda)\chi_ju_>,\chi_j u_<\big)_{L^2(M_j)}|dy \\
& \leq \left(\frac{1}{R}\int\rho\Big(\frac{y}{R}\Big)\|\chi_ju_>\|^2_{L^2(M_j)}dy\right)^{1/2}
\left(\frac{C}{R}\int\rho\Big(\frac{y}{R}\Big)\|\chi_j u_<\|^2_{L^2(M_j)}dy\right)^{1/2} \\
& \leq C\left(\frac{1}{R}\int\rho\Big(\frac{y}{R}\Big)\|\chi_ju_>\|^2_{L^2(M_j)}dy\right)^{1/2},
\nonumber
\end{split}
\end{equation}
since $\chi_ju_< \in \mathcal B^{\ast}$. By (\ref{eq:limRchiuPnasi>}), this converges to 0. Similarly, we can prove (\ref{eq:limRchiu<chiu>}) with $u_<$ and $u_>$ exchanged.

On the other hand, $(P_j(\lambda)\chi_ju_<,\chi_ju_<) \geq C\|\chi_ju_<\|^2_{L^2(\Omega_j)}$ for a constant $C > 0$, which depends on $\lambda$. Therefore by (\ref{eq:limRchiu<})
\begin{equation}
\lim_{R\to\infty}\frac{1}{R}\int_0^{\infty}\rho\left(\frac{y}{R}\right)
\|\chi_ju_<\|^2_{L^2(M_j)}dy = 0.
\label{eq:limRchiuPnasi<}
\end{equation}
By (\ref{eq:limRchiuPnasi>}) and (\ref{eq:limRchiuPnasi<}), we complete the proof of the lemma. \qed


\begin{lemma}
Suppose $u \in \mathcal B^{\ast}$ satisfies $(H - \lambda)u = f$ for $\lambda 
\in (0,\infty)\setminus{\mathcal T}(H)$ and $\partial_{\nu}u = 0$ on $\partial\Omega$. Assume also for some $1 \leq j \leq N$, $f \in L^{2,s}(\Omega_j)$  for any $s > 0$, and
\begin{equation}
\lim_{R\to\infty}\frac{1}{R}\int_{0}^R\|\chi_j u\|^2_{L^2(M_j)}dy = 0.
\label{eq:decayinOmegaj}
\end{equation}
Then $u \in L^{2,s}(\Omega_j)$ for any $s > 0$. Moreover for any $s > 0$ and any compact interval $I \subset (0,\infty) \setminus{\mathcal T(H)}$, there exists a constant $C_{s,I} > 0$ such that
\begin{equation}
\|\chi_ju\|_{L^{2.s}(\Omega_j)} \leq C_{s,I}(\|u\|_{\mathcal B^{\ast}(\Omega_j)} + \|f\|_{L^{2,s}(\Omega_j)}), \quad  \forall \lambda \in I.\label{S3chijuUnifEstimate}
\end{equation}
\end{lemma}
\noindent {\bf Proof.}
We construct counterparts of $\mathcal E_j$ and $B_j$ when the roles of 
$G$ and $G_j^{(0)}$ are interchanged. Namely, there exists an operator of extension $\widetilde{\mathcal E_j}^{(0)}$ such that
for $m \geq 1/2$ and $\psi \in H^{m}(\partial\Omega_j(R))$
\begin{equation}
\partial_{\nu_j^{(0)}}\widetilde{\mathcal E_j}^{(0)}\psi = 
\left\{
\begin{split}
& \psi \quad {\rm on} \quad \partial\Omega_j(R), \\
& 0 \quad {\rm on} \quad \Omega\setminus\left(\Omega_j\cap\Omega_{\infty}(R-1/2)\right),
\end{split} 
\right.
\nonumber
\end{equation} 
\begin{equation}
{\rm supp}\,(\widetilde{\mathcal E_j}^{(0)}\psi) \subset \Omega_j\cap\Omega_{\infty}(R-1),
\nonumber
\end{equation}
\begin{equation}
\widetilde{\mathcal E_j}^{(0)} \in {\bf B}(H^{m,s}(\partial\Omega_j(R));H^{m+3/2,s}(\Omega_j)), \quad m \geq 1/2, \quad s \geq 0. 
\nonumber
\end{equation}
If $\partial_{\nu}u = 0$ on $\partial\Omega$, we have
$$
\partial_{\nu^{(0)}_j}(\chi_ju) = y^{-1-\epsilon_0}B_j^{(0)}u \quad
{\rm on} \quad \partial\Omega_j(R),
$$
where
$$
B_j^{(0)} = y^{1+\epsilon_0}\chi_j(\partial_{\nu^{(0)}_j} - \partial_{\nu}).
$$
 We put
$$
\mathcal E_j^{(0)} = y^{-1-\epsilon_0}\widetilde{\mathcal E}_j^{(0)}.
$$
Then 
$$
\partial_{\nu^{(0)}_j}\mathcal E_j^{(0)}B_j^{(0)}u = \partial_{\nu^{(0)}_j}
(\chi_ju).
$$

Suppose $u \in \mathcal B^{\ast}$ satisfies $(H - \lambda)u = f$, $\lambda \in (0,\infty)\setminus{\mathcal T}(H)$, and $\partial_{\nu}u = 0$ on $\partial\Omega$. We put $v_j = \chi_ju - \mathcal E_j^{(0)}B_j^{(0)}u$. Then $v_j$ satisfies 
\begin{equation}
\left\{
\begin{split}
& (- \partial_y^2 - \Delta_{h_j} - \lambda)v_j = \chi_jf + L_{j}u + 
(\partial_y^2 + \Delta_{h_j})\mathcal E_j^{(0)}B_j^{(0)}u \quad 
{\rm in} \quad \Omega_j, \\
& \partial_{\nu_j^{(0)}}v_j = 0 \quad {\rm on} \quad \partial\Omega_j.
\end{split}
\right.
\label{Lemma36chijueq} 
\end{equation}
Here $L_{j}$ is a 2nd order differential operator with coefficients decaying like $O(y^{-1-\epsilon_0})$. 
Note that $ f_j := \chi_jf + L_{j,n}u + 
(\partial_y^2 + \Delta_{h_j})\mathcal E_j^{(0)}B_j^{(0)}u \in L^{2,1+\epsilon_0}(\Omega_j)$. 

 Let $v_{j,n} = (v_j(\cdot,y),\psi_n(\cdot))_{L^2(M_j)}$, where $\psi_n(x)$ is the normalized eigenvector associated with the eigenvalue $\lambda_n$ of $- \Delta_{h_j}$. Then we have
\begin{equation}
(- \partial_y^2  - \mu_n)v_{j,n} = g_{j,n}, \quad \mu_n = \lambda - \lambda_n,
\label{S3:1dimeq}
\end{equation}
where $g_{j,n} \in L^{2,(1+ \epsilon_0)/2}((-\infty,\infty))$, and $v_{j,n}(y) = g_{j,n}(y) = 0$ for $y < 0$. 
Let $r_0(z) = (- \partial_y^2 - z)^{-1}$ in $L^2({\R})$, i.e.
$$
\left(r_0(z)g\right)(y) = \frac{i}{2\sqrt{z}}\int_{-\infty}^{\infty}e^{i\sqrt{z}|y-y'|}
g(y')dy',
$$
where ${\rm Im}\,\sqrt{z} \geq 0$.
Then as can be checked easily for any $s > 0$ and $\delta > 0$, there exists a constant $C_{s,\delta} > 0$ such that
$$
\|(1 + |y|)^{s}r_0(- a)(1 + |y|)^{-s}\|_{{\bf B}(L^2({\R});L^2({\R}))} \leq C_{s,\delta}, \quad \forall a > \delta.
$$
Therefore by (\ref{S3:1dimeq}), one can show that
\begin{equation}
E_j((\lambda,\infty))v_j \in L^{2,(1+ \epsilon_0)/2}(\Omega_j),
\label{S3:estimates>lambda}
\end{equation}
where $E_j(\cdot)$ is the spectral projection associated with $- \Delta_{h_j}$.

For $\lambda_n < \lambda$, we study $v_{j,n}$ separately. 
By (\ref{eq:decayinOmegaj}),
\begin{equation}
\lim_{R\to\infty}\frac{1}{R}\int_{0}^{R}|v_{j,n}(y)|^2dy = 0.
\label{S3:1dimdecay}
\end{equation}
 In view of (\ref{S3:1dimdecay}), we see that $v_{j,n}$ satisfies both of the outgoing and incoming radiation conditions. Adopting the outgoing radiation condition, we see that $v_{j,n}$ is written as $v_{j,n} = r_0(\mu_n + i0)g_{j,n}$, i.e.
$$
v_{j,n}(y) = \frac{i}{2\sqrt{\mu_n}}\Big(\int_0^{y}e^{i\sqrt{\mu_n}(y-y')}g_{j,n}(y')dy' + \int_y^{\infty}e^{i\sqrt{\mu_n}(y'-y)}g_{j,n}(y')dy'\Big).
$$
Note that $g_{j,n} \in L^1((0,\infty))$, since $g_{j,n} \in L^{2,(1 + \epsilon_0)/2}((0,\infty))$.
Therefore
$$
\lim_{y\to\infty}v_{j,n}(y) = \frac{i}{2\sqrt{\mu_n}}\int_0^{\infty}e^{i\sqrt{\mu_n}(y-y')}g_{j,n}(y')dy'.
$$
The condition (\ref{S3:1dimdecay}) implies that this limit is equal to 0.
which implies
$$
v_{j,n}(y) = \frac{i}{2\sqrt{\mu_n}}\Big(- \int_y^{\infty}e^{i\sqrt{\mu_n}(y -y')}g_{j,n}(y')dy' + \int_y^{\infty}e^{i\sqrt{\mu_n}(y'-y)}g_{j,n}(y')dy'\Big).
$$
Using the following Lemma 3.7 (Hardy's inequality), we have $(1 + y)^{(\epsilon_0 - 1)/2}v_{j,n} \in L^2((0,\infty))$. Using (\ref{S3:estimates>lambda}), we the have $v_j \in L^{2,(-1 + 2\epsilon_0)/2}((0,\infty))$. By Lemma 3.2 and the formula $\chi_ju = v_j + \mathcal E_j^{(0)}B_j^{(0)}u$,  we have $u \in L^{2,(-1+\epsilon_0)/2}(\Omega_j)$. Thus we have seen that $u$ gains the decay of order $\epsilon_0$ in $\Omega_j$. Then in (\ref{S3:1dimeq}), $g_{j,n} \in L^{2,(1 + 2\epsilon_0)/2}((0,\infty)$. Hence $v_{j} \in L^{2,(-1+ 3\epsilon/2})(\Omega_j)$.  Repeating this procedure, we obtain $\chi_ju \in L^{2,m\epsilon_0}(\Omega_j), \ \forall m > 0$. The estimate (\ref{S3chijuUnifEstimate}) can be proven by re-examining the above arguments.
\qed


\begin{lemma}
Let $f(y) \in L^1((0,\infty))$ and put
$$
u(y) = \int_y^{\infty}f(t)dt.
$$
The for $s > 1/2$
$$
\int_0^{\infty}y^{2(s-1)}|u(y)|^2dy \leq \frac{4}{(2s-1)^2}\int_0^{\infty}
y^{2s}|f(y)|^2dy.
$$
\end{lemma}
For the proof, see \cite{IsKu08} Chap. 3, Lemma 3.3. \qed


\begin{lemma}
Let $\sigma_{rad}(H)$ be the set of $\lambda \not\in \mathcal T(H)$ for which there exists a non-trivial solution $u \in \mathcal B^{\ast}$ of the equation $(H - \lambda)u = 0$ satisfying the radiation condition. Then $\sigma_{rad}(H) = \sigma_p(H)\setminus\mathcal T(H)$. Moreover it is a discrete subset of ${\R}\setminus\mathcal T(H)$ with possible accumulation points in $\mathcal T(H)$ and $\infty$.
\end{lemma}
\noindent {\bf Proof.} The first part of the lemma is proved by Lemmas 3.5 and 3.6. Let $I$ be a compact interval in ${\R}\setminus\mathcal T(H)$, and suppose there exists an infinite number of eigenvalues (counting multiplicities) in $I$. Let $u_n, n = 1, 2, \cdots$ be the associated orthonormal eigenvectors.

Take any $R > 0$ and let $\chi_0 = \chi_0^R$ be the function introduced in the beginning of this section. We decompose
$$
u_n = \chi_0^Ru_n + \sum_{j=0}^N(1 - \chi_0^R)\chi_ju_n.
$$
 Then by (\ref{S3chijuUnifEstimate}), for any $\epsilon > 0$ there exists $R > 0$ such that $\|(1 - \chi_0^R)u_n\|_{L^2(\Omega)} < \epsilon$ uniformly in $n$. By the compactness of the imbedding of $H^1_{loc}(\Omega)$ to $L^2_{loc}(\Omega)$, $\{\chi_Ru_n\}_n$ is compact in $L^2(\Omega)$. Therefore $\{u_n\}_n$ contains a convergent subsequence, which is a contradiction. \qed

\bigskip
{\it The set of exceptional points} $\mathcal E(H)$ is now defined by
\begin{equation}
\mathcal E(H) = \mathcal T(H)\bigcup\sigma_p(H).
\label{S3:ExcepPoints}
\end{equation}
{\newtekst Weyl's formula for the asymptotic distribution of eigenvalues on compact manifolds and Lemma 3.8 imply that $\mathcal T(H)$ is discrete and $\mathcal E(H)$ has only finite number of accumulation points on any compact interval in $\R$. }


\subsection{Limiting absorption principle} 
For a self-adjoint operator $H$ defined in a Hilbert space $\mathcal H$, 
the limit 
$$
\lim_{\epsilon\to 0}(H - \lambda \mp i\epsilon)^{-1} =: (H - \lambda \mp i0)^{-1}, \quad \lambda \in \sigma(H),
$$
does not exist as a bounded operator on $\mathcal H$. However if $\lambda$ is in the continuous spectrum of $H$, it is sometimes possible to guarantee the existence of the above limit in ${\bf B}(\mathcal X;\mathcal X^{\ast})$, where $\mathcal X, \mathcal X^{\ast}$ are Banach spaces such that $\mathcal X \subset \mathcal H = H^{\ast} \subset \mathcal X^{\ast}$, and $\mathcal H$ is identified with its dual space via Riesz' theorem. This fact, called the limiting absorption principle, is central in studying the absolutely continuous spectrum, and many works are devoted to it. We employ in this paper the classical method of integration by parts pioneered by Eidus \cite{Eid65}. The crucial step is to establish a-priori estimates as in \S 2 of this paper, and to show the uniqueness of solutions to the reduced wave equation satisfying the radiation condition. After this hard analysis part, the remaining arguments are almost routine. 
 
We take any compact interval $I \subset (0,\infty)\setminus\mathcal E(H)$ 
and let
\begin{equation}
J = \{z \in {\C} \, ; \, {\rm Re}\,z \in I, \ {\rm Im}\,z \neq 0\}.
\nonumber
\end{equation}
We first note that Lemma 2.3 also holds for the solution to the equation
\begin{equation}
\left\{
\begin{split}
&(H - z)u = f \quad {\rm in} \quad \Omega, \\
& \partial_{\nu}u = 0 \quad {\rm on} \quad \partial\Omega,
\end{split}
\right.
\nonumber
\end{equation}
by the standard elliptic regularity estimates.
We put $u = R(z)f$ and $v_j = \chi_ju - \mathcal E_j^{(0)}B_j^{(0)}u$ as in the proof of Lemma 3.6, Then $v_j$ satisfies
(\ref{Lemma36chijueq}) with $\lambda$ replaced by $z$.
We can then apply Theorem 2.7 to see that
\begin{equation}
\|\chi_ju\|_{\mathcal B^{\ast}} \leq C_s\Big(\|f\|_{\mathcal B} + 
\|u\|_{-s}\Big),
\label{eq:Perturbestimate}
\end{equation}
for any $s > 1/2$, where $C$ is independent of $z \in J$. Once (\ref{eq:Perturbestimate}) is proved, we can repeat the arguments in Chap 2, \S 2 of \cite{IsKu08} or those of Ikebe-Saito \cite{IkSa72} without any essential change. Note that here and in the sequel, we use $(\;,\,)$ to denote the inner product
\begin{equation}
(u,v) = \int_{\Omega}u\overline{v} dV
\nonumber
\end{equation}
of $L^2(\Omega)$ as well as the coupling between ${\mathcal B}$ and $\mathcal B^{\ast}$, or $L^{2,s}$ and $L^{2,-s}$.


\begin{lemma}
Take $s > 1/2$ sufficiently close to $1/2$. \\
\noindent
(1) There exists a constant $C > 0$ such that
$$
\sup_{z\in J}\|R(z)f\|_{-s} \leq C\|f\|_{\mathcal B}.
$$
(2) For any $\lambda \in I$ and $f \in \mathcal B$, the strong limit $\lim_{\epsilon\to0}R(\lambda \pm i\epsilon)f$exists in $L^{2,-s}$. \\
\noindent 
(3) $I \ni \lambda \to R(\lambda \pm i0)f \in L^{2,-s}$ is continuous.
\end{lemma}

\noindent
{\it Sketch of the proof}. 
Suppose the uniform bound (1) is not true. Then there exist a sequence $z_n \in J$ and $f_n \in \mathcal B$ such that $z_n = R(z_n)f_n$ satisfies $\|u_n\|_{-s} = 1$ and $\|f_n\|_{\mathcal B} \to 0$. Without loss of generality, we can assume that $z_n \to \lambda \in I$. Using (\ref{eq:Perturbestimate}) with $0 < s' < s$ and the compactness of the embedding of $H^2_{loc}$ into $L^2_{loc}$, one can assume that $u_n$ converges to some $u \in \mathcal B^{\ast}$, and $u$ satisfies the equation $(H - \lambda)u = 0$ and the radiation condition (see Corollary 2.6). Therefore $u = 0$ by Lemma 3.8. However this contradicts $\|u_n\|_{-s} = 1$. The assertions (2) and (3) are proved in a similar manner. \qed

\medskip
Using this lemma one can prove the following theorem.


\begin{theorem}
(1) For any $\lambda \in I$, $\lim_{\epsilon \to 0} R(\lambda \pm i\epsilon)f$ exists in the weak-${\ast}$ sense:
\begin{equation}
\exists \lim_{\epsilon\to0}(R(\lambda\pm i\epsilon)f,g) =: (R(\lambda \pm i0)f,g), \quad \forall f, g \in \mathcal B.
\nonumber
\end{equation}
(2) There exists a constant $C > 0$ such that
\begin{equation}
\|R(\lambda \pm i0)f\|_{\mathcal B^{\ast}} \leq C\|f\|_{\mathcal B}, \quad 
\lambda \in I.
\nonumber
\end{equation}
Moreover $R(\lambda \pm i0)f$ satisfies the outgoing  radiation condition for $\lambda + i0$ and incoming radiation condition for $\lambda - i0$.\\
\noindent
(3) For any $f, g \in \mathcal B$, $I \ni \lambda \to (R(\lambda \pm i0)f,g)$ is continuous. \\
\noindent
(4) Let $E(\cdot)$ be the spectral decomposition of $H$. Then $E([0,\infty)\setminus\mathcal E(H))L^2(\Omega)$ $ = \mathcal H_{ac}(H)$, and we have the following orthogonal decomposition
$$
L^2(\Omega) = \mathcal H_{ac}(H) \oplus \mathcal H_p(H).
\nonumber
$$
\end{theorem}

{\it Sketch of the proof}. Since $L^{2,-s} \ (s > 1/2)$ is dense in ${\mathcal B}^{\ast}$, (1) follow from Lemma 3.9 (2) and (\ref{eq:Perturbestimate}). The assertion (2) follows from Lemma 3.9 (1) and (\ref{eq:Perturbestimate}). The remaining assertions are proved in the same way as in Chap 2, \S 2 of \cite{IsKu08} or  Ikebe-Saito \cite{IkSa72}. \qed

\medskip
Let us recall that for a self-adjoint operator $H = \int_{-\infty}^{\infty}\lambda dE(\lambda)$, the absolutely continuous subspace for $H$, $\mathcal H_{ac}(H)$, is the set of $u$ such that $(E(\lambda)u,u)$ is absolutely continuous with respect to $d\lambda$, and the point spectral subspace, $\mathcal H_p(H)$, is the closure of the linear hull of eigenvectors of $H$.


\section{Forward problem}


\subsection{Unperturbed spectral representations} 
Let $\{\chi_j\}_{j=0}^N$ be the partition of unity defined in \S 3.
Recall the spaces $\mathcal B$ and $\mathcal B^{\ast}$ introduced in \S 2.
 For two functions $f, g$ on $\Omega$, $f \simeq g$ means that 
 \begin{equation}
\lim_{R \to \infty}\frac{1}{R}\int_{0}^R\|\chi_j(y)\left(f(\cdot,y) - g(\cdot,y)\right)\|^2_{L^2(M_j)}dy = 0, \quad 1 \leq \forall j \leq N.
\nonumber
\end{equation}
We also use the same notation $f \simeq g$ for $f, g$ defined on $\Omega_j$.

Green's function of $- d^2/dy^2 - \zeta$ on $(0,\infty)$ with Neumann boundary condition at $y = 0$ is
\begin{equation}
G(y,y';\zeta) = \frac{i}{\sqrt{\zeta}}\left\{
\begin{split}
\cos(\sqrt{\zeta}y)e^{i\sqrt{\zeta}y'}, \quad 0 < y < y', \\
e^{i\sqrt{\zeta}y}\cos(\sqrt{\zeta}y'), \quad 0 < y' < y.
\end{split}
\right.
\nonumber
\end{equation}
Let $\lambda_{j,1} < \lambda_{j,2} \leq \cdots$ be the eigenvalues of $- \Delta_{h_j}$ with normalized eigenvectors $\varphi_{j,n}(x)$, $n = 1,2,\cdots$. Without loss of generality, we assume that $\varphi_{j,n}(x)$'s are real-valued. Let $H_j^{(0)} = - \partial_y^2 - \Delta_{h_j}$ with Neumann boundary condition. Then $R_j^{(0)}(z) = (H_j^{(0)} - z)^{-1}$ is written as
\begin{equation}
\left(R_j^{(0)}(z)f\right)(x,y) = \sum_{n=1}^{\infty}\int_0^{\infty}G(y,y';z - \lambda_{j,n})\left(P_{j,n}f\right)(x,y')dy',
\label{eq:R0ExplicitForm}
\end{equation}
\begin{equation}
\left(P_{j,n}f\right)(x,y) = \langle f(\cdot,y),\varphi_{j,n}\rangle\varphi_{j,n}(x),
\nonumber
\end{equation}
$\langle\; , \,\rangle$ being the inner product of $L^2(M_j\,;\sqrt{\det\,(h_j)}\,dx)$. Note that $\det\,(h_{ij}) = \det \, G_j^{(0)}$. For $f(x,y) \in L^2\big(M_j\times (0,\infty);\big(\det\,G_j^{(0)}\big)^{1/2}dxdy\big)$, we define its cosine transform by
\begin{equation}
\mathcal F_{cos}(\lambda)f(x) = \pi^{-1/2}\lambda^{-1/4}
\int_0^{\infty}\cos\Big(y\sqrt{\lambda}\Big)f(x,y)dy.
\nonumber
\end{equation}


\begin{lemma}
For $f \in {\mathcal B}$, and  $\lambda \in (0,\infty)\setminus\sigma_p(- \Delta_{h_j})$, we have
\begin{equation}
R_j^{(0)}(\lambda \pm i0)f \simeq \pm i \sqrt{\pi}
\sum_{\lambda_{j,n}<\lambda}
(\lambda - \lambda_{j,n})^{-1/4}e^{\pm iy\sqrt{\lambda - \lambda_{j,n}}}
\mathcal F_{cos}(\lambda - \lambda_{j,n})P_{j,n}f(x).
\nonumber
\end{equation}
\end{lemma}
\noindent {\bf Proof.} We first show that the right-hand side of (\ref{eq:R0ExplicitForm}) is a bounded operator from $\mathcal B$ to $\mathcal B^{\ast}$. The sum over the terms in which $\lambda_{j,n} > \lambda$ is rewritten as
\begin{equation}
\begin{split}
& A_j(\lambda)f \\
& :=  \sum_{\lambda_{j,n} > \lambda}\frac{1}{2k_n}\int_{0}^{\infty}
\left(e^{-k_n\,|y-y'|} 
+ e^{-k_n\,(y + y')}\right)f_{j,n}(x,y')dy'
\end{split}
\nonumber
\end{equation}
where $f_{j,n} = P_{j,n}f$ and $k_n = \sqrt{\lambda_{j,n}-\lambda}$. Then we have
\begin{equation}
\begin{split}
& \|A_j(\lambda)f(\cdot,y)\|^2_{L^2(M)} \\
& =  \sum_{\lambda_{j,n}> \lambda}
\frac{1}{4k_n^2}\left|\int_{0}^{\infty}
\left(e^{-k_n\,|y-y'|} 
+ e^{-k_n\,(y + y')}\right)
\langle f(\cdot,y'),\varphi_{j,n}\rangle_{L^2(M_j)}dy'\right|^2 \\
& \leq C_{\lambda}\sum_{\lambda_{j,n}>\lambda}\int_{0}^{\infty}
|\langle f(\cdot,y),\varphi_{j,n}\rangle\big|^2dy \\
& \leq C_{\lambda}\|f\|_{L^2(M_j\times(0,\infty))}^2.
\end{split}
\nonumber
\end{equation}
Hence $A_j(\lambda) \in {\bf B}(L^2;L^{\infty}(\R_+;L^2(M_j)) \subset {\bf B}(\mathcal B;\mathcal B^{\ast})$. To estimate the term in which $\lambda_{j,n} < \lambda$, we put
$$
u_{j,n}(x) = \int_0^{\infty}G(y,y';\lambda \pm i0 - \lambda_{j,n})f_{j,n}(x,y')dy'.
$$ 
Then we have
$$
|u_{j,n}(x)| \leq C_{\lambda}\int_0^{\infty}|f_{j,n}(x,y)|dy.
$$
Since
$$
\|u_{j,n}\|_{{\mathcal B}^{\ast}} \leq C\|u_{j,n}\|_{L^{\infty}}, \quad
\|f_{j,n}\|_{L^1} \leq C\|f_{j,n}\|_{\mathcal B},
$$
We have proven that $R_j^{(0)}(\lambda \pm i0) \in {\bf B}({\mathcal B};{\mathcal B}^{\ast})$.

Now the assertion of the lemma is easy to prove if there exists $n_0 > 0$ such 
that $f_{j,n} = 0$ for $n \geq n_0$, and $f_{j,n}$ is compactly supported for $n < n_0$. Since such an $f$ is dense in $\mathcal B$, we have proven the lemma.
 \qed

\bigskip
The generalized eigenfunction of $H_j^{(0)}$ is defined for $\lambda > \lambda_{j,n}$
\begin{equation}
\Psi_{j,n}^{(0)}(x,y;\lambda) = \pi^{-1/2}(\lambda - \lambda_{j,n})^{-1/4}
\cos\Big(y\sqrt{\lambda - \lambda_{j,n}}\Big)\varphi_{j,n}(x).
\label{eq:Psim0}
\end{equation}
This $\Psi_{j,n}^{(0)}(x,y;\lambda)$ is often denoted by $\Psi_{j,n}^{(0)}(\lambda)$ in the sequel. 
It satisfies
\begin{equation}
\left\{
\begin{split}
& (- \Delta_{G^{(0)}_j} - \lambda)\Psi_{j,n}^{(0)}(\lambda) = 0 \quad {\rm in} \quad \Omega_j, \\
&\partial_{\nu_j^{(0)}}\Psi_{j,n}^{(0)}(\lambda) = 0 \quad {\rm on} \quad \partial\Omega_j.
\end{split}
\right.
\label{eq:Psim0equation}
\end{equation}
The Fourier transformation associated with $H_j^{(0)}$ is defined by
\begin{equation}
\mathcal F^{(0)}_j(\lambda)f = \sum_{n=1}^{\infty}c_{\lambda_{j,n}}(\lambda)\mathcal F_{j,n}^{(0)}(\lambda)f,
\label{eq:Ujdefine}
\end{equation}
where $c_{j,n}(\lambda)$ is the characteristic function of the interval $(\lambda_{j,n},\infty)$, and
\begin{equation}
\begin{split}
\left(\mathcal F_{j,n}^{(0)}(\lambda)f\right)(x) & = \left(\int_{\Omega_j}
\overline{\Psi_{j,n}^{(0)}(\lambda)}fdV_j^{(0)}\right)\varphi_{j,n}(x) \\
 &= \mathcal F_{cos}(\lambda - \lambda_{j,n})P_{j,n}f(x),
\end{split}
\label{eq:DefineUmj}
\end{equation}
where $dV_j^{(0)} = \big(\det\,G_j^{(0)}\big)^{1/2}dxdy$. 
Define a subspace of $L^2((0,\infty);L^2(M_j);d\lambda)$ by
\begin{equation}
\begin{split}
\widehat{\mathcal H}_j &= \sum_{n=1}^{\infty}L^2\big((\lambda_{j,n},\infty);d\lambda\big)\otimes\varphi_{j,n}(x) \\
& = \left\{\sum_{n=1}^{\infty}f_n(\lambda)\varphi_{j,n}(x) \, ; \, f_n \in L^2\big((\lambda_{j,n},\infty)\,;\,d\lambda\big)\right\}.
\end{split}
\label{S4:mathcalH}
\end{equation}
Then $\mathcal F^{(0)}_j$ defined by $\big(\mathcal F^{(0)}_jf)(\lambda) = \mathcal F^{(0)}_j(\lambda)f$ for $f \in C_0^{\infty}(\Omega_j)$ is uniquely extended to a unitary operator
\begin{equation}
\mathcal F_j^{(0)} : L^2(\Omega_j) \to \widehat{\mathcal H}_j.
\nonumber
\end{equation}
We put
\begin{equation}
{\bf h} = {\mathop\bigoplus_{j=1}^{N} L^2(M_j)},
\label{eq:Defineboldh}
\end{equation}
where $L^2(M_j) = L^2(M_j;\sqrt{\det(h_j)}\,dx)$, and
and also
\begin{equation}
\mathcal F^{(0)} = (\mathcal F^{(0)}_1,\cdots,\mathcal F^{(0)}_N).
\label{defineU}
\end{equation}

By the computation similar to the one to be given in the proof of Lemma 4.3 below, one can show that
\begin{equation}
\frac{1}{2\pi i}\left(\left[R_j^{(0)}(\lambda + i0) - R_j^{(0)}(\lambda - i0)\right]f,f\right) = \|\mathcal F_j^{(0)}(\lambda)f\|_{L^2(M_j)}^2.
\nonumber
\end{equation}
Therefore, $\mathcal F^{(0)}_j(\lambda) \in {\bf B}(\mathcal B\,;L^2(M_j))$, and  $\mathcal F^{(0)}_j(\lambda)^{\ast} \in {\bf B}(L^2(M_j);\mathcal B^{\ast})$.

Here we must pay attention to the following remarks. The first one is that in (\ref{eq:Ujdefine}), $\mathcal F_j^{(0)}(\lambda)$ is a finite sum:
\begin{equation}
\mathcal F_j^{(0)}(\lambda) = \sum_{\lambda_{j,n}<\lambda}\mathcal F_{j,n}^{(0)}(\lambda).
\label{S4Fjlambdafinitesum}
\end{equation}
The second remark is that the adjoint $\ast$ is taken in the following sense:
\begin{equation}
(\mathcal F_j^{(0)}(\lambda)f,h)_{L^2(M_j)} = 
(f,\mathcal F_j^{(0)}(\lambda)^{\ast}h)_{L^2(\Omega_j)}
= \int_{\Omega_j}f\,\overline{\mathcal F_j^{(0)}(\lambda)^{\ast}h}\,dV_j^{(0)},
\label{astl2Omegaj}
\end{equation}
$(h \in L^2(M_j)$). 
Therefore 
\begin{equation}
\mathcal F_j^{(0)}(\lambda)^{\ast} = \sum_{\lambda_{j,n}<\lambda}\mathcal F_{j,n}^{(0)}(\lambda)^{\ast},
\label{S4Fj0ast}
\end{equation}
and for $h \in L^2(M_j)$
\begin{equation}
\left(\mathcal F_{j,n}^{(0)}(\lambda)^{\ast}h\right)(x,y) = 
\Psi_{j,n}^{(0)}(x,y;\lambda)(h,\varphi_{j,n})_{L^2(M_j)}.
\label{S4Fjn0ast}
\end{equation}
Since $\mathcal F^{(0)}_j(\lambda)^{\ast}$ satisfies $(H_j^{0)} - \lambda)\mathcal F^{(0)}_j(\lambda)^{\ast} = 0$, we have
\begin{equation*}
\mathcal F^{(0)}_j(\lambda)^{\ast} \in {\bf B}(L^2(M_j);H^{2,-s}), \quad s > 1/2,
\end{equation*}
hence
\begin{equation*}
\mathcal F^{(0)}_j(\lambda) \in {\bf B}(H^{-2,s};L^2(M_j)), \quad s > 1/2.
\end{equation*}


\subsection{Perturbed spectral representations}
Using $\mathcal E_j$, $B_j$ and $V_j(z)$ in Subsection 3.1, for $\lambda > \lambda_{j,n}$ we define the generalized eigenfunction for $H$ by
\begin{equation}
\Psi_{j,n,\pm}(\lambda) = \left(\chi_j - {\mathcal E}_jB_j\right)\Psi_{j,n}^{(0)}(\lambda) - 
R(\lambda \mp i0)V_j(\lambda)\Psi_{j,n}^{(0)}(\lambda).
\label{eq:Psimj}
\end{equation}
Here putting $s = (1 + \epsilon_0)/2$, we regard $\mathcal E_jB_j$ and $V_j(\lambda)$ in ${\bf B}(H^{2,-s};L^{2,s})$. 
Note that $\Psi_{j,n,\pm}(\lambda) \in {\mathcal B}^{\ast}$. This definition easily implies
\begin{equation}
\left\{
\begin{split}
&(- \Delta_{G} - \lambda)\Psi_{j,n,\pm}(\lambda) = 0 \quad {\rm in} \quad \Omega, \\
& \partial_{\nu}\Psi_{j,n,\pm}(\lambda) = 0 \quad {\rm on} \quad 
\partial\Omega.
\end{split}
\right.
\end{equation}

The generalized Fourier transformation for $H$ is defined by perturbing $\mathcal F_j^{(0)}$. We put for $\lambda > \lambda_{j,n}$
\begin{equation}
\mathcal F_{j,n,\pm}(\lambda) =  \mathcal F_{j,n}^{(0)}(\lambda)J_j\big(\chi_j - ({\mathcal E}_jB_j)^{\ast} - V_j(\lambda)^{\ast}\big)R(\lambda \pm i0),
\label{eq:Fouriermj}
\end{equation}
where $J_j = \big(\det\, G/\det\, G_j^{(0)}\big)^{1/2}$. Note that $({\mathcal E}_jB_j)^{\ast}, V_j(\lambda)^{\ast} \in {\bf B}(L^{2,-s};H^{-2,s})$, and $R(\lambda \pm i0) \in {\bf B}(L^{2,s};H^{2,-s})\cap{\bf B}(H^{-2,s};L^{2,-s})$, hence (\ref{eq:Fouriermj}) is well-defined.


\begin{lemma} 
For $f \in C_0^{\infty}(\Omega)$, we have
\begin{equation}
\big(\mathcal F_{j,n,\pm}(\lambda)f\big)(x)  = \left(\int_{\Omega}\overline{\Psi_{j,n,\pm}(\lambda)}f\,dV\right)\varphi_{j,n}(x), 
\end{equation}
where $dV = \left({\rm det}\,(G)\right)^{1/2}dxdy$.
\end{lemma}
Proof. We put $u = \big(\chi_j - (\mathcal E_jB_j)^{\ast} - V_j(\lambda)^{\ast}\big)R(\lambda \pm i0)f$. Then by using (\ref{astl2Omegaj})
\begin{equation}
\begin{split}
(\mathcal F_{j,n,\pm}(\lambda)f,h)_{L^2(M_j)} &= 
(\mathcal F_{j,n}^{(0)}(\lambda)J_ju,h)_{L^2(M_j)} \\
&= \int_{\Omega_j}u\,\overline{\mathcal F_{j,n}^{(0)}(\lambda)^{\ast}h}\,J_jdV_j^{(0)}.
\nonumber
\end{split}
\end{equation}
We then use (\ref{S4Fjn0ast}) to see that the right-hand side is equal to 
\begin{equation}
\begin{split}
& \int_{\Omega_j}u\,\overline{\Psi_{j,n}^{(0)}(\lambda)}\,dV\,\overline{(h,\varphi_{j,n})}_{L^2(M_j)} \\
&=  
\left(\big(\chi_j - (\mathcal E_jB_j)^{\ast} - V_j(\lambda)^{\ast}\big)R(\lambda \pm i0)f,\Psi_{j,n}^{(0)}(\lambda)\right)\,\overline{(h,\varphi_{j,n})}_{L^2(M_j)} \\
&= \left(f,\Psi_{j,n,\pm}(\lambda)\right)(\varphi_{j,n},h)_{L^2(M_j)},
\end{split}
\nonumber
\end{equation}
which proves the lemma. \qed

The adjoint operator $\mathcal F_{j,n,\pm}(\lambda)^{\ast}$ is defined by the following formula:
\begin{equation}
(\mathcal F_{j,n,\pm}(\lambda)f,h)_{L^2(M_j)} = (f,\mathcal F_{j,n,\pm}(\lambda)h^{\ast})_{L^2(\Omega)}, \quad h \in L^2(M_j).
\label{S4fjnipminnerproduct}
\end{equation}


\begin{lemma} The adjoint operator $\mathcal F_{j,n,\pm}(\lambda)^{\ast}$ has the following expression:
\begin{equation}
\mathcal F_{j,n,\pm}(\lambda)^{\ast} = \left(\chi_j - \mathcal E_jB_j - R(\lambda \mp i0)V_j(\lambda)\right)\mathcal F_{j,n}^{(0)}(\lambda)^{\ast},
\label{S4Fjnpmast}
\end{equation}
where the adjoint $\mathcal F_{j,n}^{(0)}(\lambda)^{\ast}$ is taken in the sense of (\ref{astl2Omegaj}).
\end{lemma}
Proof. Let $u = \big(\chi_j - (\mathcal E_jB_j)^{\ast} - V_j(\lambda)^{\ast}\big)R(\lambda \pm i0)f$. Then as is shown in the proof of Lemma 4.2, 
\begin{equation}
\begin{split}
(\mathcal F_{j,n,\pm}(\lambda)f,h)_{L^2(M_j)} 
&= \int_{\Omega_j}u\,\overline{\mathcal F_{j,n}^{(0)}(\lambda)^{\ast}h}\,J_jdV_j^{(0)} \\
&= (u,\mathcal F_{j,n}^{(0)}(\lambda)^{\ast}h)_{L^2(\Omega)}.
\nonumber
\end{split}
\end{equation}
Plugging the form of $u$, we see that the right-hand side is equal to
$$
(f,(\chi_j - \mathcal E_jB_j - R(\lambda \mp i0)V_j(\lambda))\mathcal F_{j,n}^{(0)}(\lambda)^{\ast}h)_{L^2(\Omega)},
$$
which proves the lemma. \qed

We define
\begin{equation}
\mathcal F_{j,\pm}(\lambda) = \sum_{n=1}^{\infty}c_{j,n}(\lambda)\mathcal F_{j,n,\pm}(\lambda) = \sum_{\lambda_{j,n}<\lambda}\mathcal F_{j,n,\pm}(\lambda),
\label{eq:Fourierj}
\end{equation}
\begin{equation}
\mathcal F_{\pm}(\lambda) = (\mathcal F_{1,\pm}(\lambda),\cdots,
\mathcal F_{N,\pm}(\lambda)).
\label{eq:Fouriertotal}
\end{equation}


\begin{lemma}
For any $\lambda \in (0,\infty)\setminus\mathcal E(H)$ and $f \in \mathcal B$, we have on $\Omega_j$ \\
\begin{equation}
R(\lambda \pm i0)f \simeq \pm i \,\sqrt{\pi}\sum_{\lambda_{j,n}<\lambda}(\lambda - \lambda_{j,n})^{-1/4}e^{\pm iy\sqrt{\lambda-\lambda_{j,n}}}\mathcal F_{j,n,\pm}(\lambda)f.
\label{eq:LimitResolv}
\end{equation}
\end{lemma}
\noindent {\bf Proof.} This follows from (\ref{eq:Resolventeq}),  Lemma 4.1 and the definition (\ref{eq:Fouriermj}). \qed


\begin{lemma}
For any $\lambda \in (0,\infty)\setminus\mathcal E(H)$ and $f \in \mathcal B$, we have 
\begin{equation}
\frac{1}{2\pi i}\left((R(\lambda + i0) - R(\lambda - i0))f,f\right) = 
\|\mathcal F_{\pm}(\lambda)f\|^2_{\bf h}.
\nonumber
\end{equation}
\end{lemma}
\noindent {\bf Proof.} We prove the case for $\mathcal F_{+}(\lambda)$. We have only to prove the lemma when $f \in C_0^{\infty}(\Omega)$. We compute in a way similar to that in Lemma 3.5. Take $\rho(t) \in C_0^{\infty}((0,\infty))$ such that $\int_0^{\infty}\rho(t)dt = 1$, and put $\chi(t) = \int_t^{\infty}\rho(s)ds$. Let $u = R(\lambda + i0)f$ and 
$$
\psi_R = \chi_0 + \sum_{j=1}^N\chi\big(\frac{y}{R}\big)\chi_j(y),
$$
where $\{\chi_j\}_{j=0}^N$ is the partition of unity on $\Omega$, and $y$ in $\chi_j(y)$ is the local coordinate on $\Omega_j$.
We then have
$$
\left([H - \lambda,\psi_R]u,u\right) = (\psi_Ru,f) - (f,\psi_Ru).
$$
As $u \in \mathcal B^{\ast}$, by computing the commutator $[H,\psi_R]$, we then have
$$
\lim_{R\to\infty}\sum_{j=1}^N\frac{2}{R}\left(\rho\big(\frac{y}{R}\big)\chi_j(y)\partial_yu,u\right) = (u,f) - (f,u).
$$
Since $u = R(\lambda + i0)f$ satisfies the radiation condition (see Theorem 3.10 (2)), $(\partial_y - iP_j(\lambda))\chi_ju \simeq 0$.  Therefore
$$
\lim_{R\to\infty}\sum_{j=1}^N\frac{2i}{R}\left(\rho\big(\frac{y}{R}\big)\chi_j(y)P_j(\lambda)u,u\right) = (u,f) - (f,u).
$$
Now we note that
\begin{equation}
\lim_{R\to\infty}\frac{1}{R}\left(\rho\big(\frac{y}{R}\big)\chi_j(y)P_j(\lambda)u,u\right) =
\lim_{R\to\infty}\frac{1}{R}\int_0^{\infty}\rho\big(\frac{y}{R}\big)\left(P_j(\lambda)u,u\right)_{L^2(M_j)}dy.
\nonumber
\end{equation}
Let $v_{\pm}$ be the term in the right-hand side of (\ref{eq:LimitResolv}).
Using Lemma 4.4, we first replace $u$ of the right-hand side of $\left(P_j(\lambda)u,u\right)_{L^2(M_j)}$ by $v_{\pm}$. We next move $P_j(\lambda)$ to the right-hand side of the inner product, and replace $u$ by $v_{\pm}$. Since $P_{j,n}(\lambda)\varphi_{j,n}^{(0)} = \sqrt{\lambda-\lambda_{j,n}}\,\varphi_{j,n}^{(0)}$, we have $P_j(\lambda)\mathcal F_{j,n,+}(\lambda) = \sqrt{\lambda - \lambda_{j,n}}\,\mathcal F_{j,n,+}(\lambda)$.
The lemma then follows from a direct computation. \qed

\medskip
The formula in Lemma 4.5, when integrated with respect to $\lambda$ over $(0,\infty)$, is a counterpart of the Parseval formula in the Fourier transformation, and a crucial step for the spectral representation.
Using $\widehat{\mathcal H}_j$ in (\ref{S4:mathcalH}), we put 
\begin{equation}
\widehat{\mathcal H} = {\mathop\bigoplus_{j=1}^N} \, \widehat{\mathcal H}_j. 
\label{S4WidehatH}
\end{equation}
The following theorem can be proved in the same way as in \cite{Ik75} or Chap. 3 of \cite{IsKu08}.
 

\begin{theorem}
(1) For $\lambda \not\in \mathcal T(H)$, $\mathcal F_{\pm}(\lambda) \in {\bf B}(\mathcal B;{\bf h})$. \\
\noindent
(2) The operator $(\mathcal F_{\pm}f)(\lambda) = \mathcal F_{\pm}(\lambda)f$ defined for $f \in {\mathcal B}$ is uniquely extended to a partial isometry with initial set  $\mathcal H_{ac}(H)$ and final set $\widehat{\mathcal H}$. \\
\noindent
(3) $\ (\mathcal F_{\pm}Hf)(\lambda) = \lambda\, (\mathcal F_{\pm}f)(\lambda), \quad \forall \lambda \in (0,\infty)\setminus{\mathcal E}(H), \quad \forall f \in D(H)$. \\
\noindent
(4) $\ \mathcal F_{\pm}(\lambda)^{\ast} \in {\bf B}({\bf h};{\mathcal B}^{\ast})$ is an eigenoperator of $H$ with eigenvalue $\lambda$ in the sense that
\begin{equation}
(H - \lambda)\mathcal F_{\pm}(\lambda)^{\ast}\psi = 0, \quad \forall \psi \in {\bf h}.
\nonumber 
\end{equation}
(5)  For any compact interval $I \subset (0,\infty)\setminus{\mathcal T}(H)$  and $g \in \widehat{\mathcal H}$, we have
$$
\int_I \mathcal F_{\pm}(\lambda)^{\ast}g(\lambda)\,d\lambda \in L^2(\Omega).
$$
Let $I_n$ be a finite union of compact intervals in $(0,\infty)\setminus{\mathcal E}(H)$ such that $I_n \subset I_{n+1}$, $\cup_{n=1}^{\infty}I_n = (0,\infty)\setminus\mathcal E(H)$. Then for any $f \in \mathcal H_{ac}(H)$, the inversion formula holds:
\begin{equation}
f = {\mathop{\rm s-lim}_{n\to\infty}}\int_{I_n}\mathcal F_{\pm}(\lambda)^{\ast}(\mathcal F_{\pm}f)(\lambda)d\lambda.
\nonumber
\end{equation}
\end{theorem}


\subsection{Time-dependent scattering theory} Let $H_{j}^{(0)} = - \partial_y^2 - \Delta_{h_j}$ be the unperturbed Laplacian in the end $\Omega_j$.  


\begin{theorem} The wave operator $W_{\pm} : {\mathop\bigoplus_{j=1}^N}L^2(\Omega_j) \to L^2(\Omega)$ defined by
\begin{equation}
W_{\pm} = {\mathop{\rm s-lim}_{t\to\pm\infty}}\sum_{j=1}^Ne^{it\sqrt{H}}\chi_je^{-it\sqrt{H_{j}^{(0)}}} 
\nonumber
\end{equation}
exists and is complete, i.e. {\rm Ran}$\,W_{\pm} = \mathcal H_{ac}(H)$. Moreover
\begin{equation}
W_{\pm} = \big(\mathcal F_{\pm}\big)^{\ast}\mathcal F^{(0)},
\label{eq:WaveOpStationary}
\end{equation}
where $\mathcal F^{(0)}$ is the Fourier transformation  defined by (\ref{defineU}) for the system of Laplacians $(H_1^{(0)},\cdots,H^{(0)}_N)$.
\end{theorem}

{\it Sketch of the proof}. We argue in the same way as in Chap. 2, Theorem 8.9 of \cite{IsKu08}.
Take $f \in \mathcal H_{ac}(H)$ such that $\big(\mathcal F_{j,n,+}f\big)(\lambda) \in C_0^{\infty}((\lambda_{j,n},\infty))$ and $\mathcal F_{j,n,+}f = 0$ except for a finite number of $n$. Then by Theorem 4.6 and Lemma 4.3
\begin{equation}
\begin{split}
e^{-it\sqrt{H}}f =& \int_0^{\infty}e^{-it\sqrt{\lambda}}\mathcal F_{+}(\lambda)^{\ast}\big(\mathcal F_{+}f\big)(\lambda)d\lambda \\
=& \sum_{j,n}\int_0^{\infty}e^{-it\sqrt{\lambda}}
\chi_j\left(\mathcal F_{j,n}^{(0)}(\lambda)\right)^{\ast}
\big(\mathcal F_{j,n,+}f\big)(\lambda)d\lambda \\
& -  \sum_{j,m}\int_0^{\infty}e^{-it\sqrt{\lambda}}
\mathcal E_jB_j\left(\mathcal F_{j,n}^{(0)}(\lambda)\right)^{\ast}
\big(\mathcal F_{j,n,+}f\big)(\lambda)d\lambda \\
&- \sum_{j,n}\int_0^{\infty}e^{-it\sqrt{\lambda}}R(\lambda - i0)V_j(\lambda)
\left(\mathcal F_{j,n}^{(0)}(\lambda)\right)^{\ast}
\big(\mathcal F_{j,n,+}f\big)(\lambda)d\lambda.
\end{split}
\label{eq:eitHspectral}
\end{equation}
Because of the decay of $\mathcal E_j$, the 2nd term of the right-hand side tends to $0$ in $L^2(\Omega)$.
Letting $A = \sqrt{H}$, we have
$$
(H - k^2 + i0)^{-1} = (A - k + i0)^{-1}(A + k)^{-1}.
$$
We then put
$$
g(k) = 2k(A + k)^{-1}V_j(k^2)\left(\mathcal F_{j,n}^{(0)}(k^2)\right)^{\ast}\big(\mathcal F_{j,n,+}f\big)(k^2),
$$
$$
\widetilde g(t) = \int_0^{\infty}e^{-itk}g(k)dk.
$$
We show that
\begin{equation}
\|\widetilde g(t)\| \leq C(1 + t)^{-1-\epsilon}, \quad t > 0.
\label{eq:widetildegtintegrable}
\end{equation}
In fact, take $h \in L^2(\Omega)$ and consider
\begin{equation}
\begin{split}
(\widetilde g(t),h) &= \int_0^{\infty}2ke^{-itk}\big(\left(\mathcal F_{j,n}^{(0)}(k^2)\right)^{\ast}\big(\mathcal F_{j,n,+}f\big)(k^2),V_j(k^2)(A+k)^{-1}h\big)dk \\
&= \sum_n\int dV_j^{(0)}\int_0^{\infty}\Big(e^{-i(tk + y\sqrt{k^2-\lambda_{j,n}})} 
+ e^{-i(tk - y\sqrt{k^2 - \lambda_{j,n}})}\Big)\cdots dk.
\end{split}
\nonumber
\end{equation}
Since $V_j(k^2)$ contains a factor $(1 + y)^{-1-\epsilon}$, splitting the integral suitably and integrating by parts, we have
$$
|(\widetilde g(t),h)| \leq C(1 + t)^{-1-\epsilon}\|h\|,
$$
which proves (\ref{eq:widetildegtintegrable}). We use the notation $f(t) \sim g(t)$ if $\|f(t) - g(t)\| \to 0$ as $t \to \infty$. In view of the following Lemma 4.8, we obtain as $t \to \infty$
\begin{equation}
\begin{split}
e^{-it\sqrt{H}}f \sim & \sum_{j,n}\chi_j\int_0^{\infty}e^{-it\sqrt{\lambda}}
\left(\mathcal F_{j,n}^{(0)}(\lambda)\right)^{\ast}
\big(\mathcal F_{j,n,+}f\big)(\lambda)d\lambda \\
=& \sum_{j,n} \chi_je^{-it\sqrt{H_j^{(0)}}}\big(\mathcal F_{j,n}^{(0)}\big)^{\ast}
\mathcal F_{j,n,+}f,
\end{split}
\nonumber
\end{equation}
in $L^2(\Omega)$.  
This implies the existence of the limit
\begin{equation}
{\mathop{\rm s-lim}_{t\to\infty}}\sum_{j=1}^Ne^{it\sqrt{H_j^{(0)}}}\chi_je^{-it\sqrt{H}} = \left(\mathcal F^{(0)}\right)^{\ast}\mathcal F_+.
\label{eq:InverseWaveOp}
\end{equation}
 Since $\left(\mathcal F^{(0)}\right)^{\ast}\mathcal F_+$ is a partial isometry with initial set $\mathcal H_{ac}(H)$ and final set $L^2(\Omega)$, (\ref{eq:InverseWaveOp}) also implies for $g = (g_1,\cdots,g_N) \in {\mathop\bigoplus_{j=1}^N}L^2(\Omega_j)$
\begin{equation}
\|e^{it\sqrt{H}}\sum_{j=1}^N\chi_je^{-it\sqrt{H_j^{(0)}}}g - \big(\mathcal F_+\big)^{\ast}\mathcal F^{(0)}g\| \to 0.
\label{Waveop}
\end{equation}

Let us prove this fact. We put $U(t) = \sum_{j=1}^Ne^{it\sqrt{H_j^{(0)}}}\chi_je^{-it\sqrt{H}}$. Then (\ref{eq:InverseWaveOp}) implies that $U(t) \to \big(\mathcal F^0\big)^{\ast}\mathcal F_+ =: U$ strongly, which implies 
\begin{equation}
U(t)^{\ast} \to U^{\ast} \quad {\rm weakly}.
\label{S4Weakconv}
\end{equation}
 We show that
\begin{equation}
\|U(t)^{\ast}g\| \to \|g\| = \|U^{\ast}g\|, \quad g = (g_1,\cdots,g_N) \in {\mathop\bigoplus_{j=1}^N}L^2(\Omega_j).
\label{normconv}
\end{equation}
In fact, we have
$$
\|U(t)^{\ast}g\|^2 = \|\sum_{j=1}^N\chi_je^{-it\sqrt{H_j^{(0)}}}g_jf\|^2 = 
\sum_{j=1}^N\|\chi_je^{-it\sqrt{H_j^{(0)}}}g_j\|^2.
$$
By the scattering property of $e^{-it\sqrt{H_j^{(0)}}}$, $\|(1- \chi_j)e^{-it\sqrt{H_j^{(0)}}}g_j\| \to 0$, which proves
$$
\sum_{j=1}^N\|\chi_je^{-it\sqrt{H_j^{(0)}}}g_j\|^2 \to \sum_{j=1}^N\|g_j\|^2 = \|g\|^2.
$$
Now, (\ref{S4Weakconv}) and (\ref{normconv}) yield $\|U(t)^{\ast}g - U^{\ast}g\| \to 0$.
This completes the proof of Theorem 4.7 for $W_+$. The assertion for $W_-$ is proved similarly. \qed


\begin{lemma}
Let $A$ be a self-adjoint operator on a Hilbert space $\mathcal H$. For $f(k) \in C_0((0,\infty);\mathcal H)$, we put
$$
\widetilde f_{\pm}(t) = \int_0^{\infty}e^{\pm ik t}f(k)dk.
$$
Then for any $\epsilon > 0$
$$
\left\|\int_0^{\infty}(A - k \mp i\epsilon)^{-1}e^{\pm ikt}f(k)dk\right\| 
\leq \int_t^{\infty}\|\widetilde f_{\pm}(s)\|ds.
$$
\end{lemma}

\noindent {\bf Proof.} This is proved in \cite{IsKu08}, Chap. 2, Lemma 8.10. For the reader's convenience, we reproduce the proof. By virtue of the identity
$$
(A - k \mp i\epsilon)^{-1} = \pm i\int_0^{\infty}
e^{\mp is(A - k \mp i\epsilon)}ds,
$$
we have
$$
\int_0^{\infty}(A - k \mp i\epsilon)^{-1}e^{\pm ikt}f(k)dk = 
\pm i \int_0^{\infty}e^{\mp is(A \mp i\epsilon)}
\widetilde f_{\pm}(s + t)ds,
$$
which proves the lemma. \qed


\subsection{S-matrix}
The scattering operator is defined by $S = \big(W_+\big)^{\ast}W_-$. 
We consider its Fourier transform : $\widehat S = \mathcal F^{(0)}S\left(\mathcal F^{(0)}\right)^{\ast}$.


\begin{lemma}
We have a direct integral representation:
$$
(\widehat Sf)(\lambda) = \widehat S(\lambda)f(\lambda), \quad 
\forall f \in \widehat{\mathcal H}, \quad \forall \lambda > 0,
$$
where $\widehat S(\lambda) = \big(\widehat S_{jk}(\lambda)\big)_{1\leq j,k \leq N}$ is a bounded operator on $\bf h$ called the S-matrix, and is written as follows
$$
\widehat S_{jk}(\lambda) = \delta_{jk} - 2\pi i {\mathcal F}_{j,+}(\lambda)V_k(\lambda)\left(\mathcal F_{k}^{(0)}(\lambda)\right)^{\ast}.
$$
\end{lemma}
\noindent {\bf Proof.} Lemma 4.5 implies
$$
\frac{1}{2\pi i}\left(R(\lambda + i0) - R(\lambda - i0)\right) = 
\mathcal F_{\pm}(\lambda)^{\ast}\mathcal F_{\pm}(\lambda).
$$
By Lemma 4.3, we then have
\begin{equation}
\mathcal F_{k,+}(\lambda)^{\ast} - \mathcal F_{k,-}(\lambda)^{\ast} 
= 2\pi i \mathcal F_+(\lambda)^{\ast}\mathcal F_+(\lambda)V_k(\lambda)
\mathcal F_{k}^{(0)}(\lambda)^{\ast}.
\nonumber
\end{equation}
Then we have by Theorem 4.6 (2), for $f, g \in \widehat{\mathcal H}$
\begin{equation}
\big((\mathcal F_+ - \mathcal F_-)(\mathcal F_+)^{\ast}f,g) = 
 - 2\pi i\sum_{k=1}^N\int_0^{\infty}\left(f(\lambda),\mathcal F_+(\lambda)V_k(\lambda)\left(\mathcal F_{k}^{(0)}(\lambda)\right)^{\ast}g(\lambda)\right)_{\bf h}d\lambda. 
\nonumber
\end{equation}
By (\ref{eq:WaveOpStationary}), $\widehat S = \mathcal F_+\big(\mathcal F_-\big)^{\ast}$. Hence the lemma follows. \qed

\medskip
Let ${\bf h}_j(\lambda)$ be the linear subspace of $L^2(M_j)$ spanned by $\varphi_{j,n}$ such that $\lambda_{j,n} < \lambda$ and put
\begin{equation}
{\bf h}(\lambda) = \bigoplus_{j=1}^N {\bf h}_j(\lambda).
\nonumber
\end{equation}
Then $\widehat S(\lambda)$ is a partial isometry on $\bf h$ with initial and final set ${\bf h}(\lambda)$.
The scattering amplitude is defined by
\begin{equation}
A_{\,jk}(\lambda) = \mathcal F_{j,+}(\lambda)V_k(\lambda)\left(\mathcal F_{k}^{(0)}(\lambda)\right)^{\ast}.
\nonumber
\end{equation}
{\newtekst Let $A_{jm,kn}(\lambda) : L^2(M_k) \to L^2(M_j)$ be given by}
\begin{equation}
A_{jm,kn}(\lambda) = \mathcal F_{j,m,+}(\lambda)V_k(\lambda)
\left(\mathcal F_{k,n}^{(0)}(\lambda)\right)^{\ast}.
\label{eq:Ajmkn}
\end{equation}
We then have
\begin{equation}
\widehat S_{jk}(\lambda) - \delta_{jk}I_j = - 2\pi i
\sum_{\lambda_{j,m}<\lambda,\, \lambda_{k,n}<\lambda}
A_{jm,kn}(\lambda),
\nonumber
\end{equation}
{where $I_j$ is the identity operator on $L^2(M_j)$.}
When $j, k$ and the energy $\lambda > 0$ is fixed, $\big(A_{jm,kn}(\lambda)\big)$ is a finite matrix of size $(d_j,d_k)$, where $d_j = \, \#\{(j,m) ; \lambda_{j,m}< \lambda\}$. 
Let $\mathcal A_{jm,kn}(\lambda)$ be defined by
\begin{equation}
\mathcal A_{jm,kn}(\lambda) = \left(A_{jm,kn}(\lambda)\varphi_{k,n},\varphi_{j,m}\right)_{L^2(M_j)}.
\label{S4MathcalAjmkn}
\end{equation}
Then we have
\begin{equation}
A_{jm,kn}(\lambda)h = \mathcal A_{jm,kn}(\lambda)(h,\varphi_{k,n})_{L^2(M_k)}\varphi_{j,m}, \quad 
\forall h \in L^2(M_k).
\end{equation}

The scattering amplitude is computed from the asymptotic expansion of the generalized eigenfunction in the following way.

\begin{lemma}
\begin{equation}
 P_{j,m}\left(\Psi_{k,n,-}(\lambda) - \chi_j\Psi_{k,n}^{(0)}(\lambda)\right) 
 \simeq - \frac{i\sqrt{\pi}e^{iy\sqrt{\lambda - \lambda_{j,m}}}}{(\lambda - \lambda_{j,m})^{1/4}}
\mathcal A_{jm,kn}(\lambda)\varphi_{j,m}.
\nonumber
\end{equation} 
\end{lemma}
\noindent {\bf Proof.} This directly follows from (\ref{eq:Psimj}) and Lemma 4.4.  \qed


\section{{\newtekst From scattering data to boundary data}}


\subsection{Non-physical scattering amplitude}
In this section, we observe waves coming in from and going out of the end $\Omega_1$ assuming that 
\begin{equation}
G_1 = (dy)^2 + h_1(x,dx) \quad {\rm on} \quad \Omega_1.\label{eq:Omega1Gfree}
\end{equation}
This amounts to studying the scattering amplitude $A_{1m,1n}(\lambda)$ of (\ref{eq:Ajmkn}), which is rewritten as
\begin{equation}
\begin{split}
A_{1m,1n}(\lambda) & = 
\mathcal F_{cos}(\lambda - \lambda_{1,m})P_{1,m}J_1\big(\chi_1  - V_1(\lambda)^{\ast}R(\lambda + i0)\big) \\
& \ \ \ \ \cdot V_1(\lambda)\left(\mathcal F_{cos}(\lambda - \lambda_{1,n})\right)^{\ast}P_{1,n}.
\end{split}
\label{eq:Smnlambda}
\end{equation}
Note that $B_1 = 0$, because of the assumption (\ref{eq:Omega1Gfree}).
By the expression (\ref{eq:Vjzformula}), $V_1(\lambda)$ and $V_1(\lambda)^{\ast}$ are {independent of $\lambda$ and compactly supported in the $y$-variable.}
Therefore, $A_{1m,1n}(\lambda)$ defined for $\lambda > {\rm max}\,\{\lambda_{1,m},\lambda_{1,n}\}$ is analytically continued to the upper half plane ${\C}_+ = \{{\rm Im}\,\lambda > 0\}$. {\newtext This analytic continuation
can be extended to a continuous function on 
$\C_+\cup (\R\setminus {\mathcal E}(H))$. We denote the obtained function
for $\{\lambda < \max\{\lambda_{1,m},\lambda_{1,n}\}\}\setminus {\mathcal E}(H)$
by $A_{1m,1n}^{(nph)}(\lambda)$ and call it the} {\it non-physical} scattering amplitude. 
{These functions can be} represented by (\ref{eq:Smnlambda}), 
where $\mathcal F_{cos}(\lambda - \lambda_{1,m})$ 
and $\mathcal F_{cos}(\lambda - \lambda_{1,n}) $ are replaced by their analytic continuations. Let
\begin{equation}
\Phi_{1,n}^{(0)}(x,y;\lambda) = 
\pi^{-1/2}e^{-\pi i/4}(\lambda_{1,n} - \lambda)^{-1/4}
\cosh\Big(y\sqrt{\lambda_{1,n} - \lambda}\Big)\varphi_{1,n}(x),
\label{eq:Nonphyfreeeigenfunc}
\end{equation}
and put, {similarly to (\ref{eq:DefineUmj})}
\begin{equation}
\left(\mathcal F_{cosh}(\lambda_{1,n} - \lambda)P_{1,n}f\right)(x) = 
\left(\int_{\Omega_1}\overline{\Phi_{1,n}^{(0)}(\lambda)}f\,dV_1^{(0)}\right)\varphi_{1,n}(x).
\nonumber
\end{equation}
In the following, we always assume that $\lambda \not\in \mathcal E(H)$.
The explicit form of $ A_{1m,1n}^{(nph)}(\lambda)$ is  
{\newtext given by the following lemma. Recall that the non-physical scattering amplitude
${A}_{1m,1n}^{(nph)}(\lambda)$ coincides with the physical scattering
amplitude ${ A_{1m,1n}}(\lambda)$ for $\lambda > {\rm max}\,\{\lambda_{1,m},\lambda_{1,n}\}$.} 

\begin{lemma}
(1) If $\lambda_{1,m} < \lambda < \lambda_{1,n}$,
\begin{equation}
\begin{split}
A_{1m,1n}^{(nph)}(\lambda) & = 
\mathcal F_{cos}(\lambda - \lambda_{1,m})P_{1,n}J_1\big(\chi_1  - V_1(\lambda)^{\ast}R(\lambda + i0)\big) \\
& \ \ \ \ \cdot V_1(\lambda)\left(\mathcal F_{cosh}(\lambda_{1,n}-\lambda)\right)^{\ast}P_{1,n}.
\end{split}
\nonumber
\end{equation}
(2) If $\lambda_{1,n} < \lambda < \lambda_{1,m}$,
\begin{equation}
\begin{split}
A_{1m,1n}^{(nph)}(\lambda) & = 
\mathcal F_{cosh}(\lambda_{1,m}-\lambda)P_{1,m}J_1\big(\chi_1  - V_1(\lambda)^{\ast}R(\lambda + i0)\big) \\
& \ \ \ \ \cdot V_1(\lambda)\left(\mathcal F_{cos}(\lambda - \lambda_{1,n})\right)^{\ast}P_{1,n}.
\end{split}
\nonumber
\end{equation}
(3) If $\lambda < {\rm min}\{\lambda_{1,m}, \lambda_{1,n}\}$,
\begin{equation}
\begin{split}
A_{1m,1n}^{(nph)}(\lambda) & = 
\mathcal F_{cosh}(\lambda_{1,m}-\lambda)P_{1,m}J_1\big(\chi_1  - V_1(\lambda)^{\ast}R(\lambda + i0)\big) \\
& \ \ \ \ \cdot V_1(\lambda)\left(\mathcal F_{cosh}(\lambda_{1,n}-\lambda)\right)^{\ast}P_{1,n}.
\end{split}
\nonumber
\end{equation}
\end{lemma}

In accordance with (\ref{eq:Psimj}), we define non-physical eigenfunction by
\begin{equation}
\Phi_{1,m,\pm}(\lambda) = \chi_1\Phi_{1,m}^{(0)}(\lambda) - R(\lambda \mp i0)V_1(\lambda)\Phi_{1,m}^{(0)}(\lambda).
\label{e:Nonphyseigenfn}
\end{equation}
Note that the {\it physical} eigenfunction $\Psi_{1,m,-}(\lambda)$ defined for $\lambda > \lambda_{1,m}$ is 
analytically continued {\newtext through the upper half space $\C_+$} to the {\it nonphysical} eigenfunction $\Phi_{1,m,-}(\lambda)$ defined for $\lambda < \lambda_{1,m}$. The non-physical scattering amplitude is computed from the asymptotic behavior of non-physical eigenfunction in the following way.

We put 
\ba
{\mathcal A}_{1m,1n}^{(nph)}(\lambda)=(A_{1m,1n}^{(nph)}(\lambda)\varphi_{1,n},\varphi_{1,m})_{L^2(M_1)}.
\ea
Then we have for $h\in L^2(M_1)$
\ba
A_{1m,1n}^{(nph)}(\lambda)h={\mathcal A}_{1m,1n}^{(nph)}(\lambda)\,(h,\varphi_{1,n})_{L^2(M_1)}\,\varphi_{1,m}.
\ea


\begin{lemma}
(1) If $\lambda_{1,m} < \lambda < \lambda_{1,n}$, we have as $y \to \infty$,
\begin{equation}
P_{1,m}\Big(\Phi_{1,n,-}(\lambda) - \Phi_{1,n}^{(0)}(\lambda)\Big)
\simeq -  \frac{i\sqrt{\pi}\, e^{iy\sqrt{\lambda - \lambda_{1,m}}}}{(\lambda - \lambda_{1,m})^{1/4}}
{\mathcal A}_{1m,1n}^{(nph)}(\lambda)\varphi_{1,n}
\nonumber
\end{equation}
(2) If $\lambda < \max\{\lambda_{1,m}, \lambda_{1,n}\}$, we have as $y \to \infty$,
\begin{equation}
P_{1,m}\Big(\Phi_{1,n,-}(\lambda) - \Phi_{1,n}^{(0)}(\lambda)\Big)
\sim - \frac{e^{\pi i/4}\, \sqrt{\pi}e^{-y\sqrt{\lambda_{1,m} - \lambda}}}{(\lambda_{1,m} - \lambda)^{1/4}}   {\mathcal A}_{1m,1n}^{(nph)}(\lambda)\varphi_{1,n},
\nonumber
\end{equation} 
with a super exponentially decreasing error, that
is, with the error  {\newtext $r(y)$ satisfying $|r(y)|\leq C_Ne^{-Ny}$ 
for any $N>0$}.
\end{lemma}

\noindent {\bf Proof.} The assertion (1) is proved in the same way as in Lemma 4.8. By (\ref{eq:R0ExplicitForm}), letting $\zeta = \lambda - \lambda_{1,m}$, we have as $y \to \infty$
\begin{equation}
P_{1,m}R_1(\lambda + i0)f(x,y) \sim
\frac{ie^{i\sqrt{\zeta}y}}{\sqrt{\zeta}}\int_0^{\infty}\cos\left(\sqrt{\zeta}y'\right)P_{1,m}f(x,y')dy'
\nonumber
\end{equation}
with a super exponentially decaying error. This, together with (\ref{eq:Resolventeq}) and Lemma 5.1, proves (2).
\qed

\subsection{Splitting the manifold}
 We take a compact hypersurface $\Gamma \subset \Omega_1$ having the following property. 
\medskip

\noindent
(C-1) $\Gamma$ splits $\Omega$ into a  union: $\Omega = \Omega_{ext}\cup\Omega_{int}$ 
so that $\Omega_{ext}\cap\Omega_{int} =\Gamma$,  $\Omega_{int}$ is a manifold
with smooth boundary,
 and  $\Omega_{ext}\subset \Omega_1$. (See figure 2.)
\medskip

\begin{figure}[htbp]
\begin{center}
\psfrag{1}{$\Omega_{int}$}
\psfrag{2}{$\Omega_{ext}$}
\psfrag{3}{$\Gamma$}
\psfrag{4}{$\Omega_1$}
\includegraphics[width=8cm]{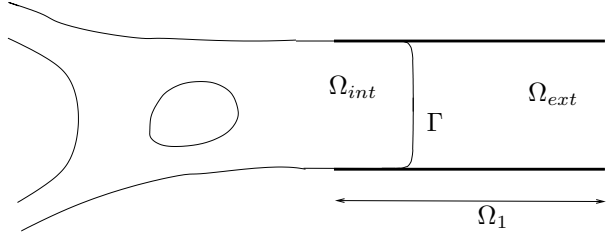}
 \label{pic 2}
\end{center}
\caption{Surface $\Gamma$ splits $\Omega$ to two parts, manifold
$\Omega_{int}$ with a smooth boundary and its complement
$\Omega_{ext}\subset \Omega_1$.
}
\end{figure}

Let ${\mathcal O}\subset  \Omega_{int}$ be an open, relatively compact
set such that it has a smooth boundary not intersecting $\p \Omega_{int}$ 
and that $\Omega_{int}\setminus {\mathcal O}$ is connected. 
Denote $\Omega_{\mathcal O}=\Omega_{int}\setminus {\mathcal O}$
and 
\begin{equation}
\Gamma_{\mathcal O} = \left\{
\begin{split}
& \Gamma \quad {\rm if} \quad \mathcal O = \emptyset, \\
& \partial\mathcal O \quad {\rm if} \quad \mathcal O \neq \emptyset.
\end{split}
\right.
\nonumber
\end{equation}

\medskip

We put for 
$f, g \in L^2(\Gamma_{\mathcal O})$
\begin{equation}
(f,g)_{\Gamma_{\mathcal O}} = \int_{\Gamma_{\mathcal O}}f(x)\overline{g(x)}dS_x,
\nonumber
\end{equation}
$dS_x$ being the measure induced from the metric $G$ on $\Gamma_{\mathcal O}$. We put $H_{\mathcal O} = - \Delta_G$ in $\Omega_{\mathcal O}$ endowed with the Neumann boundary condition:
\begin{equation}
\partial_{\nu}v = 0 \quad {\rm on} \quad \partial\Omega_{\mathcal O},
\label{eq:OmegainBC}
\end{equation}
$\nu$ being the unit normal to the boundary. 
If $\Omega$ has only one end, $\Omega_{int}$ is a bounded region. If $\Omega$ has more than one end, $\Omega_{int}$ is unbounded and the spectral theory developed for $H$ applies also to $H_{\mathcal O}$. To see this, we have only to replace $\mathcal K$ by $\mathcal K\cup\left(\Omega_1\cap\Omega_{int}\right)$, and to argue in the same way as in \S 3 and \S 4. Let $\mathcal E({\mathcal H}_{\mathcal O})$ be $\sigma_p(H_{\mathcal O})$ when $\Omega_{\mathcal O}$ is bounded, and the set of exceptional points when $\Omega_{\mathcal O}$ is unbounded.

Next we consider the case  ${\mathcal O}=\emptyset$ so that $\Gamma_{\mathcal O}=\Gamma$.


\begin{lemma}
Suppose $\lambda \not\in {\mathcal E}(H)\cup{\mathcal E}(H_{\emptyset})$, and let $\Psi_{1,n,-}(\lambda)$ and $\Phi_{1,n,-}(\lambda)$ be physical and non-physical eigenfunctions for $H$. Then the linear subspace spanned by $\partial_{\nu}\Psi_{1,n,-}(\lambda)\big|_{\Gamma}$, $\partial_{\nu}\Phi_{1,n,-}(\lambda)\big|_{\Gamma}$, $n = 1, 2, \dots$, is dense in $L^2(\Gamma)$.
\end{lemma} 
 
 \noindent {\bf Proof.} We show  that, if $f \in L^2(\Gamma)$ satisfies 
\begin{equation}
(f,\partial_{\nu}\Psi_{1,n,-}(\lambda))_{\Gamma} = (f,\partial_{\nu}\Phi_{1,n,-}(\lambda))_{\Gamma} = 0, \quad \forall n \geq 1,
\label{eq:vanishassump}
\end{equation}
then  $f = 0$. We define an operator $\delta_{\Gamma}'\in {\bf B}((H^{1/2}(\Gamma))' ; H^{-2}(\Omega))$,  
where $(H^{1/2}(\Gamma))' $ is the dual space of $H^{1/2}(\Gamma)$, by
$$
(\delta_{\Gamma}'f,w)=  (f,\partial_{\nu}w)_{\Gamma}, \quad \forall w \in H^{2}(\Omega),
$$
and put $u = R(\lambda - i0)\delta_{\Gamma}'f$ by duality. This means that ,
{\newtekstt  if  $G_-(\lambda;  X, X')$ is Green's function, i.e.  the integral (Schwartz) kernel } of $R(\lambda - i0)$,
$$
u(X)=\left(R(\lambda - i0)\delta_{\Gamma}'f\right)(X) = 
\int_{\Gamma}\partial_{\nu'}G_-(\lambda;X,X')f(X')dS_{X'},
$$
where $\partial_{\nu'}$ means the conormal differentiation with respect to the variable $X'$.
Then $u \in \mathcal B^{\ast}$,
and by (\ref{eq:Resolventeq}), we have the following asymptotic expansion on $\Omega_1$
\begin{equation}
u \simeq \sum_{\lambda_{1,n}<\lambda}C_n(\lambda)e^{-iy\sqrt{\lambda - \lambda_{1,n}}}\left(f,\partial_{\nu}\Psi_{1,n,-}(\lambda)\right)_{\Gamma}\varphi_{1,n}(x).
\nonumber
\end{equation}  
In particular, if $\lambda_{1,n} < \lambda$
\begin{equation}
(u,\varphi_{1,n}) \simeq C_n(\lambda)e^{-iy\sqrt{\lambda - \lambda_{1,n}}}
\left(f,\partial_{\nu}\Psi_{1,n,-}(\lambda)\right)_{\Gamma},
\label{eq:uvarphim<lambda}
\end{equation}
$C_n(\lambda)$ being a constant.
In a similar way, we have for $\lambda_{1,n} > \lambda$ 
\begin{equation}
(u,\varphi_{1,n}) \sim C_n'(\lambda)e^{-y\sqrt{\lambda_{1,n} - \lambda}}
\left(f,\partial_{\nu}\Phi_{1,n,-}(\lambda)\right)_{\Gamma}
\label{eq:uvarphim>lambda}
\end{equation}
modulo a super exponentially decaying term.
Note that $u_n = \big(u,\varphi_{1,n}\big)$ satisfies the equation $(- \partial_y^2 + \lambda_{1,n} - \lambda)u_n = 0$ for $y > a$, $a$ being a sufficiently large constant. In view of the assumption of (\ref{eq:vanishassump}) and (\ref{eq:uvarphim<lambda}), (\ref{eq:uvarphim>lambda}), we then have $\big(u,\varphi_{1,n}\big) = 0$ for $y > a$, hence $u(x,y) = 0$ for $y > a$. The unique continuation theorem then implies $u = 0$ on $\Omega_{ext}$.
By the property of classical double layer potential, $\partial_{\nu}u$ is continuous across $\Gamma$,
{\newtekstt so that $\p_\nu u| _{\Gamma}=0$. }

 Next we show that $u = 0$ in $\Omega_{int}$. In the region $\Omega_{int}$, we have $(- \Delta_G- \lambda)u = 0$. If $\Omega_{int}$ is bounded, then $u = 0$ since $\lambda$ is not a Neumann eigenvalue. If $\Omega_{int}$ is not bounded, $u$ satisfies the incoming radiation condition, since so does $u$ in $\Omega$. Then $u = 0$ in $\Omega_{int}$ by Lemma 3.4. {As $u = R(\lambda - i0)\delta_{\Gamma}'f \in L^2_{loc}(\Omega)$, it follows from the above that $u = 0$ in $\Omega$. } Applying $H - \lambda$, we have $\delta_{\Gamma}'f = 0$ as a distribution, hence $f = 0$ on $\Gamma$. \qed

\subsection{Interior boundary value problem}
For $z \in {\C}\setminus\mathcal E(H_{\mathcal O})$, we consider the following boundary value problem
\begin{equation}
\left\{
\begin{split}
& (H_{\mathcal O} - z)u = 0 \quad {\rm in} \quad \Omega_{{\mathcal O}}, \\
& \partial_{\nu}u = 0 \, \quad {\rm on} \quad \partial\Omega_{{\mathcal O}}\setminus\Gamma_{\mathcal O}, \\
& \partial_{\nu}u = f \in H_0^{1/2}(\Gamma_{\mathcal O}) \quad {\rm on} \quad \Gamma_{\mathcal O}.
\end{split}
\right.
\label{eq:NeumanninOmegain}
\end{equation} 
The incoming radiation condition is also imposed, if $\Omega_{int}$ is unbounded and $z \in \R$. The Neumann-Dirichlet map (N-D map) is then defined by
\begin{equation}
\Lambda_{\mathcal O}(z)f = u\big|_{\Gamma_{\mathcal O}},
\label{eq: N-D map}
\end{equation}
where $u$ is the solution to (\ref{eq:NeumanninOmegain}). When ${\mathcal O}=\emptyset$, we use for 
the N-D map of the operator $H_{\emptyset}$ the notation $\Lambda_{\mathcal O}(z)=\Lambda(z)$.

Now we consider the operator theoretical meaning of the N-D map. Note that from now on ${\mathcal O}$
may be a non-empty set.
We put $\mathcal F = (\mathcal F_c,\mathcal F_p)$, where $\mathcal F_c$ is the generalized Fourier transform for $H_{\mathcal O}$ (which is absent when $\Omega_{int}$ is bounded) and $\mathcal F_p$  is defined by
\begin{equation}
\mathcal F_p : {\mathcal H}_{p}(H_{\mathcal O}) \ni u \mapsto ((u,\psi_1), (u,\psi_2),\cdots),
\nonumber
\end{equation}
and where $\mathcal H_{p}(H_{\mathcal O})$ is the point spectral subspace for $H_{\mathcal O}$ and
$\psi_i$ is the eigenfunction associated with the eigenvalue $\lambda_{i}$ of $H_{\mathcal O}$. There are two kinds of generalized Fourier transformation, $\mathcal F_+$ and $\mathcal F_-$. Both choices will do as $\mathcal F_c$.
Then $\mathcal F$ is a 
unitary
\begin{equation}
{\mathcal F} : L^2(\Omega_{int}) \to \widehat{\mathcal H} \oplus {\C}^d,
\label{S5:SpecRepre}
\end{equation}
where $d = {\dim}\,\mathcal H_p(H_{\mathcal O})$. If $d = \infty$, ${\C}^d = \ell^2$.
Moreover, we have 
\begin{equation}
(H_{\mathcal O} - z)^{-1} = \int_0^{\infty}\frac{\mathcal F_c(\lambda)^{\ast}\mathcal F_c(\lambda)}{\lambda - z}d\lambda + 
\sum_{i=1}^{d}\frac{P_i}{\lambda_i - z},
\label{eq:Sect5ResolventandFourier}
\end{equation} 
where {\newtext $P_i$ are the eigenprojections associated with  eigenvalues $\lambda_i$,
numbered counting multiplicities by $i=1,2,\dots,d$, and the right-hand side converges in the sense of strong limit in $L^2(\Omega_{\mathcal O})$. 

Let $r_{\Gamma_{\mathcal O}} \in {\bf B}(H^1(\Omega_{\mathcal O});H^{1/2}(\Gamma_{\mathcal O}))$ be the trace operator to $\Gamma_{\mathcal O}$, 
$$
r_{\Gamma_{\mathcal O}} : H^1(\Omega_{\mathcal O}) \ni f \to f\big|_{\Gamma_{\mathcal O}} \in H^{1/2}(\Gamma_{\mathcal O}). 
$$
We define $\delta_{\Gamma_{\mathcal O}} \in {\bf B}((H^{1/2}(\Gamma_{\mathcal O}))';(H^{1}(\Omega_{\mathcal O}))')$ as the adjoint of $r_{\Gamma_{\mathcal O}}$:
\begin{equation}
(\delta_{\Gamma_{\mathcal O}}f,w)_{L^2(\Omega_{int})} = (f,r_{\Gamma_{\mathcal O}}w)_{L^2(\Gamma_{\mathcal O})}, \quad 
f \in (H^{1/2}(\Gamma_{\mathcal O}))', \quad w \in H^1(\Omega_{\mathcal O}).
\nonumber
\end{equation}
With this in mind we write 
$$
r_{\Gamma_{\mathcal O}} = \delta_{\Gamma_{\mathcal O}}^{\ast}.
$$


\begin{lemma}\label{lemma 5.4}
For $z \not\in {\mathcal E}(H_{\mathcal O})$, the N-D map has the following representation
\begin{equation}
\begin{split}
\Lambda_{\mathcal O}(z) &= \delta_{\Gamma_{\mathcal O}}^{\ast}(H_{\mathcal O} - z)^{-1}\delta_{\Gamma_{\mathcal O}} \\
& = \int_0^{\infty}\frac{\delta_{\Gamma_{\mathcal O}}^{\ast}\mathcal F_c(\lambda)^{\ast}\mathcal F_c(\lambda)\delta_{\Gamma_{\mathcal O}}}{\lambda - z}d\lambda + \sum_{i=1}^{d}\frac{\delta_{\Gamma_{\mathcal O}}^{\ast}P_i\delta_{\Gamma_{\mathcal O}}}{\lambda_i - z}.
\end{split}
\nonumber
\end{equation}
\end{lemma}

 \noindent {\bf Proof.} For $f \in H_0^{1/2}(\Gamma_{\mathcal O})$, take $\widetilde f \in H^{2}(\Omega_{\mathcal O})$ such that $\partial_{\nu}\widetilde f = f$ on $\Gamma_{\mathcal O}$ and $\widetilde f$ has compact support in $\overline{\Omega_{int}}$. Then the solution $u$ of (\ref{eq:NeumanninOmegain}) is written as $u = \widetilde f - (H_{\mathcal O} - z)^{-1}(- \Delta_{G} - z)\widetilde f$. Let $g = \mathcal F_c(\lambda)(\Delta_{G} + z)\widetilde f$.
Then for any {$h \in \widehat{\mathcal H}_{int}$, where $\widehat{\mathcal H}_{int}$ is defined by (\ref{S4WidehatH}) with $j = 2,\dots,N$.}
\begin{equation}
\begin{split}
(\mathcal F_c(\lambda)(\Delta_{G}+z)\widetilde f,h) &= 
((\Delta_G + z)\widetilde f,\mathcal F_c(\lambda)^{\ast}h) \\
&= (\partial_{\nu}\widetilde f,r_{\Gamma_{\mathcal O}}\mathcal F_c(\lambda)^{\ast}h)_{L^2(\Gamma_{\mathcal O})} + (\widetilde f,(\Delta_G + \overline{z})\mathcal F_c(\lambda)^{\ast}h) \\
&= (f,r_{\Gamma_{\mathcal O}}\mathcal F_c(\lambda)^{\ast}h)_{L^2(\Gamma_{\mathcal O})} + (\widetilde f, (- \lambda + \overline{z})\mathcal F_c(\lambda)^{\ast}h).
\end{split}
\nonumber
\end{equation}
This implies
$$
\mathcal F_c(\lambda)(\Delta_G + z)\widetilde f = \mathcal F_c(\lambda)\delta_{\Gamma_{\mathcal O}}f + (- \lambda + z)\mathcal F_c(\lambda)\widetilde f,
$$
Hence
$$
\int_0^{\infty}\frac{\mathcal F_c(\lambda)^{\ast}\mathcal F_c(\lambda)(\Delta_G + z)\widetilde f}{\lambda - z}d\lambda = \int_0^{\infty}\frac{\mathcal F_c(\lambda)^{\ast}\mathcal F_c(\lambda)\delta_{\Gamma_{\mathcal O}}f}{\lambda - z}d\lambda - \mathcal F_c^{\ast}\mathcal F_c\widetilde f.
$$
Similarly,
$$
\sum_{i=1}^{d}\frac{P_i(\Delta_G + z)\widetilde f}{\lambda_i - z} = 
\sum_{i=1}^{d}\frac{P_i\delta_{\Gamma_{\mathcal O}} f}{\lambda_i - z} - \sum_{i=1}^{d} P_i\widetilde f.
$$
Since $\mathcal F_c^{\ast}\mathcal F_c\widetilde f+\sum_{i=1}^{d} P_i\widetilde f=
\widetilde f$,
by  (\ref{eq:Sect5ResolventandFourier}), these imply that
$$
u = (H_{\mathcal O} - z)^{-1}\delta_{\Gamma_{\mathcal O}}f,
$$
which proves the lemma. \qed

\bigskip
Let us call the set
\begin{equation}
\Big\{\delta_{\Gamma_{\mathcal O}}^{\ast}\mathcal F_c(\lambda)^{\ast}\mathcal F_c(\lambda)\delta_{\Gamma_{\mathcal O}} ; \lambda \in (0,\infty)\setminus \mathcal E(H_{\mathcal O})\Big\} \cup\Big\{\Big(\lambda_i, \delta_{\Gamma_{\mathcal O}}^{\ast}P_i\delta_{\Gamma_{\mathcal O}}\Big)\Big\}_{i=1}^{d},
\label{eq:Sec5BSP}
\end{equation}
where $d = \dim{\mathcal H}_p(H_{\mathcal O})$, the boundary spectral projection ({\bf BSP}) for $H_{\mathcal O}$
on $\Gamma_{\mathcal O}$. On the other hand, the set 
\begin{equation}
\Big\{\mathcal F_c(\lambda)\delta_{\Gamma_{\mathcal O}} ; \lambda \in (0,\infty)\setminus \mathcal E(H_{\mathcal O})\Big\} \cup\Big\{\Big(\lambda_i, \psi_i(x)\big|_{\Gamma_{\mathcal O}}\Big)\Big\}_{i=1}^{d}
\label{eq:SectBSD}
\end{equation}
is called the boundary spectral data ({\bf BSD}) on $\Gamma_{\mathcal O}$.

By using the formula (\ref{S3:HelffeSjostrand}), we have the following lemma.
 

\begin{lemma} For a bounded Borel function $\varphi(\lambda)$ with support in 
${\R}\setminus\mathcal T(H_{\mathcal O})$, {where $\mathcal T(H_{\mathcal O})$ is defined by (\ref{S3Thresholds}) with $j = 2, \cdots, N$}, we have
\begin{equation}
\delta_{\Gamma_{\mathcal O}}^{\ast}\varphi(H_{\mathcal O})\delta_{\Gamma_{\mathcal O}} = 
\int_0^{\infty}\varphi(\lambda)\delta_{\Gamma_{\mathcal O}}^{\ast}{\mathcal F_c}(\lambda)^{\ast}{\mathcal F_c}(\lambda)\delta_{\Gamma_{\mathcal O}}d\lambda + 
\sum_{i=1}^{d}\varphi(\lambda_i)\delta_{\Gamma_{\mathcal O}}^{\ast}P_i\delta_{\Gamma_{\mathcal O}}.
\nonumber
\end{equation}
\end{lemma}

 \noindent {\bf Proof.} By the formulae (\ref{S3:HelffeSjostrand}) and (\ref{eq:Sect5ResolventandFourier}), this lemma holds for any $\varphi(\lambda) \in C_0^{\infty}({\R}\setminus\mathcal T(H_{\mathcal O}))$. The general case the follows from the approximation. \qed

\medskip
Usually BSD is referred as given data in the inverse boundary value problems. What is actually used in {\newtext our reconstruction
for the manifold is the} BSP.


\begin{lemma}\label{lemma 5.6} Let $\O \subset \Omega_{int}$. Then knowing the N-D map $\Lambda_{\mathcal O}(z)$ for all $z \not\in \sigma(H_{\mathcal O})$ is equivalent to knowing the BSP for  $H_{\mathcal O}$.
\end{lemma}
\noindent {\bf Proof.} By Lemma \ref{lemma 5.4}, one can compute the N-D map by using BSP.
Taking $\varphi(\lambda)$ as the characteristic function of the interval $[a,t)$ and {taking note of the remark after (\ref{S3:ExcepPoints}), we differentiate} the formula in Lemma 5.5 with respect to $t$ to recover 
$\delta_{\Gamma_{\mathcal O}}^{\ast}\mathcal F_c(t)^{\ast}\mathcal F_c(t)\delta_{\Gamma_{\mathcal O}}$ for $t \in \R\setminus{\mathcal E}(H_{\mathcal O})$. Since
$$
\sum_{i=1}^{d}\frac{\delta_{\Gamma_{\mathcal O}}^{\ast}P_i\delta_{\Gamma_{\mathcal O}}}{\lambda_i - z} = 
\Lambda_{\mathcal O}(z) - \int_0^{\infty}\frac{\delta_{\Gamma_{\mathcal O}}^{\ast}\mathcal F_c(\lambda)^{\ast}\mathcal F_c(\lambda)\delta_{\Gamma_{\mathcal O}}}{\lambda - z}d\lambda,
$$
one can obtain eigenvalues $\lambda_i$ as the poles of the right-hand side. {The residues in these poles provide us with $\delta_{\Gamma_{\mathcal O}}^{\ast}\sum_{\lambda_j=\lambda_i}P_j\delta_{\Gamma_{\mathcal O}}$. This determines the terms $\delta_{\Gamma_{\mathcal O}}^{\ast}P_j\delta_{\Gamma_{\mathcal O}}$ 
for indexes $j$ such that $\lambda_j = \lambda_i$, up to an 
{orthogonal transformation of the eigenspace associated
to the eigenvalue $ \lambda_i$}, see \cite[Lem.\ 4.9]{KKL01} or \cite{KKL04b}. Thus we can determine the 
BSP for  $H_{\mathcal O}$.} \qed
\medskip

We complete this section by the following result used later to prove Theorem \ref{main thm}.}
Let $\Omega^{(r)}$, $r = 1,2$, be as in Theorem \ref{main thm} We take $\Gamma$ as above, which moreover has the following property: $G_1^{(1)} = G_1^{(2)}$ on
 $\Omega_{ext} = \Omega_{ext} ^{(1)} = \Omega_{ext} ^{(2)}$.
We put the superscript $^{(r)}$ for all relevant operators and functions explained above.
Let $\Lambda^{(r)}(\lambda)$, {$r=1,2$
 be the N-D map for $H_{\emptyset}^{(r)}$,  that is, when
 ${\mathcal O}=\emptyset$}. The basic idea of the following Lemma is due to Eidus
\cite{Eid65}.


\begin{lemma}\label{lemma 5.7} Under the assumptions of Theorem \ref{main thm},  we have $\Lambda^{(1)}(\lambda) = \Lambda^{(2)}(\lambda)$ for $\lambda \in (0,\infty)\setminus\cup_{r=1,2}(\mathcal E(H^{(r)})\cup\mathcal E({H_{\emptyset}^{(r)}}))$, {and 
BSP's for $H_{\emptyset}^{(1)}$ and $H_{\emptyset}^{(2)}$ coincide on $\Gamma$.}
\end{lemma}
\noindent \noindent {\bf Proof.} Since $\widehat S_{11}^{(1)}(\lambda) = \widehat S_{11}^{(2)}(\lambda)$, the physical scattering amplitudes coincide, hence so do non-physical scattering amplitudes by analytic continuation.
 Let $u = \Psi_{1,n,-}^{(1)}(\lambda) - \Psi_{1,n,-}^{(2)}(\lambda)$ and 
$v = \Phi_{1,n,-}^{(1)}(\lambda) - \Phi_{1,n,-}^{(2)}(\lambda)$. Then since $H^{(1)} = H^{(2)} = - \partial_y^2 - \Delta_{h_1}$ on $\Omega_{ext}$, {$u$ and $v$} satisfy $(- \partial_y^2 - \Delta_{h_1} - \lambda)u = 0$ and $(- \partial_y^2 - \Delta_{h_1} - \lambda)v = 0$ in $\Omega_{ext}$. Using Lemma 5.2 and arguing in the same way as in the proof of Lemma 5.3, we have $u = v = 0$ in $\Omega_{ext}$. Therefore, $\Psi_{1,n,-}^{(r)}$ and $\Phi_{1,n,-}^{(r)}$ as well as their  normal derivatives coincide for $r = 1, 2$ and for all $n \in {\bf Z}_+$. Since they satisfy the equation (\ref{eq:NeumanninOmegain}) for $H_{\emptyset} = H_{\emptyset}^{(r)}$, 
we have $\Lambda^{(1)}(\lambda) = \Lambda^{(2)}(\lambda)$ due to Lemma 5.3. {The last statement now follows immediately from Lemma \ref{lemma 5.6}. }\qed


\section{Boundary control method for manifolds with asymptotically cylindrical ends}
In this section we reconstruct the isometry type of the manifold  $(\Omega_{int},G)$
using given data.


\begin{figure}[htbp]
\begin{center}
\psfrag{1}{$\Omega_{ext}$}
\psfrag{2}{$\Sigma$}
\psfrag{3}{$X_1$}
\psfrag{4}{$\Omega_{int}$}
\psfrag{5}{$U_1$}
\includegraphics[width=12cm]{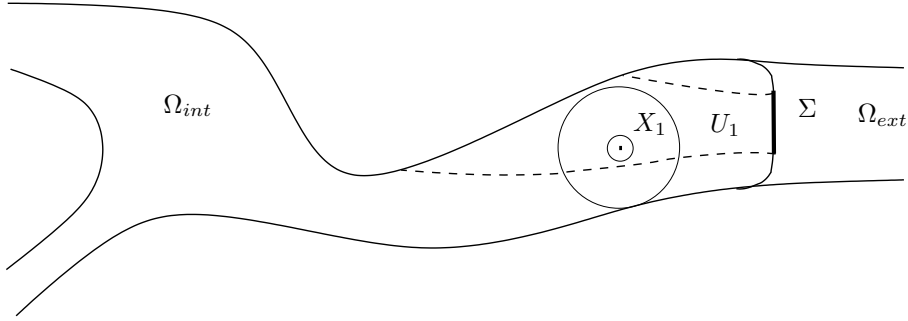} 
\label{pic 4}
\end{center}
\caption{We will construct manifold $\Omega_{int}$ by iterating
local constructions.
First, a neighborhood $U_1\subset \Omega_{int}$ of $\Sigma\subset \Gamma$ is reconstructed.
Next, a ball $\O=B(X_1,\rho)\subset U_1$ is removed
from the manifold and data analogous to measurements on $\p \O$ are
constructed. After that, the metric 
is reconstructed in a larger ball $B(X_1,\tau)$,
 and the procedure is iterated to reconstruct
the whole manifold $\Omega_{int}$.
}
\end{figure}

\begin{theorem} \label{BC thm} 
Assume that we are given
the set $\Gamma$ as a differentiable manifold,
the metric $G$ on $\Gamma$, and the BSP for $H_{\emptyset}$. These data determine the manifold $(\Omega_{int},G)$
up to an isometry.
\end{theorem}

For proving this theorem, we use  the boundary control method
described for compact manifolds in e.g. \cite{Be97,KKL01,KKL08}.
The reconstruction of non-compact manifolds is considered previously 
in the conference proceedings \cite{KKL04} and in \cite{BKLS} with different kind of data, 
using iterated time reversal for solutions of the wave equation.  
 
The proof of Theorem \ref{BC thm} is divided into a series of lemmas.
Our reconstruction of $(\Omega_{int},G)$ is of recurrent nature.
We will begin with the case when ${\mathcal O}=\emptyset$ so that
we are given just the set $\Gamma_{\mathcal O}=\Gamma$ as a differentiable manifold,
the metric on it, and the  BSP for the operator $H_{\emptyset}$ on $\Gamma$. We apply the boundary control method to reconstruct the metric $G$ on some neighborhood $U_1$ of $\Gamma$. Then, we will take a point 
$X_1\in U\setminus \Gamma$ and $\rho>0$
such that $B(X_1,2\rho)\subset U_1$,
where $B(X_1,r)$ denotes the ball of radius $r$ with center at $X_1$. 
We take ${\mathcal O}=B(X_1,\rho)$ and show that we can find
 the BSP for the operator $H_{\mathcal O}$ on $\Gamma_{\mathcal O}=\p \O$. Then we  apply the boundary control 
 method starting from $\Gamma_{\mathcal O}$,
 which would allow us to recover $(\Omega_{int},G)$ in a larger neighborhood $U_2 \supset U_1$ of $\Gamma$. Proceeding in this way, we will eventually recover the whole of $(\Omega_{int},G)$. Therefore, our further considerations deal with arbitrary $\mathcal O \subset \Omega_{int}$ including the case $\mathcal O = \emptyset$.

\subsection{Blagovestchenskii's identity}
Let us first consider the initial boundary value problem
\begin{equation}
\left\{
\begin{split}
& \partial_t^2u = \Delta_Gu, \quad {\rm in} \quad \Omega_{\mathcal O}\times \R_+, \\
& u\big|_{t=0} = \partial_tu\big|_{t=0} = 0, \quad {\rm in } \quad \Omega_{\mathcal O},\\
& \partial_{\nu}u = f, \quad {\rm in} \quad \partial\Omega_{\mathcal O}\times \R_+, \quad {\rm supp}\,f \subset \Gamma_{\mathcal O}\times{\R}_+.
\end{split}
\right.
\label{S5:IBVPWave}
\end{equation}


\begin{lemma} \label{lem: step 2}
Assume that we are given the set  $\Gamma_{\mathcal O}$ as a differentiable manifold, the metric $G$ on $\Gamma_{\mathcal O}$ and the BSP for $H_{\mathcal O}$ on $\Gamma_{\mathcal O}$. 
Then for any given $f,h\in C^\infty_0(\Gamma_{\mathcal O}\times \R_+)$ and $t,s>0$ these data uniquely determine
\ba
(u^f(t),u^h(s))=\int_{\Omega_{\mathcal O}} u^f(x,t)\,\overline {u^h(x,s)}\,dV_x.
\ea
\end{lemma} 

\noindent {\bf Proof.} 
Let
$$
S(t,\lambda) = \frac{\sin(\sqrt{\lambda}t)}{\sqrt\lambda}.
\nonumber
$$
Then the solution $u^f(t)$ is written as 
\begin{equation}
\begin{split}
u^f(t) = & \int_0^tds
\int_0^\infty d\lambda
S(t-s,\lambda)
\mathcal F_c(\lambda)^{\ast}\mathcal F_c(\lambda)\delta_{\Gamma_{\mathcal O}}f(s) \\
& + \int_0^tds\sum_{i=1}^{d}S(t-s,\lambda_i)P_i\delta_{\Gamma_{\mathcal O}}f(s).
\end{split}
\nonumber
\end{equation}
{Using the similar decomposition for $u^h(s)$ and the fact that $\mathcal F_c(\mu)\mathcal F_c(\lambda)^{\ast} = \delta(\mu-\lambda)$, we  obtain the following formula:}
\begin{equation}
\begin{split}
& (u^f(t),u^h(s)) \\
= & \int_0^tdt'\int_0^sds' \int_0^{\infty}d\lambda\,\widetilde S(t-t',s-s',\lambda)
\left(\delta_{\Gamma_{\mathcal O}}^{\ast}\mathcal F_c(\lambda)^{\ast}\mathcal F_c(\lambda)\delta_{\Gamma_{\mathcal O}}f(t'),h(s')\right)_{L^2(\Gamma_{\mathcal O})} \\
& + \int_0^tdt'\int_0^sds'\sum_{i=1}^{d} \widetilde S(t-t',s-s',\lambda_i)\left(\delta_{\Gamma_{\mathcal O}}^{\ast}P_i\delta_{\Gamma_{\mathcal O}}f(t'),h(s')\right)_{L^2(\Gamma_{\mathcal O})},
\end{split}
\label{eq:Sect5Blagovecont}
\end{equation}
where $\widetilde S(t,s,\lambda) = S(t,\lambda)S(s,\lambda)$. Observe that the right-hand side 
depends only on BSP and the metric on $\Gamma_{\mathcal O}$. 
\qed

\medskip
Above, the formula (\ref{eq:Sect5Blagovecont}) is a generalization of Blagovestchenskii identity (see 
\cite[Theorem 3.7]{KKL01}) for non-compact manifolds.

\subsection{Finite propagation property of waves}
Let us next introduce  some notations. 
For $t\geq 0$ and $\Sigma\subset \Gamma_{\mathcal O}$ arbitrary, let
\ba
\Omega_{\mathcal O}(\Sigma,t)=\{X \in \Omega_{\mathcal O}\ ; \ d_{\mathcal O}(X,\Sigma)\leq t\}
\ea
be the domain of influence of $\Sigma$ at time $t$. Here, 
$d_{\mathcal O}(X,Y)$ is the distance between $X$ and $Y$ in $\Omega_{\mathcal O}$. We use also the notation $\Omega_{\mathcal O}(Y,t)=\Omega_{\mathcal O}(\{Y\},t)$.
More generally,
when $I=\{(\Sigma_j,t_j)\}_{j=1}^J$ 
is a finite collection of pairs $(\Sigma_j,t_j)$, where $\Sigma_j\subset \Gamma_{\mathcal O}$ 
and $t_j>0$, we denote
\ba
\Omega_{\mathcal O}(I)=\bigcup_{j=1}^J \Omega_{\mathcal O}(\Sigma_j,t_j)
=\{X \in \Omega_{\mathcal O}\ ;\ d_{\mathcal O}(X,\Sigma_j)\leq t_j\hbox{ for some }j=1,\dots,J\},
\ea
For any measurable set $B\subset \Omega_{\mathcal O}$, 
we denote $ L^2(B)=\{v\in
 L^2(\Omega_{\mathcal O});\ v|_{\Omega_{\mathcal O}\setminus B}=0\}$,
identifying functions and their zero continuations.


\begin{lemma} \label{lem: step 4}
Assume that we are given the set  $\Gamma_{\mathcal O}$ as a differentiable manifold, 
 the metric on $\Gamma_{\mathcal O}$ and the BSP for $H_{\mathcal O}$ on $\Gamma_{\mathcal O}$. 
Then, for any given $f\in C^\infty_0(\Gamma_{\mathcal O}\times \R_+)$, $T>0$, and $I=\{(\Sigma_j,t_j)\}_{j=1}^J$,
where $\Sigma_j\subset \Gamma_{\mathcal O}$ are open sets or single points, and $t_j<T$,
we can determine
\beq\label{eq: a formula 2versionA}
a_{I,T}(f)=\int_{\Omega_{\mathcal O}\setminus \Omega_{\mathcal O}(I)}|u^{f}(T)|^2\,dV.
\eeq
\end{lemma} 

\noindent {\bf Proof.}
When  $\Sigma\subset \Gamma_{\mathcal O}$ is an open set and
  $h\in C^\infty_0(\Sigma\times \R_+)$, it follows from the finite velocity
of wave propagation  (see
e.g.\ \cite[Sec.\ 4.2]{Lad85}, see also \cite[Ch.\ 6]{IsKu08}) that the wave $u^h(t)=u^{h}(\cdotp,t)$ is supported in the
domain $\Omega_{\mathcal O}(\Sigma,t)$ at time $t>0$.
It follows from Tataru's seminal 
unique continuation result, see \cite{Ta1,Ta3}, that the set
\beq\label{eq: densitity}
\{u^{h}(t);\ h\in C^\infty_0(\Sigma\times \R_+)\}
\eeq
is dense in $L^2(\Omega_{\mathcal O}(\Sigma,t))$, see e.g.\ \cite[Theorem 3.10] {KKL01}.
This clearly implies that, when $T>0$ and 
$I=\{(\Sigma_j,t_j)\}_{j=1}^J$, where $\Sigma_j$ are open and  $t_j<T$,  the set
\ba
X_I^T&:=&\{u^{h}(T);\ h=h_1+\dots+h_J,\ h_j\in C^\infty_0(\Sigma_j\times [T-t_j,T])\}
\\ &=&\hbox{span}_{j=1,\dots,J}\, \{u^{h}(t_j);\ h\in C^\infty_0(\Sigma_j\times [0,t_j])\}
\ea
is dense in $L^2(\Omega_{\mathcal O}(I))$.

Next, we consider the non-linear functional
\ba
a_{I,T}(f)=\inf\{\|u^{f-h}(T)\|_{L^2(\Omega_{\mathcal O})}^2;\ 
 h=h_1+\dots+h_J,\ h_j\in C^\infty_0(\Sigma_j\times [T-t_j,T])\}
\ea
where $f\in C^\infty_0(\Gamma_{\mathcal O}\times \R_+)$, $T>0$, and
$I=\{(\Sigma_j,t_j)\}_{j=1}^J$, $\Sigma_j\subset \Gamma_{\mathcal O}$ are open, and $t_j<T$. 
By the formula (\ref{eq:Sect5Blagovecont}),
the BSP and the metric on $\Gamma_\O$ determine the value $a_{I,T}(f)$ for any $f$.
Moreover, as
$u^{f-h}(T)=
u^{f}(T)-u^h(T)$ and $X_I^T$ is dense in $L^2(\Omega_{\mathcal O}(I))$,
we see that
\beq\label{eq: a formula}
a_{I,T}(f)=\|(1-\chi_{\Omega_{\mathcal O}(I)})u^{f}(T)\|_{L^2(\Omega_{\mathcal O})}^2,
\eeq
where $\chi_{\Omega_{\mathcal O}(I)}(x)$ is the characteristic function of the set
$\Omega_{\mathcal O}(I)$ on $\Omega_{\mathcal O}$. 
{This proves the lemma for the case when all $\Sigma_j$ are open. 

 If for some $j$, the set $\Sigma_j$ is just a point $X_j \in \Gamma_{\mathcal O}$, 
we define for those $j$'s $\Sigma_j^{(k)}\subset \Gamma_{\mathcal O}$, $k = 1,2,\cdots$
to be open neighborhoods of $X_j$ such that $\overline{\Sigma_j^{(k+1)}} \subset \Sigma_j^{(k)}$ and $\bigcap_k\Sigma_j^{(k)} = \{X_j\}$.
For those $j$'s for which  $\Sigma_j$ is open, we define 
$\Sigma_j^{(k)}=\Sigma_j$.
Denote the corresponding finite collection of $(\Sigma_j^{(k)},t_j)$ by $I(k)$. Then 
$$
\Omega_{\mathcal O}(I(k+1)) \subset \Omega_{\mathcal O}(I(k)), \quad
\Omega_{\mathcal O}(I) = \bigcap_{k=1}^{\infty}\Omega_{\mathcal O}(I(k)),
$$
and for any $b \in L^2(\Omega_{\mathcal O})$,
$$
\big(1 - \chi_{\Omega_{\mathcal O}(I(k))}\big) b \to 
\big(1 - \chi_{\Omega_{\mathcal O}(I)}\big) b, \quad 
{\rm a.e.} \quad {\rm as} \quad k \to \infty.
$$
As $|\big(1 - \chi_{\Omega_{\mathcal O}(I(k))}\big)b(\cdot)| \leq 
|\big(1 - \chi_{\Omega_{\mathcal O}(I)}\big)b(\cdot)|$, a.e., using the monotone convergence theorem, we see that
$$
a_{I(k),T}(f) \to \|\big(1 - \chi_{\Omega_{\mathcal O}(I)}\big)u^f(T)\|^2_{L^2(\Omega_{\mathcal O})} = a_{I,T}(f).
$$
Thus, the BSP and the metric  on $\Gamma_{\mathcal O}$ determine $a_{I,T}(f)$ for such $I's$.}
\qed


\begin{definition}
Let  $I=\{(\Sigma_j,t_j)\}_{j=1}^J$, $I'=\{(\Sigma_j',t_j')\}_{j=1}^J$ 
and $T>0$, where $\Sigma_j,\Sigma'_j\subset \Gamma_{\mathcal O}$ and $t_j,t_j'<T$. 
We say that the relation $I\geq I'$ is valid on manifold 
$\Omega_{{\mathcal O}}$ if
\beq\label{eq: omega test}
\Omega_{\mathcal O}(I')\setminus \Omega_{\mathcal O}(I)\hbox{ has measure zero}.
\eeq
\end{definition}


\begin{lemma} \label{lem: step 5}
Let  $I=\{(\Sigma_j,t_j)\}_{j=1}^J$, $I'=\{(\Sigma_j',t_j')\}_{j=1}^J$ 
and $T>0$, where $\Sigma_j,\Sigma'_j\subset \Gamma_{\mathcal O}$ are open sets or single points and $t_j,t_j'<T$. 
Assume that we are given the set  $\Gamma_{\mathcal O}$ as a differentiable manifold, 
 the metric on $\Gamma_{\mathcal O}$, the BSP for $H_{\mathcal O}$ on $\Gamma_{\mathcal O}$, and the
collections $I$ and $I'$. 
Then we can determine whether the relation $I\geq I'$ is valid on manifold 
$\Omega_{{\mathcal O}}$  or not.
\end{lemma} 

\noindent {\bf Proof.} The relation $I\geq I'$ is valid on manifold 
$\Omega_{{\mathcal O}}$ 
if and only if
\beq\label{eq: a test}
a_{I,T}(f)\leq a_{I',T}(f)\quad\hbox{for all } f\in C^\infty_0(\Gamma_{\mathcal O}\times \R_+).
\eeq
Indeed, the equivalence of (\ref{eq: omega test}) and  (\ref{eq: a test})
follows from (\ref{eq: a formula}) and the fact that, by Tataru's 
density result (\ref{eq: densitity}), the functions $u^f(T)$,
 $f\in C^\infty_0(\Gamma_{\mathcal O}\times \R_+)$, are dense in $L^2(\Omega_{\mathcal O}(\Gamma_{\mathcal O},T))$.
As for given $f$, by Lemma \ref{lem: step 4}, we can evaluate both sides of (\ref{eq: a test}), 
using the BSP and the metric on $\Gamma_{\mathcal O}$, these data determine, for any pair $(I,I')$, if the relation  $I\geq I'$ is valid or not.
\qed

\medskip
For any $X_0 \in \Omega_{int}\setminus\partial\Omega_{int}$, introduce the exponential map 
$$
\exp_{X_0}:(\xi,t) \mapsto \gamma_{(X_0,\xi)}(t),
$$
 where $\xi \in S_{X_0}(\Omega_{int}) = \{\eta \in T_{X_0}(\Omega_{int});|\eta| = 1\}$ and $0 \leq t \leq s(X_0,\xi)$. Here $\gamma_{(X_0,\xi)}(t)$ is the geodesic on $\Omega$, parametrized by the arclength, with $\gamma_{(X_0,\xi)}(0) = X_0$, $\dot\gamma_{(X_0,\xi)}(0) = \xi$, and $[0,s(X_0,\xi))$ is the maximal interval of $t$, when $\gamma_{(X_0,\xi)}(t)$ stays in $\Omega_{int}$, that is,
 $s(X_0,\xi) = \sup\{t\, ; \, \gamma_{(X_0,\xi)}([0,t)) \subset \Omega_{int}\setminus\partial\Omega_{int}\}$.
Denote by
\begin{equation}
 s(X_0) = \inf_{\xi \in S_{X_0}(\Omega)}s(X_0,\xi)
\label{defines}
\end{equation}
so that
$$
B(X_0,s(X_0)) \subset \Omega_{int}\setminus\partial\Omega_{int}.
$$
Define now
$$
\tau(X_0,\xi) = \sup_{0<t<s(X_0)}\{t \,;\, d_{\emptyset}(\gamma_{(X_0,\xi)}(t),X_0) = t\}.
$$
At last, define
\begin{equation}
\tau(X_0) = \inf_{\xi\in S_{X_0}(\Omega_{int})} \tau(X_0,\xi).
\label{definetau}
\end{equation}
In geometric terms, the above definition of $\tau(X_0)$ means that in the ball $B(X_0,\tau(X_0)) \subset \Omega_{int}\setminus\partial\Omega_{int}$, it is possible to introduce the Riemannian normal coordinates 
$$
X\mapsto (\xi,t) : \xi \in S_{X_0}(\Omega_{int}), 0 \leq t < \tau(X_0)
$$
which satisfy $\gamma_{(X_0,\xi)}(t) = X$. 

We also need the boundary exponential map
$$
\exp_{\Gamma_{\mathcal O}}: 
\{(Z,t)\in 
\Gamma_{\mathcal O}\times\R_+\, ; \, 0 \leq t < s_\O(Z)\} \ni 
 (Z,t) \to \gamma_{(Z,\nu)}(t)\in  \Omega_{\mathcal O}. 
$$
Here $\nu$ is the interior unit normal (with respect to $\Omega_{\mathcal O}$) 
 to $\Gamma_{\mathcal O}$ and
\begin{equation}
s_{\mathcal O}(Z) = \sup\{t>0 \,;\, \gamma_{(Z,\nu)}((0,t)) \subset 
\Omega_{\mathcal O}\setminus \p \Omega_{\mathcal O}\}.
\label{defines1}
\end{equation}
For any $Z \in \Gamma_{\mathcal O}$, let
\begin{equation}
\tau_{\mathcal O}(Z) = \sup_{0\leq t \leq s_{\mathcal O}(Z)}\{t\,;\, d_{\mathcal O}(\gamma_{(Z,\nu)}(t),\Gamma_{\mathcal O}) = t\}.
\label{definetau1}
\end{equation}

In the following, we impose the following condition (C-2) on 
$\Sigma$.
\medskip

\noindent (C-2)         For $\mathcal O = \emptyset$ , $\Sigma$ 
is an open subset of $\Gamma$ such that
                  $d_{\emptyset}(\Sigma,\partial\Gamma) > 0$, and for $\mathcal O \neq 
\emptyset,$ $\Sigma =\partial\mathcal O$.
\medskip

 We define
\begin{equation}
\tau_{\mathcal O}(\Sigma) = \inf_{Z \in \Sigma}\tau_{\mathcal O}(Z).
\label{definetau2}
\end{equation}
In geometric terms, the above definition of $\tau_{\mathcal O}(\Sigma)$ means that, in the set 
$$
L(\Sigma,\tau_{\mathcal O}(\Sigma)) = \{\gamma_{(Z,\nu)}(t)\, ; \, Z \in \Sigma, 0 \leq t < \tau_{\mathcal O}(\Sigma)\} \subset \big(\Omega_{\mathcal O}\setminus\partial\Omega_{\mathcal O}\big) \cup \Sigma,
$$
it is possible to introduce the boundary normal coordinates
$$
X\mapsto  (Z,t), Z \in \Sigma, 0 \leq t < \tau_{\mathcal O}(\Sigma)
$$
satisfying $X = \gamma_{(Z,\nu)}(t)$.   Observe that when $\mathcal O = B(X,\rho)$, $X \in \Omega_{int}\setminus\Omega_{int}$ and $\rho > 0$ is small enough, then
$$
\tau_{\mathcal O}(\p {\mathcal O}) = \tau(X) - \rho.
$$


\begin{lemma} \label{lem: step 6A} 
 Assume that $\Sigma\subset \Gamma_{\mathcal O}$ satisfies
condition (C-2).
Let $Y \in \Sigma$, $Z\in\Gamma_{\mathcal O}$, $t< \tau_{\mathcal O}(\Sigma)$,  and  $X =\gamma_{(Y,\nu)}(t)$. 
Assume that we are given the set  $\Gamma_{\mathcal O}$ as a differentiable manifold, the metric on $\Gamma_{\mathcal O}$ and the BSP for $H_{\mathcal O}$ on $\Gamma_{\mathcal O}$.
 Then we can determine the distance $d_{\mathcal O}(X,Z)$
 {\newtekstt on $\Omega_{\mathcal O}$}. 
 \end{lemma} 

\noindent {\bf Proof.} 
Note that as $t< \tau_{\mathcal O}(\Sigma)$, the set $\Omega_{\mathcal O}(Y,t)\setminus \Omega_{\mathcal O}(\Gamma_{\mathcal O},t-\e)$ contains
a non-empty open set for all $\e>0$.  
For $s,\e>0$, let us  denote (see Fig.\ 4)
\ba
I_{\epsilon}(t) = \big\{(Y, t), (\Gamma_{\mathcal O},t-\e)\big\}, \quad
I'_\e(t, s) =\big\{(Z, s),\,(\Gamma_{\mathcal O},t-\e) \big\}.
\ea

\begin{figure}[htbp]
\begin{center}
\psfrag{1}{$Y$}
\psfrag{2}{$\Gamma$}
\psfrag{3}{$Z$}
\psfrag{4}{$X$}
\includegraphics[width=8cm]{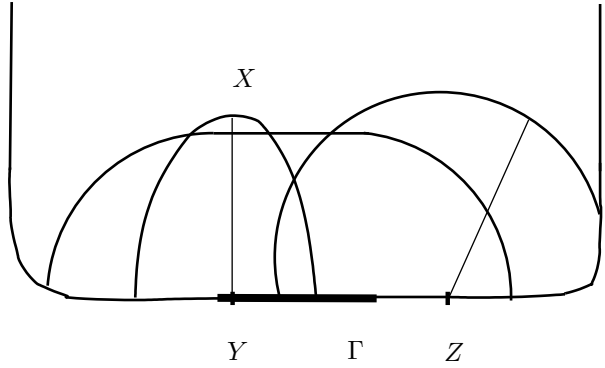} 
\label{pic 3}
\end{center}
\caption{In the figure ${\mathcal O}=\emptyset$ and $s$ is {\newtekstt  small enough so that
the set $\Omega_{\mathcal O}(Y,t)\setminus \Omega_{\mathcal O}(\Gamma_{\mathcal O},t-\e)$ is not contained in}
 $\Omega_{\mathcal O}(Z,s)$. This is the situation when 
$I'_\e(t, s)\not \geq I_{\epsilon}(t)
$.
}
\end{figure}

Let us next show 
that for any $r>0$ there is $\e_0>0$ such that 
$$
\Omega_{\mathcal O}(Y,t)\setminus \Omega_{\mathcal O}(\Gamma_{\mathcal O},t-\e)\subset B(X,r), \quad {\rm  when} \quad 
\e<\e_0.
$$
If this is not true, there are $r>0$, a sequence $\e_j\to 0$, and
$X_j\in \Omega_{\mathcal O}(Y,t)\setminus \Omega_{\mathcal O}(\Gamma_{\mathcal O},t-\e_j)$
such that $d_{\mathcal O}(X_j,X)\geq r$. As $\Omega_{\mathcal O}(Y,t)$ is compact, by considering a subsequence, we can 
assume that $X_j$ converge to $\tilde X\in \Omega_{\mathcal O}(Y,t)$. Then
\ba
& &d_{\mathcal O}(\tilde X,Y) =  \lim_{j\to \infty} d_{\mathcal O}(X_j,Y)\leq t,\\
& &d_{\mathcal O}(\tilde X,\Gamma_{\mathcal O}) = \lim_{j\to \infty} d_{\mathcal O}(X_j,\Gamma_{\mathcal O})\geq t,
\ea
implying that $Y$ is a closest point of $\Gamma_{\mathcal O}$ to $\tilde X$ and $d_{\mathcal O}(\tilde X,Y)=t$.
Let us recall  that the shortest curve from a point in $\Omega_{\mathcal O}$ to 
$\Gamma_{\mathcal O}$, 
which end point is an interior point of $\Gamma_{\mathcal O}$, is a normal  geodesic to $\Gamma_{\mathcal O}$. 
Thus, we see that $\tilde X =\gamma_{(Y,\nu)}(t)= X$, which is in contradiction
to  $d(\tilde X,X)\geq r$. Thus, the existence of $\e_0$ for any $r$ is proven.

The above implies that when $s>d(X,Z)$, the set 
 $\Omega_{\mathcal O}(Y,t)\setminus \Omega_{\mathcal O}(\Gamma_{\mathcal O},t-\e)$ is contained in $\Omega_{\mathcal O}(Z,s)$
for all sufficiently small $\e>0$ and therefore,
\beq\label{eq: Cond a1}
\hbox{there is $\e_1>0$ such that $I'_\e(t, s)\geq I_\e(t)$ for all $0<\e<\e_1$.}
\eeq
On the other hand, for $s<d(X,Z)$, 
the set 
 $\Omega_{\mathcal O}(Y,t)\setminus \Omega_{\mathcal O}(\Gamma_{\mathcal O},t-\e) \neq \emptyset$ 
 do not intersect with
 $\Omega_{\mathcal O}(Z,s)$ at any $\e>0$ small enough and thus (\ref{eq: Cond a1}) does not hold. 
Thus, by 
 Lemma \ref{lem: step 5},  we can find $d_{\mathcal O}(X,Z)$
 for any $ Z \in \Gamma_{\mathcal O}$
as the infimum of all $s>0$ for which (\ref{eq: Cond a1}) holds.
\qed

\medskip
 For  $\Sigma\subset \Gamma_{\mathcal O}$ satisfying (C-2) and 
$0<T<\tau_{\mathcal O}(\overline \Sigma)$, let $N_{\Sigma,T}$ and  $M_{\Sigma,T}$ be the sets
\begin{equation}
\begin{split}
&N_{\Sigma,T}=\{X \in \Omega_{\mathcal O};\ X=\gamma_{(Y,\nu)}(t),\
0\leq t \leq T,\ Y \in \overline \Sigma\},\\
&M_{\Sigma,T}=\{X \in \Omega_{\mathcal O};\ X =\gamma_{(Y,\nu)}(t),\
0<t <T,\ Y\in\Sigma\} \subset N_{\Sigma,T} = \overline{M_{\Sigma,T}}.
\label{defineMN}
\end{split}
\end{equation}
Note that $M_{\Sigma,T}$ is open in $\Omega_{\mathcal O}$.

\subsection{Boundary distance functions and reconstruction of topology}
Let us next consider the collection of the boundary
distance functions associated with $\Gamma_{\mathcal O}$. 
For each $X \in \Omega_{\mathcal O}$,
the corresponding {\it restricted boundary
distance function}, $r_X \in C( \Gamma_{\mathcal O})$  (note that 
${ \overline \Gamma_{\mathcal O}}$ is compact) is given by
\[
 r_X: \Gamma_{\mathcal O}\to \R_+,\quad r_X(Z)=d_{\mathcal O}(X,Z), \quad
Z \in  \Gamma_{\mathcal O}.
\]
The restricted boundary
distance functions define {\it the
 boundary distance map} $R_{\mathcal O}:\Omega_{\mathcal O}\to C( \Gamma_{\mathcal O})$,
$R_{\mathcal O}(X)=r_X$. 
The boundary distance representation of $N_{\Sigma,T}\subset \Omega_{\mathcal O}$ is the set
\[
R_{\mathcal O}(N_{\Sigma,T})=\{r_X \in C( \Gamma_{\mathcal O}); \, X \in N_{\Sigma,T}\},
\]
that is, the image of $N_{\Sigma,T}$ in $R_{\mathcal O}$. Clearly $R_{\Gamma_{\mathcal O}}:\Omega_{\mathcal O}\to C( \Gamma_{\mathcal O})$ is continuous.


\begin{lemma} \label{lem: step 7}
Assume that we are given the set  $\Gamma_{\mathcal O}$ as a differentiable manifold, 
 the metric on $\Gamma_{\mathcal O}$, the BSP for $H_{\mathcal O}$ on $\Gamma_{\mathcal O}$, 
 an open set $\Sigma\subset \Gamma_{\mathcal O}$ satisfying
condition (C-2),
  and 
 $0<T<\tau_{\mathcal O}(\overline \Sigma)$.
Then we can determine the set
 \ba R_{\mathcal O}(N_{\Sigma,T})=R_{\mathcal O}(\{\gamma_{(Y,\nu)}(t);\ Y \in \overline\Sigma ,\ 0\leq t\leq T\}).
 \ea
\end{lemma} 

\noindent {\bf Proof.} 
By Lemma \ref{lem: step 6A},  for $Y \in  \Sigma$, $t<T$
and $Z \in \Gamma_{\mathcal O}$, we can find $d_{\mathcal O}(X,Z)$
where $X =\gamma_{(Y,\nu)}(t)$ from BSP. This gives us the function
$r_X(Z),$ $Z \in \Gamma_{\mathcal O}$, and for such $X$'s. Thus, BSP 
and the metric on $\Gamma_\O$ determine the set $R_{\mathcal O}(M_{\Sigma,T})$.
{\newtekstt Using  (\ref{defineMN}), we obtain $R_{\mathcal O}(N_{\Sigma,T})$ by closure of 
$R_{\mathcal O}(M_{\Sigma,T})$ in $C(\Gamma_{\mathcal O})$.}
\qed
\medskip

Consider properties of $R_{\mathcal O}$.
Assume that $r_X=r_Y$ for some
$X, Y\in N_{\Sigma,T}$. Let $Z \in  \Gamma_{\mathcal O}$ be the point where the function $r_X$ attains its minimum. Then, it is the closest point of
$ \Gamma_{\mathcal O}$ to $X$.
Thus, the shortest geodesic from $X$ to $Z$ is normal 
to $\Gamma_{\mathcal O}$, i.e. $X =\gamma_{(Z,\nu)}(t)$ with $t=r_X(Z)$. The same arguments show that $Z$ is also the closest
point of $\Gamma_{\mathcal O}$ to $Y$ and $t=r_Y(Z)$, and hence
$Y =\gamma_{(Z,\nu)}(t)$. Thus $X = Y$ and $R_{\mathcal O}$ is injective on $N_{\Sigma,T}$.

Thus,  map 
$R_{\mathcal O}:N_{\Sigma,T}\to R_{\mathcal O}(N_{\Sigma,T})$ is
a bijective continuous map defined on a
compact set, implying that it is a homeomorphism.
This implies that
the map 
$R_{\mathcal O}:M_{\Sigma,T}\to R_{\mathcal O}(M_{\Sigma,T})$ is
a homeomorphism. As BSP and the metric on $\Gamma_{\mathcal O}$ determine the manifold $R_{\mathcal O}(M_{\Sigma,T})$ with
its topological structure inherited from $C(\Gamma_{\mathcal O})$, we see that
these data determine
the manifold $M_{\Sigma,T}$  as a topological space.


\begin{lemma} \label{lem: step 8}
The set
 $R_{\mathcal O}(M_{\Sigma,T})\subset C(\Gamma_{\mathcal O})$
can be endowed, in a constructive way, with a differentiable structure and a metric tensor $\tilde G$,
so that $(R_{\mathcal O}(M_{\Sigma,T}),\tilde G)$ becomes a manifold which is
 isometric
to $(M_{\Sigma,T},G)$ with $R_{\mathcal O}$ being an isometry. 
\end{lemma} 

For compact manifolds, the result analogous to Lemma  \ref{lem: step 8}  is presented in detail in \cite[Sect. 3.8]{KKL01}. Since the proof
is based on local constructions, it works for non-compact manifolds without any change. However, for the convenience of the reader, we present this construction.

{\bf Proof}. Let us define
the evaluation functions,
$E_{Z}, \, Z \in \Gamma_{\mathcal O}$, 
\ba
E_{Z}:R_{\mathcal O}(M_{\Sigma,T}) \to \R,
\quad
E_{Z}(r_X)=r_X(Z)
=d_{\mathcal O}(X,Z).
\ea
For $r(\cdot) \in R_{\mathcal O}(M_{\Sigma,T})$ corresponding to a point
$X \in M_{\Sigma,T}$,
i.e. $r(\cdot)= r_X(\cdot)$,  we can choose points $Z_1, \dots, Z_n \in \Gamma_{\mathcal O}$ close
to the nearest point of $\Gamma_{\mathcal O}$ to $X$ so that $X \mapsto (d_{\mathcal O}(X,Z_j))_{j=1}^n$ forms a system of coordinates on $\Omega_{\mathcal O}$
 near $X$, see \cite[Lem.\ 2.14]{KKL01}. Similarly, the functions
$E_{Z_j}, \, j=1,\dots,n,$ form a system of coordinates in $R_{\mathcal O}(M_{\Sigma,T})$ near
$r_{X}$. These coordinates provide for $R_{\mathcal O}(M_{\Sigma,T})$ a differential structure which
makes it diffeomorphic to manifold $M_{\Sigma,T}$.

Let us denote by $\tilde{G}$ the metric on $R_{\mathcal O}(M_{\Sigma,T})$ which makes it
isometric to
$(M_{\Sigma,T},G)$, that is, $\tilde G =((R_{\mathcal O})^{-1})^{\ast}G$.
Let $r\in R_{\mathcal O}(M_{\Sigma,T})$ and $X\in M_{\Sigma,T}$ be such 
that $r=r_X$. Let $Z_0$ is a point where $r$ obtains its minimum, that is,
the closest point of $\Gamma_{\mathcal O}$ to $X$. When $Z$ is close
to $Z_0$, the
differentials of functions $E_Z$ 
are covectors of length $1$ on $(R_{\mathcal O}(M_{\Sigma,T}),\tilde G)$,
see \cite[Lem.\ 2.15]{KKL01}. This is equivalent
to the fact that the gradients of the distance functions $X \mapsto d_{\mathcal O}(X,Z)$
have length one.
By this observation, it is possible to find infinitely many
covectors $dE_Z$, $Z \in \Gamma_{\mathcal O}$ of length 1 at any point
$r$ of $R_{\mathcal O}(M_{\Sigma,T})$. Using such vectors,
one can reconstruct the metric tensor $\tilde{G}$ at $r$. By the above considerations, 
BSP determines 
the manifold  $(M_{\Sigma,T},G)$ up to an isometry. 
\qed

\subsection{Continuation of the data}
Let us now consider the case when ${\mathcal O}=\emptyset$ and we are given
the set $\Gamma$ as a differentiable manifold,
the metric $G$ on $\Gamma$, and the BSP for $H_{\emptyset}$.
Assume that there are two manifolds $\Omega_{int}^{(1)}$ and $\Omega_{int}^{(2)}$ 
such that $\Gamma$ is isometric to subsets $\Gamma^{(j)}\subset \p \Omega_{int}^{(j)}$ for $j=1,2$
and that the BSP for $H_{\emptyset}^{(j)}$, $j=1,2$, coincides with the given data.
Let now $\Sigma\subset \Gamma$ satisfy condition (C-2) and 
\ba
0<T<\min(\tau^{(1)}_{\emptyset}(\overline \Sigma),\tau^{(2)}_{\emptyset}(\overline \Sigma)).
\ea 
Then the above constructions show that the manifolds
\ba
M^{(j)}_{\Sigma,T}=\{X \in \Omega_{int}^{(j)};\  X = \gamma_{(Y,\nu)}(t),\
0<t <T,\ Y\in  \Sigma\}
\ea
with $j=1$ and $j=2$, are isometric. Thus, we can consider the set $M^{(1)}_{\Sigma,T}$,
denoted by $U_1$ as a subset of both manifolds $\Omega_{int}^{(1)}$ and $\Omega_{int}^{(2)}$,
and, by the previous considerations, we can construct a metric $\tilde G$
on it which makes $(U_1, \tilde G)$ isometric to $(M^{(j)}_{\Sigma,T}, G^{(j)})$,
$j=1,2$.

We continue the construction by continuation of the data using
Green's functions, cf.\ \cite{LTU03,LU01}. 
To this end, let
$z \in \C\setminus\R_+$ and consider the Schwartz kernel 
$G_{\mathcal O}(z; Y,Y')$ of the operator $(H_{\mathcal O} -z)^{-1}$. It satisfies the equation
\beq\label{eq: elliptic}
& &(H_{\mathcal O}-z)G_{\mathcal O}(z; \cdotp,Y')=\delta_{Y'},\quad Y,Y'\in \Omega_\O=\Omega_{int}\setminus \O,\\
& &\nonumber
\p_\nu G_{\mathcal O}(z; \cdotp,Y')|_{\p \Omega_\O}=0.
\eeq
We denote $ G(z; Y,Y')= G_{\O}(z: Y,Y')$ when $\O=\emptyset$.


\begin{lemma} \label{lem: step 1}
Let $U\subset \Omega_{int}$ be a connected neighborhood of an open 
set  $\Sigma\subset \Gamma$, where $\Sigma$  satisfies condition (C-2) with $\O=\emptyset$. Let
$X_0 \in U\setminus\partial\Omega_{int}$ and $\rho>0$
be such that $\mathcal O = B(X_0,\rho)\subset U\setminus\partial\Omega_{int}$.
Assume that  we are given the metric tensor $G$ in $U$.
Then  BSP on $\Gamma $ for the operator $H_{\emptyset}$ determines 
$G(z; Y,Y')$ for $Y, Y'\in U$ and $z\in \C\setminus {\mathcal E}(H_{\emptyset})$. 
 Moreover, these data determine
 BSP on $\Gamma_{\mathcal O}$ for the operator $H_{\mathcal O}$.
\end{lemma} 

\noindent {\bf Proof.} By Lemma \ref{lemma 5.6},
 BSP on $\Gamma $ determines the N-D map $\Lambda(z)$ at $\Gamma\times \Gamma$.
By Lemma \ref{lemma 5.4}, the Schwartz kernel of the N-D map $\Lambda(z)$ at $\Gamma\times \Gamma$ coincides with $G(z; Y,Y')$. Thus we know the function  $G(z; Y,Y')$ for $Y, Y'\in \Sigma$. As the Neumann boundary values of $Y \mapsto G(z; Y,Y')$ on $\Gamma\setminus \{Y'\}$ vanish,
using the Unique Continuation Principle for the elliptic equation (\ref{eq: elliptic}) in the $Y$ variable, 
we see that the values of $G(z; Y,Y')$ are uniquely determined for $Y'\in \Sigma$ and  
$Y \in U\setminus\{Y'\}$.
Using the symmetry $\overline{G(z; Y,Y')}=G(\overline z; Y',Y)$ and again the Unique Continuation Principle, now in the $Y' $ variable,  we can  determine the values of
$G(z; Y,Y')$ in  $\{(Y,Y')\in U\times U;\ Y \not=Y'\}$. Considering
$G(z; Y,Y')$ as a locally integrable function, we see that it is
defined a.e.\ in  $U\times U.$

For $Y' \in \left(\Omega_{\mathcal O}\cap U\right)\setminus 
\p\Omega_{\mathcal O}$, denote by $G^{ext}_{\mathcal O}(z; Y,Y')$ a smooth extension of 
$G_{\mathcal O}(z; Y,Y')$ into $\mathcal O$. Then 
$$
(- \Delta_G - z)G^{ext}_{\mathcal O}(z; Y,Y') - \delta(Y,Y') = 
F(Y,Y') \in C^{\infty}(\Omega_{int}),
$$
where ${\rm supp}\, F(\cdot,Y') \subset {\overline{\mathcal O}}$. 
Therefore,
$$
G_{\mathcal O}(z; Y,Y') = G(z; Y,Y') + 
\int_{\mathcal O}G(z; Y,Y'')F(Y'',Y')dV_{Y''}.
$$
In particular, 
\begin{equation}
\partial_{\nu(Y)}G(z; Y,Y') + \int_{\O}\partial_{\nu(Y)}G(z; Y,Y'')F(Y'',Y')dV_{Y''} = 0, \quad Y \in \partial\mathcal O,
\label{Gznormal}
\end{equation}
where $\nu(Y)$ is the unit normal to $\mathcal O$ at $Y$. On the other hand, if $F(\cdotp, Y')
 \in C^{\infty}(U)$, ${\rm supp}\, F(\cdotp, Y') \subset \overline{\mathcal O}$, 
 satisfies (\ref{Gznormal}), the function
\begin{equation}
G(z; Y,Y') + \int_{\O}G(z; Y,Y'')F(Y'',Y')dV_{Y''}, \quad 
Y,Y' \in U\setminus\overline  \O,
\label{extra}
\end{equation}
is $G_{\mathcal O}(z; Y,Y')$. As we have in our disposal $G(z; Y,Y')$ for $Y, Y' \in U$, 
we can verify for a given $F$, condition (\ref{Gznormal}).
 
Now, we return to $\Omega_{int}^{(1)},\, \Omega_{int}^{(2)}$ with $\Gamma$ 
and BSP on $\Gamma$ being the same.
We denote the associated functions appearing above by adding the superscript $(j)$. Let
(\ref{Gznormal}) holds with
$G(z; Y,Y')$, $F(Y'',Y')$ replaced by $G^{(1)}(z; Y,Y')$, $F^{(1)}(Y'',Y')$, respectively. Since
$G^{(1)}(z; Y,Y') = G^{(2)}(z; Y,Y')$ on $U\times U$, (\ref{Gznormal})
 also holds with $G(z; Y,Y')$,
$F(Y'',Y')$ replaced by
$G^{(2)}(z; Y,Y')$, $F^{(1)}(Y'',Y')$, respectively. Thus,  for  $Y, Y' \in
U\setminus\overline{ \mathcal O}$, we have
\HOX{Proof is slightly changed since somehow I feel that reference to L. 5.4 is a bit  difficult to follow}
$$
G_{\mathcal O}^{(j)}(z; Y,Y') = G^{(j)}(z; Y,Y') + \int_{\mathcal
O}G^{(j)}(z; Y,Y'')F^{(1)}(Y'',Y')dV_{Y''}, \quad j=1, 2,
$$  
{\newtekstt so that
$$
G_{\mathcal O}^{(1)}(z; Y,Y') =G_{\mathcal O}^{(2)}(z; Y,Y'),
\quad z \in \C \setminus \R,\,\, Y, Y' \in
U\setminus\overline{ \mathcal O}.
$$
In particular, this implies that 
$\Lambda_{\mathcal O}^{(1)}(z) = \Lambda_{\mathcal O}^{(2)}(z)$, $z \in {\C}\setminus{\R}$. }
Then by Lemma \ref{lemma 5.6}, BSP's for $H_{\mathcal O}^{(1)}$ and $H_{\mathcal O}^{(2)}$ coincide.
 \qed

\medskip
Next we show that we can use these data to determine the critical distance which we use
in the step-by-step construction of the manifold. 


\begin{lemma} \label{lem: step 6} Let  $X_0\in \Omega_{int}\setminus\partial\Omega_{int}$ and $0<\rho<\tau(X_0)/2$.
Let  ${\mathcal O}=B(X_0,\rho)$ and $\Gamma_\O=\p {\mathcal O}$. 
Assume that we are given the set  $\Gamma_{\mathcal O}$ as a differentiable manifold, 
the metric $G\big|_{\Gamma_{\mathcal O}}$ on $\Gamma_{\mathcal O}$, and 
the BSP for $H_{\mathcal O}$ on $\Gamma_{\mathcal O}$. Then 
these data determine $\tau_{\mathcal O}(\Gamma_\O)=\tau(X_0)-\rho$.
  \end{lemma} 

\noindent {\bf Proof.}  
 Let us assume that
 $t_0 < \tau(X_0) - \rho$. Then,
for any $Y \in \Gamma_\O$, the set 
$\Omega_{\mathcal O}(Y,t_0)\setminus\Omega_{\mathcal O}(\Gamma_\O,t_0-\e)$ 
contains an open neighborhood of $\gamma_{(Y,\nu)}(t_0-\e/2)$ and, therefore,
has  positive measure. 
Hence, if  $t < \tau(X_0) - \rho$, then  the condition
\begin{equation}
\forall Y \in \Gamma_\O\, \,\forall\e>0: \ 
I_{\epsilon, t} := \left\{(\Gamma_\O,t-\e)\right\} \not \geq 
I'_{Y,t} := \left\{(Y,t)\right\}
\label{A1}
\end{equation}
is valid.

Let us next assume that condition (\ref{A1}) is valid
and consider its consequences.

First, observe that by (\ref{defines}) and (\ref{definetau}), we have either 
\smallskip

\noindent
(a)  $s(X_0)=\tau(X_0)$ and there is $Y \in \Gamma_\O$ such that  
$X = \gamma_{(Y,\nu)}(\tau(X_0)-\rho) \in \partial\Omega_{int}$, 
\smallskip

\noindent
or
\smallskip

\noindent
(b)  $s(X_0) > \tau(X_0)$ and there are
$Y \in  \Gamma_\O$ and $s$ such that $s(X_0) > s > \tau(X_0)-\rho$ and
$d_{\mathcal O}(\gamma_{(Y,\nu)}(s), \Gamma_\O) < s$.
\smallskip

Let us consider these two cases separately.

(a) It follows from (\ref{defines})  and (\ref{definetau}) that $X$ is a closest point to
$X_0$  on $\partial\Omega_{int}$. 
Therefore, the geodesic $\gamma_{(Y,\nu)}$ intersects 
$\partial\Omega_{int}$ normally at 
$X=\gamma_{(Y,\nu)}(s),$  $s = \tau(X_0)-\rho$.

Assume next that $t>0$ is such that
\beq\label{Matti A 1}
\forall \e>0: I_{t,\epsilon}  \not \geq I'_{Y,t}.
\eeq
Then for any $\e>0$ there is 
\begin{equation}\label{AaA}
X_\e \in \Omega_{\mathcal O}(Y,t)\setminus\Omega_{\mathcal O}( \Gamma_\O,t-\epsilon).
\end{equation}
As $\Omega_{\mathcal O}(Y,t)$ is relatively compact, there are $\e_n\to 0$ and $X_n=X_{\e_n}$
such that $X_n \to X'\in \Omega_{int}$ as $n\to \infty$. Then
\begin{equation}
d_{\mathcal O}(X',Y) = t, \quad d_{\mathcal O}(X', \Gamma_\O) = t.
\label{A6}
\end{equation}
This shows that $Y$ is the closest point of $ \Gamma_\O$ to $X'$ in
$\Omega_\O$.
{\newtekstt Consider a  
 shortest curve $\mu(s)$ from $Y$ to $X'$.  By \cite{Al},  a shortest curve  between
 on
a manifold with boundary is a  $C^1$-curve. Moreover, it is a geodesic on
$\Omega_\O \setminus \p\Omega_\O$. Since $\mu(s)$ is  a shortest curve
from $X'$ to $ \p\Omega_\O$, it is normal to $\p\Omega_\O$ at $Y$.
Thus $\mu(s)=\gamma_{Y, \nu}(s),\, s \leq \tau(X_0)-\rho$. 
However,  $\gamma_{(Y,\nu)}(s)$ hits 
 $\p \Omega_{int}$ normally at $s=\tau(X_0)-\rho$.
 Therefore, by the short-cut arguments, } 
  we see that the curve
$\gamma_{(Y,\nu)}([0,\tau(X_0)-\rho])\subset \Omega_\O$ can not be extended to a longer
curve which is a shortest curve between $Y$ and its other 
end point. Thus $\mu \subset \gamma_{(Y,\nu)}([0,\tau(X_0)-\rho])$,
implying that $t=d_\O(Y,X')\leq \tau(X_0)-\rho$.
Hence in the case (a) the condition (\ref{A1}) implies that 
$t\leq \tau(X_0)-\rho$.

(b) In this case arguments are similar but slightly simpler. 
Again, assume that $t>0$ is such that (\ref{Matti A 1})
is satisfied. Again, there are  $\epsilon_n>0$ and  $X_n=X_{\e_n}$
satisfying (\ref{AaA}),
such that  $X_n\to X'$ and $X'\in \Omega_{int}$ satisfies (\ref{A6}).
Moreover, a shortest curve $\mu(s)$ from $Y$ to $X'$ coincides with 
the normal geodesic $\gamma_{(Y,\nu)}(s)$ for small values of
$s$. Since the geodesic 
$\gamma_{(Y,\nu)}([0,s'])$ is a shortest curve between its   
end points for $s'\leq \tau(X_0)-\rho$ but not for
$s(X_0)-\rho> s'>\tau(X_0)-\rho$, we see that $\mu \subset \gamma_{(Y,\nu)}([0,\tau(X_0)-\rho])$ and thus $t\leq \tau(X_0)-\rho$.

Therefore, in both cases (a) and (b), 
the condition (\ref{A1}) implies that 
$t\leq \tau(X_0)-\rho$. 
Combining  these facts, we see that
\ba
\tau(X_0) - \rho = \sup\{t>0 \, ; \ \hbox{condition (\ref{A1})  is  satisfied for $t$}\}.
\ea 
The lemma then follows from this and Lemma \ref{lem: step 5}.
\qed


\subsection{Proof of Theorem \ref{BC thm}}
We are now in a position to complete the proof of Theorem \ref{BC thm}.

\subsubsection{Local reconstruction of Riemannian structure}
We start our considerations with  ${\mathcal O}=\emptyset$. 
Let $\Sigma \subset \Gamma$ satisfies condition (C-2)
and $T>0$ be a sufficiently small. In fact, we can consider any 
$0<T<\tau_\emptyset(\overline \Sigma)$. Using 
Lemma \ref{lem: step 7} 
we see that the set 
$R_{\emptyset}(M_{\Sigma,T}) \subset C(\Gamma)$ is
uniquely determined. On this set we 
introduce the boundary normal coordinates,
$$
r(\cdot) \mapsto (Z,t), \quad t = \min_{Z' \in \Sigma}r(Z'),
$$
where $Z$ is the unique point on $\Sigma$ on which $r(\cdot)$ attains its minimum. 
Observe that these coordinates on $R_{\emptyset}(M_{\Sigma,T})$ coincide with the boundary normal
 coordinates of the point $X\in \Omega_{int}$ such that
$$
r(\cdot) = r_X(\cdot).
$$
Thus, $R_{\emptyset}(M_{\Sigma,T})$ with the above coordinates is diffeomorphic to $M_{\Sigma,T}$. 

Next we use Lemma \ref{lem: step 8} to endow $R_{\emptyset}(M_{\Sigma,T})$ with Riemannian metric, $\tilde G$, 
so that $(R_{\emptyset}(M_{\Sigma,T}),\tilde G)$ is isometric to manifolds
$(M_{\Sigma,T},G)$.

\bigskip
\noindent
{\it Remark}. For the inverse scattering problem considered
 in the introduction, Section 6.5.1 is not necessary, 
because we know {\it a priori} the Riemannian structure of the open set 
$\left(\Omega_{int}\setminus\partial\Omega_{int}\right)\cap \Omega_1$. 
However, to make the results of \S 6 appropriate for general 
non-compact manifolds with asymptotically cylindrical ends, we have included this step.

\subsubsection{Iteration of local reconstruction}
To describe the procedure which we will iterate, 
let us  assume that $U_1\subset \Omega_{int}$ is
a connected neighborhood $\Sigma \subset \Gamma$ which
satisfies condition (C-2) with $\O=\emptyset$
and that we know the
Riemannian manifold $(U_1,G)$ up to an isometry. Since
the set $(R_{\emptyset}(M_{\Sigma,T},\tilde G))$ is already determined,
we can take $U_1=M_{\Sigma,T}$, where $T>0$ is sufficiently small . 

Choose $X_1 \in U_1$ and $\rho > 0$ such that $\mathcal O = B(X_1,\rho) \subset U_1$. 
By Lemma \ref{lem: step 1} we can determine $G_z(Y,Y')$ for all
$Y, Y' \in U_1$ and $z \in \C\setminus\R$.
Then Lemma \ref{lem: step 1} gives us BSP on $\partial{\mathcal O}$. Therefore by Lemma \ref{lem: step 6}, these
data determine $\tau_{\mathcal O}(\Gamma_{\mathcal O})$, hence $\tau(X_1) =
\tau_{\mathcal O}(\Gamma_{\mathcal O}) + \rho$. Take any $X \in
B(X_1,\tau)\setminus{\mathcal O}$, where $\tau = \tau(X_1)$, and let $Y$ be the
intersection of $\partial{\mathcal O}$ and the geodesic with end points $X_1$ and $X$.
Taking any $Z \in \partial{\mathcal O}$ and applying Lemma \ref{lem: step 6A}, we can then find $d_{\mathcal O}(X,Z)$.

Using, similarly to the above, Lemmas \ref{lem: step 7} and \ref{lem: step 8}, we 
can find the image of the embedding
$R_{\mathcal O}:B(X_1,\tau)\setminus\mathcal O\to C(\p{\mathcal O})$. We then recover, 
in the boundary normal coordinates associated with $\p{\mathcal O}$, i.e. the Riemannian normal coordinates 
centered at $X_1$, the metric tensor $G$ on $B(X_1,\tau)\setminus B(X_1,\rho)$,
{\newtekstt and, since $G$ on $B(X, \rho)$ is known, on the whole $B(X, \tau)$.}
This construction makes it possible to introduce the structure of the differentiable manifold 
on $U_1\bigsqcup B(X_1,\tau)$ which we considered, by now, as a disjoint union of two Riemannian manifolds. 
Next we glue these two components together.
To this end we observe that, since $\mathcal O \subset U_1$, we have in our disposal Green's function $G(z; Y,Y')$ for 
$ Y,Y' \in \mathcal O$ and $z \in \C\setminus\R$. The set $\O$ can be considered also as
the subset $B(X_1,\rho)$ of $B(X_1,\tau)$, and thus we know
the function $G(z; Y,Y')$ for 
$ Y,Y' \in B(X_1,\rho)$ e.g. in the Riemannian normal coordinates 
centered at $X_1$. Thus, using the Unique Continuation Principle, 
we can determine, in the Riemannian normal coordinates, the function $G(z; Y,Y')$ 
for all $Y \in B(X_1,\tau)$ and $Y' \in B(X_1,\rho)$.

 Since $Y'\mapsto G(z; Y,Y')$ is a smooth function in $\Omega_{int}\setminus \{Y\}$
and  $G(z; Y,Y') \to \infty$ as $Y'\to Y$, we see that for $Y_1,Y_2\in  \Omega_{int}$,
we have $Y_1=Y_2$ if and only if $G_z(Y_1,Y') = G_z(Y_2,Y')$ for all 
$Y'\in  \Omega_{int}, \, z \in \C \setminus \R$. Using the  
 Unique Continuation Principle,  this is equivalent to 
 $G(z; Y_1,Y') = G(z; Y_2,Y')$ for all $Y'\in  B(X_1,\rho), \, z \in \C \setminus \R$. 
 Next, let us define that  
the points $X_U \in U_1$ and $X_B \in B(X_1,\tau)$ 
are equivalent and denote $X_U\sim X_B$ if
$G(z; X_U,Y') = G(z; X_B,Y')$ for all $Y' \in  B(X_1,\rho), \, z \in \C \setminus \R$. 
Then the manifold $U_2=U_1\cup  B(X_1,\tau)\subset \Omega_{int}$
is diffeomorphic to manifold $(U_1\bigsqcup B(X_1,\tau))/\sim$,
which is  obtained by 
 glueing together the equivalent points on $U_1$ and $B(X_1,\tau)$. As we know the 
metric tensor on both $U_1$ and $B(X_1,\rho)$, we have reconstructed 
a Riemannian manifold $(U_2,G)\subset (\Omega_{int},G)$ up to an isometry.

\subsubsection{Maximal reconstruction}
Let us iterate the above process, that is,
we start from an open set $\Sigma\subset  \Gamma$ satisfying
condition (C-2) with $\O=\emptyset$, construct its neighborhood $U_1$,
 and iterate the construction by choosing at each step $j=1,2,\dots$ a point 
 $X_j\in U_j$ and constructing a Riemannian manifold isometric to
 $U_{j+1}=U_j\cup B(X_j,\tau(X_j)\subset \Omega_{int}$. 

{\newtekstt Consider the open sets in $\Omega_{int} \setminus \p\Omega_{int}$
which can be reconstructed, with the metric, when we are given the set $\Gamma$ 
with its metric 
and the BSP on $\Gamma$.
As the collection of these sets is closed with respect to taking the union, consider
maximal open set 
$U_{max}\subset \Omega_{int} \setminus \p\Omega_{int}$ which can be reconstructed, with its metric, from the set $\Gamma$ with its metric and the BSP on $\Gamma$.}
Let us  show that $U_{max} = \Omega_{int}\setminus\partial\Omega_{int}$. 
Since $\Omega_{int}\setminus\partial\Omega_{int}$ is connected, it suffices to show that $U_{max}$ 
is open and closed. By construction, $U_{max}$ is open. 
Let now $X \not\in \partial\Omega_{int}$ be a limit point of $U_{max}$, i.e., 
$X = \lim_{n\to \infty} X_n$, $X_n \in U_{max}$. Denote $a = d(X,\partial \Omega_{int})$ 
so that if $Y \in B(X, a/4)$, then $s(Y) \geq 3a/4$, see (\ref{defines}). Since the cut locus distance 
of the Riemannian normal coordinates is continuous with respect to the center, 
see e.g. \cite[Sec. 2.1]{KKL04} or \cite{GrKlMe}, there is $\delta > 0$ such that 
$\tau(Y) \geq \delta$ for all $Y \in B(X, a/4)$.

Let now $X_n \in U_{max}$ satisfy the inequality $d(X_n,X) < \sigma = \min(a/4,\delta/4)$. Let us assume that $X_n$ has a neighborhood 
  $B(X_n,\rho_n)$, with a sufficiently small $\rho_n < d(X_n,X)$, which
 can be reconstructed using $N(n)$ iteration steps, that is,
  $B(X_n,\rho_n)\subset U_{N(n)}$.
 Then $\tau(X_n) > 4\sigma$ so that $X \in B(X_n,\tau(X_n))$. By Lemma \ref{lem: step 1},
 we can find the BSP for the operator $H_{\mathcal O}$ with $\mathcal O = B(X_n,\tau(X_n))$ 
and, using one more iteration step, to reconstruct the Riemannian structure on 
$U_{N(n)}\cup B(X_n,\tau(X_n))$ which includes the point $X$. 
Therefore, the point $X$ is in $U_{max}$.
This shows that  $U_{max}$ is relatively open and closed in 
$\Omega_{int}\setminus\partial\Omega_{int}$. 
Thus, $U_{max} = \Omega_{int}\setminus\partial\Omega_{in}$.

The above shows that using an enumerable number of iteration 
steps we can construct a Riemannian manifold isometric
to  $(\Omega_{int}\setminus\partial\Omega_{int},G)$. Thus we have reconstructed the Riemannian manifold $(\Omega_{int}\setminus\partial\Omega_{int},G)$ up to an isometry.

It remains to identify the differentiable and Riemannian structures near $\partial\Omega_{int}$.
 Observe that $\Omega_{int}$ is just the closure of $\Omega_{int}\setminus\partial\Omega_{int}$ 
with respect to the distance function generated by the metric $G$ on
 $\Omega_{int}\setminus\partial\Omega_{int}$. Moreover, for any open
 relatively compact set
$\Sigma \subset \partial\Omega_{int}$, $\tau_\emptyset(\overline\Sigma) \geq \delta > 0$.

Let $0<t< \delta$ and consider the set
$$
\Sigma_{t} = \{X \in \Omega_{int}\setminus\partial\Omega_{int} \, ;\, 
d(X,\partial\Omega_{int}) = t, \ d(X,Z)=t, \ {\rm for} \ {\rm some} \ Z \in \Sigma\}.
$$
This implies that for $X\in \Sigma_t$ the closest point $Z\in \Omega_{int}$
is in $\Sigma$ and $X=\gamma_{Z,\nu}(t)$. Therefore, $\Sigma_t$
is a smooth $(n-1)$-dimensional open submanifold in $\Omega_{int}$ of points having the form $X = \gamma_{(Z,\nu)}(t)$, 
$Z \in \Sigma$.
This makes it possible to introduce the boundary normal coordinates in $M_{\Sigma,\delta}$ which provides the differentiable structure near $\Sigma$.  Writing the metric tensor $G$ in these coordinates and extending
this tensor continuously on $\Sigma$, we find the metric tensor in $\Omega_{int}$ in
the boundary normal coordinates associated to 
$\Sigma$. \qed

\subsection{Proof of Theorem \ref{main thm}}
Having Theorem \ref{BC thm} in our disposal, Theorem \ref{main thm} follows immediately from Lemma \ref{lemma 5.7}. \qed


\end{document}